\documentclass[journal, onecolumn]{IEEEtran}

\ifCLASSINFOpdf
\else
\fi
\usepackage[cmex10]{amsmath}
\usepackage{amssymb}
\usepackage{algorithmic}

\usepackage{xcolor}
\usepackage{mathtools}
\usepackage{amsfonts}
\usepackage{tensor}
\usepackage{hyperref}

\DeclareMathOperator{\tr}{tr}
\DeclareMathOperator{\rank}{rank}

\DeclareMathOperator{\diag}{diag}

\newcommand{\mprob}{\ensuremath{\mathbb{P}}}

\newtheorem{theorem}{Theorem}[section]
\newtheorem{remark}[theorem]{Remark}
\newtheorem{definition}[theorem]{Definition}
\newtheorem{lemma}[theorem]{Lemma}
\newtheorem{corollary}[theorem]{Corollary}

\DeclareMathOperator{\Exs}{\mathbb{E}}

\makeatletter
\newcommand{\removelatexerror}{\let\@latex@error\@gobble}
\makeatother
\usepackage[ruled,norelsize]{algorithm2e}

\begin{document}

\title{Distributed Sketching for Randomized Optimization: Exact Characterization, Concentration and Lower Bounds}

\author{Burak~Bartan\IEEEmembership{}
        and~Mert~Pilanci~\IEEEmembership{}
        %and~Jane~Doe,~\IEEEmembership{Life~Fellow,~IEEE}% <-this % stops a space
\thanks{B. Bartan is with the Department of Electrical Engineering, Stanford University, CA, 94305 USA (e-mail: bbartan@stanford.edu).}
\thanks{M. Pilanci is with the Department of Electrical Engineering, Stanford University, CA, 94305 USA (e-mail: pilanci@stanford.edu).}
% \thanks{M. Shell is with the Department
% of Electrical and Computer Engineering, Georgia Institute of Technology, Atlanta,
% GA, 30332 USA e-mail: (see http://www.michaelshell.org/contact.html).}% <-this % stops a space
% \thanks{J. Doe and J. Doe are with Anonymous University.}% <-this % stops a space
% \thanks{Manuscript received April 19, 2005; revised September 17, 2014.}
}

\maketitle


\begin{abstract}
    We consider distributed optimization methods for problems where forming the Hessian is computationally challenging and communication is a significant bottleneck. We leverage randomized sketches for reducing the problem dimensions as well as preserving privacy and improving straggler resilience in asynchronous distributed systems. We derive novel approximation guarantees for classical sketching methods and establish tight concentration results that serve as both upper and lower bounds on the error. We then extend our analysis to the accuracy of parameter averaging for distributed sketches. Furthermore, we develop unbiased parameter averaging methods for randomized second order optimization for regularized problems that employ sketching of the Hessian. Existing works do not take the bias of the estimators into consideration, which limits their application to massively parallel computation. We provide closed-form formulas for regularization parameters and step sizes that provably minimize the bias for sketched Newton directions. Additionally, we demonstrate the implications of our theoretical findings via large scale experiments on a serverless cloud computing platform.
\end{abstract}


% \fi
\IEEEpeerreviewmaketitle

\section{Introduction} \label{sec:introduction}


\IEEEPARstart{W}{e} investigate distributed sketching methods for solving large scale regression problems. In particular, we study parameter averaging for variance reduction and establish theoretical results on approximation performance. We consider a distributed computing model with a single central node and multiple worker nodes that run in parallel. Employing parameter averaging in distributed computing systems enables asynchronous updates, since a running average of available parameters can approximate the result without requiring all worker nodes to finish their tasks. Moreover, sketching provably preserves the privacy of the data, making it an attractive choice for massively parallel cloud computing.

We consider distributed averaging in a variety of randomized second order methods including Newton Sketch, iterative Hessian sketch (IHS), and also in direct or non-iterative methods. We focus on the communication-efficient setting where we avoid the communication of approximate Hessian matrices of size $d^2$ and communicate only the approximate solutions of size $d$. Averaging sketched solutions was proposed in the literature in certain restrictive settings \cite{wang2018distributed}. The presence of a regularization term requires additional caution, as na\"ive averaging may lead to biased estimators of the solution. In this work, bias is always with respect to the randomness of the sketching matrices; the input data are not assumed to be sampled from a probability distribution. Although bias is often overlooked in the literature, we show that one can re-calibrate the regularization coefficient of the sketched problems to obtain unbiased estimators. We show that having unbiased estimators leads to better performance without imposing any additional computational cost.

We consider both underdetermined and overdetermined linear regression problems in the regime where the data does not fit in main memory. Such linear regression problems and linear systems are commonly encountered in a multitude of problems ranging from statistics and machine learning to optimization. Being able to solve large scale linear regression problems efficiently is crucial for many applications. In this paper, the setting that we consider consists of a distributed system with $q$ workers that run in parallel and a single central node that collects and processes the outputs of the worker nodes. Applications of randomized sketches and dimension reduction to linear regression and other optimization problems have been extensively studied in the recent literature by works including \cite{drineas2011faster, mahoney2018averaging,PilWai14b,lacotte2019faster,Mahoney11,PilWai14a}. In this work, we investigate averaging the solutions of sketched sub-problems. This setting for overdetermined problems was also studied in \cite{mahoney2018averaging}. In addition, we also consider regression problems where the number of data samples is less than the dimensionality and investigate the properties of averaging for such problems.

An important advantage of randomized sketching in distributed computing is the independent and identically distributed nature of all of the computational tasks. Therefore, sketching offers a resilient computing model where node failures, stragglers as well as additions of new nodes can be easily handled, e.g., via generating more data sketches. An alternative to averaging that would offer similar benefits is the asynchronous stochastic gradient descent (SGD) algorithm \cite{bertsekas1989,recht2011hogwild}. However, the convergence rates of asynchronous SGD methods necessarily depend on the properties of the input data such as its condition number \cite{bertsekas1989,recht2011hogwild}. In contrast, as we show in this work, distributed sketching has stronger convergence guarantees that do not depend on the condition number of the data matrix.

% Our bias correction results have additional desirable properties for distributed computing systems. For heterogeneous distributed computing environments, where workers have varying computing capabilities, it might be advantageous for each worker to solve a problem of a different size \cite{avestimehr2017heterogenous_coded}. We provide formulas that specify the correct scaling for the regularization coefficient as a function of the problem size for each worker to obtain an unbiased estimator of the optimal solution.

A fundamental emphasis of the paper is on distributed computing and parameter averaging. However, novel theoretical results including exact characterizations of expected error and exponential concentration are derived for the single sketch estimator as well. 

%We present results for the single sketch estimator that are novel by themselves without the distributed averaging component.

Table \ref{table:summary_results} summarizes the types of problems that we study in this work. The last column of the table contains references to the theorems and lemmas for each result.

\begin{table*}
\caption{Summary of theoretical results. The data matrix is $A \in \mathbb{R}^{n \times d}$ and the output vector is $b \in \mathbb{R}^n$. %The second column 'Sketch' refers to the type of sketching matrix. The methods in the last two rows address a more general class of problems than the first three rows.
}
\label{table:summary_results}
\centering
\scalebox{1}{
 \begin{tabular}{|c|c|c|} 
 \hline
 Problem & Sketch Type & Method and Theorem \\ 
 \hline
 $\min_x \|Ax-b\|_2^2 $ & Gaussian &  Distributed randomized regression: Theorem \ref{Thm:AvgGaussian}, \ref{thm:concent_averaged_sketch}, \ref{thm:fisher_lower_bounds} \\ 
 & Gaussian & Distributed Iterative Hessian sketch: Theorem \ref{thm:IHS_error_decay} \\ 
 &  Other & Distributed randomized regression: Lemma \ref{bound_z_norm_ROS}, \ref{bound_z_norm_uniform}, \ref{bound_z_norm_leverage}, Theorem \ref{thm:errorbounds_other_sketches} \\
 \hline
 $\min_x \|x\|_2^2 \mbox{ s.t. } Ax=b $ & Gaussian & Right sketch: Theorem \ref{Thm:AvgGaussian_leastnorm}  \\
 \hline
 $\min_x \|Ax-b\|_2^2 + \lambda_1 \|x\|_2^2 $ & Gaussian & Distributed randomized ridge regression: Theorem \ref{thm:opt_lambda_2_newton_LS} \\
 \hline
 Convex problems with  & Gaussian & Distributed Newton sketch: Theorem \ref{thm:dist_newton_sketch} \\
 Hessian $(H_t^{1/2})^TH_t^{1/2}$ &  & (step size for unbiasedness) \\
 \hline
 Convex problems with & Gaussian & Distributed Newton sketch: Theorem \ref{opt_lambda_2_newton} \\
 Hessian $(H_t^{1/2})^TH_t^{1/2} + \lambda_1 I_d$ & & (regularization coefficient for unbiasedness) \\ [1ex] 
 \hline
\end{tabular}
}
\end{table*}

\subsection{Cloud Computing}
Serverless computing is a relatively new technology that offers computing on the cloud without requiring any server management from end users. Functions in serverless computing can be thought of worker nodes that have very limited resources and a lifetime, but also are very scalable. The algorithms that we study in this paper are particularly suitable for serverless computing platforms, since the algorithms do not require peer-to-peer communication among different worker nodes, and have low memory and compute requirements per node. In numerical simulations, we have evaluated our methods on the serverless computing platform AWS Lambda.

Data privacy is an increasingly important issue in cloud computing, one that has been studied in many recent works including \cite{zhou2009privacy,showkatbakhsh2018privacy,ZhoLafWas07,blocki12jltprivacy}. Distributed sketching comes with the benefit of privacy preservation. To clarify, let us consider a setting where the master node computes sketched data $S_kA$, $S_kb$, $k=1,\dots,q$ locally where $S_k \in \mathbb{R}^{m\times n}$ are the sketching matrices and $A \in \mathbb{R}^{n \times d}$, $b \in \mathbb{R}^{n}$ are the data matrix and the output vector, respectively, for the regression problem $\min_x\|Ax-b\|_2^2$. In the distributed sketching setting, the master node sends only the sketched data to worker nodes for computational efficiency, as well as data privacy preservation. In particular, the mutual information and differential privacy can be controlled when we reveal $S_kA$ and keep $A$ hidden. Furthermore, one can trade privacy for accuracy by choosing a suitable sketch dimension $m$. 
\subsection{Notation}
%We consider a distributed computing model where we have $q$ worker nodes which run in parallel and a single central node, i.e., master node. The worker nodes may only be allowed to communicate with the master node. The master node collects the outputs of the worker nodes and returns the averaged result. For iterative algorithms, this step serves as a synchronization point.

% Throughout the text, we provide exact formulas for the bias and variance of sketched solutions. 

We use hats $\hat{x}$ to denote the estimator for a single sketch and bars $\bar{x}$ to denote the averaged estimator $\bar{x} := \frac{1}{q} \sum_{k=1}^q \hat{x}_k$, where $\hat{x}_k$ is the estimator for the $k$'th worker node. Stars are used to denote the optimal solution $x^*$. We use $f(\cdot)$ to denote the objective of the optimization problem being considered at that point in the text. The letter $\epsilon$ is used for error while $\varepsilon$ is used as the differential privacy parameter. The notation $\sigma_{\min} (\cdot)$ denotes the smallest nonzero singular value of its argument. $O(\cdot)$ is used for big O notation.

All the expectations in the paper are with respect to the randomness over the sketching matrices, and no randomness assumptions are made for the data. We use $S \in \mathbb{R}^{m \times n}$ to denote random sketching matrices and we often refer to $m$ as the sketch size. For non-iterative distributed algorithms, we use $S_k \in \mathbb{R}^{m \times n}$ for the sketching matrix used by worker node $k$. For iterative algorithms, $S_{t,k} \in \mathbb{R}^{m \times n}$ is used to denote the sketching matrix used by worker $k$ in iteration $t$. We assume the sketching matrices are scaled such that $\Exs[S^TS]=I_n$. We omit the subscripts in $S_{k}$ and $S_{t,k}$ for simplicity whenever it does not cause confusion. The notation $m \gtrsim d$ is used to denote that there exists a positive finite constant $c$ such that $m \geq cd$. We denote the thin SVD decomposition of the data matrix as $A=U\Sigma V^T$.

For regularized problems, we use $\lambda_1$ for the regularization coefficient of the original problem, and $\lambda_2$ for the regularization coefficient of the sketched problem. For instance, in the case of the regularized least squares problem, we have
\begin{align} \label{eq:first_eq}
    x^* &= \arg\min_x \|Ax-b\|_2^2 + \lambda_1 \|x\|_2^2 \,, \nonumber \\
    \hat{x} &= \arg\min_x \|SAx-Sb\|_2^2 + \lambda_2 \|x\|_2^2 \,.
\end{align}
We will derive expressions in the sequel for the optimal selection of the coefficient $\lambda_2$ which we denote as $\lambda_2^*$.
\subsection{Related Work} \label{subsec:prior_work}

Random projections are a popular way of performing randomized dimensionality reduction, which are widely used in many computational and learning problems \cite{vempala2005random, mahoney2011randomized, woodruff2014sketching, drineas2016randnla}. Many works have studied randomized sketching methods for least squares and other optimization problems \cite{avron2010blendenpik, rokhlin2009randomized, drineas2011faster, pilanci2015randomized,pilanci2017newton,mahoney2018averaging,mahoney2018giant}.

Sarlós \cite{sarlos2006improved} showed that the relative error of the single Gaussian sketch estimator in \eqref{eq:first_eq} in the unregularized case $\lambda_1=\lambda_2=0$ is bounded as $f(\hat{x})/f(x^*) \leq (1+Cd\log(d)/m)^2 $ with probability at least $1/3$, where $C$ is a constant. The relative error was improved to $f(\hat{x})/f(x^*) \leq (1+Cd/m)$ with exponentially small failure probability in subsequent work \cite{PilWai14a}. In contrast, we derive the expectation of the error exactly and show exponential concentration around this expected value in this paper. As a result, we obtain a tight upper and lower bound for the relative error that holds with exponentially high probability (see Theorem \ref{thm:concent_averaged_sketch}).

The work \cite{mahoney2018averaging} investigates distributed sketching and averaging for regression from optimization and statistical perspectives. %Our work studies sketched model averaging from the optimization perspective. 
The most relevant result in \cite{mahoney2018averaging}, using our notation, can be stated as follows. Setting the sketch dimension $m=O(\mu d (\log d)/ \epsilon )$ for uniform sampling, where $\mu$ is row coherence $\mu:= n/\rank(A)\max_i \ell_i$ and $\ell_i$ is the $i$'th row leverage score, and $m=\tilde{O}(d/\epsilon )$ (with $\tilde{O}$ hiding logarithmic factors) for other sketches, the inequality $f(\bar{x}) - f(x^*) \leq (\epsilon/q+\epsilon^2 ) f(x^*)$ holds with high probability. According to this result, for large $q$, the cost at the averaged solution $f(\bar{x})$ will be upper bounded by $\epsilon^2 f(x^*)$. In this work we prove that for Gaussian sketch, $f(\bar{x})$ will converge to the optimal cost $f(x^*)$ as $q$ increases. We also identify the exact expected error for a given number of workers $q$. In addition, our results on regularized least squares regression improve on the results in \cite{mahoney2018averaging} for averaging multiple sketched solutions. In particular, in \cite{mahoney2018averaging}, the sketched sub-problems use the same regularization parameter as the original problem, which leads to biased solutions. We analyze the bias of the averaged solution, and provide explicit formulas for selecting the regularization parameter of the sketched sub-problems to achieve unbiasedness.

We show that the expected difference between the costs of the averaged solution and the optimal solution has two components, namely variance and squared bias (see Lemma \ref{expected_obj_val_diff}). This result implies that for the Gaussian sketch, which we prove to be unbiased (see Lemma \ref{gaussian_one_sketch}), the number of workers required for a target error $\epsilon$ scales as $1/\epsilon$. Remarkably, this result does not involve condition numbers. In contrast, for the asynchronous Hogwild algorithm \cite{recht2011hogwild}, the number of iterations required for error $\epsilon$ scales with $\log(1/\epsilon)/\epsilon$ and also depends on the condition number of the input data, although it addresses a more general class of problems.

%Our analysis approach allows the cost to be decomposed into bias and variance. 
% We derive upper bounds for the estimator biases for additional sketching matrices including randomized orthonormal systems (ROS) based sketch, uniform sampling, and leverage score sampling sketches.

%%%%%%%%%%%%%%%%%%%%%%

One of the crucial results of our work is developing bias correction for averaging sketched solutions. Namely, the setting considered in \cite{mahoney2018giant} is based on a distributed second order optimization method that involves averaging approximate update directions. However, in that work, the bias was not taken into account which degrades accuracy and limits applicability to massively parallel computation.
%and employs the same regularization parameters for the sketched sub-problems. 
We additionally provide results for the unregularized case, which corresponds to the distributed version of the Newton sketch algorithm introduced in \cite{pilanci2017newton}.
\subsection{Overview of Our Contributions} \label{subsec:contributions}
\begin{itemize}
    \item We characterize the exact expected error of the averaged estimator for Gaussian sketch in closed form as
    \begin{align}
        \frac{\Exs[f(\bar{x})]-f(x^*)}{f(x^*)} = \frac{1}{q} \frac{d}{m-d-1} 
        \label{eq:mean_avg_sketch}
    \end{align}
    where $\bar{x}=\frac{1}{q}\sum_{k=1}^q \hat{x}_k$ is the averaged solution and $\hat x_k = \arg\min_{x} \|S_kAx-S_kb\|^2_2$ is the output of the $k$'th worker node. In addition, we obtain a similar result for the error of the averaged estimator for the least-norm solution in the underdetermined case $n < d$.
    
    \item For both the single sketch and averaged estimator, we show that the relative error $\frac{f(\bar{x})-f(x^*)}{f(x^*)}$ is concentrated around its expectation given in \eqref{eq:mean_avg_sketch} with exponentially high probability. This provides both an upper and lower bound on the performance of sketching, unlike previous results in the literature. Moreover, it offers a guaranteed recipe to set the sketch size $m$ and number of workers $q$.
    
    \item We show that for Gaussian distributed sketch, the expected error of the averaged estimator $\Exs[f(\bar{x})]-f(x^*)$ matches the error lower bound for any unbiased estimator obtained via Fisher information. In addition, we provide a lower bound for general, possibly biased sketching based estimators.
    
    \item We consider the privacy preserving properties of distributed sketching methods, in which only random projections of the data matrix $A$ need to be shared with worker notes. 
    %This provides distributed sketching algorithms with privacy preservation. 
    We derive conditions on the $(\varepsilon, \delta)$-differential privacy when $S$ is the i.i.d. Gaussian sketch matrix. Combined with our results, we show that the approximation error of this distributed privacy preserving regression algorithm scales as $O(1/\varepsilon^2)$. %This is optimal as a result of our Fisher information based lower bound on distributed sketching.
    
    \item We analyze the convergence rate of a distributed version of the iterative Hessian sketch algorithm which was introduced in \cite{PilWai14b} and show that the number of iterations required to reach error $\epsilon$ with $q$ workers scales with $\log(1/\epsilon) / \log(q)$.

    \item We show that $\hat{x}= \arg\min_x \|SAx-Sb\|_2^2 + \lambda_2 \|x\|_2^2$ is not an unbiased estimator of the optimal solution, i.e. $\Exs[A(\hat{x}-x^*)] \neq 0$ when $\lambda_2 = \lambda_1$ and provide a closed form expression for the regularization coefficient $\lambda_2^*$ to make the estimator unbiased under certain assumptions.
    
    \item In addition to Gaussian sketch, we derive bias bounds for uniform sampling, randomized Hadamard sketch and leverage score sampling. Analysis of the bias term is critical in understanding how close to the optimal solution we can hope to get and establishing the dependence on the sketch dimension $m$. Moreover, we utilize the derived bias bounds and find an upper bound on the error of the averaged estimator for these other sketching methods.
    
    \item We discuss averaging for the distributed version of Newton Sketch \cite{pilanci2017newton} and show that using the same regularization coefficient as the original problem, i.e., $\lambda_2=\lambda_1$, as most works in the literature consider, is not optimal. Furthermore, we derive an expression for choosing the regularization coefficient $\lambda_2^*$ for unbiasedness.
    
    \item We provide numerical simulations that illustrate the practicality and scalability of distributed sketching methods on the serverless cloud computing platform AWS Lambda.
\end{itemize}

% \section{Preliminaries} \label{sec:prelim}

\subsection{Preliminaries on Sketching Matrices}
We consider various sketching matrices in this work including Gaussian sketch, uniform sampling, randomized Hadamard sketch, Sparse Johnson-Lindenstrauss Transform (SJLT), and \textit{hybrid} sketch. We now briefly describe each of these sketching methods:

\begin{enumerate}
    \item \textit{Gaussian sketch:} Entries of $S \in \mathbb{R}^{m\times n}$ are i.i.d. and sampled from the Gaussian distribution. Sketching a matrix $A \in \mathbb{R}^{n\times d}$ using Gaussian sketch requires computing the matrix product $SA$ which has computational complexity equal to $O(mnd)$.
    
    \item \textit{Randomized Hadamard sketch:} The sketch matrix in this case can be represented as $S=PHD$ where $P \in \mathbb{R}^{m\times n}$ is for uniform sampling of $m$ rows out of $n$ rows, $H \in \mathbb{R}^{n\times n}$ is the Hadamard matrix, and $D \in \mathbb{R}^{n\times n}$ is a diagonal matrix with diagonal entries sampled randomly from the Rademacher distribution. Multiplication by $D$ to obtain $DA$ requires $O(nd)$ scalar multiplications. Hadamard transform can be implemented as a fast transform with complexity $O(n\log(n))$ per column, and a total complexity of $O(nd\log(n))$ to sketch all $d$ columns of $DA$. %We note that because $P$ reduces the row dimension down to $m$, it might be possible to devise a more efficient way to perform sketching with lower computational complexity.
    
    \item \textit{Uniform sampling:} Uniform sampling randomly selects $m$ rows out of the $n$ rows of $A$ where the probability of any row being selected is the same.
    
    \item \textit{Leverage score sampling:} Row leverage scores of a matrix $A$ are given by $\ell_i = \|\tilde{u}_i\|_2^2$ for $i=1,\dots,n$ where $\tilde{u}_i \in \mathbb{R}^d$ denotes the $i$'th row of $U$. The matrix $U$ is the matrix whose columns are the left singular vectors of $A$, i.e., $A=U\Sigma V^T$. There is only one nonzero element in every row $s_i \in \mathbb{R}^n$ of the sketching matrix $S$ and the probability that the $j$'th entry of $s_i$ is nonzero is proportional to the leverage score $\ell_j$. More precisely, the rows $s_1,\dots,s_m$ are sampled i.i.d. such that $\mprob[s_{i}= e_k/\sqrt{m p_k}] = p_k,\,\forall i,\,\forall k$ where $p_k= \frac{\ell_k}{\sum_{j=1}^n \ell_j}$. Note that we have $\Exs [S^TS] = \Exs[ \sum_{i=1}^m s_is_i^T] = m \sum_{k=1}^n p_k e_ke_k^T/(mp_k)=I_n$.
    
    \item \textit{Sparse Johnson-Lindenstrauss Transform (SJLT) \cite{nelson2013osnap}:} The sketching matrix for SJLT is a sparse matrix where each column has exactly $s$ nonzero entries and the columns are independently distributed. The nonzero entries are sampled from the Rademacher distribution. It takes $O(snd/m)$ addition operations to sketch a data matrix using SJLT.
    
    \item \textit{Hybrid sketch:} The method that we refer to as hybrid sketch is a sequential application of two different sketching methods. In particular, it might be computationally feasible for worker nodes to sample as much data as possible, say $m^\prime$ rows, and then reduce the dimension of the available data to the final sketch dimension $m$ using another sketch with better error properties than uniform sampling such as Gaussian sketch or SJLT. For instance, a hybrid sketch of uniform sampling followed by Gaussian sketch would have computational complexity $O(mm^\prime d)$.
\end{enumerate}

\subsection{Paper Organization}

Section \ref{sec:quad_problems} deals with the application of distributed sketching to quadratic problems for Gaussian sketch. Section \ref{sec:privacy} deals with the privacy preserving property of distributed sketching. Section \ref{sec:other_sketching_matrices} gives theoretical results for randomized Hadamard sketch, uniform sampling, and leverage score sampling. The tools developed in Section \ref{sec:quad_problems} are applied to second order optimization for unconstrained convex problems in Section \ref{sec:nonlinear_problems}. Section \ref{sec:example_problems} provides an overview of applications and example problems where distributed sketching methods could be applied. Section \ref{sec:numerical_results} presents numerical results and Section \ref{sec:discussion} concludes the main part of the paper.

We give the proofs for the majority of the lemmas and theorems in the Appendix along with additional numerical results.

\section{Distributed Sketching for Quadratic Optimization problems} \label{sec:quad_problems}

In this section, we focus on quadratic optimization problems of the form
\begin{align} \label{eq:least_squares_problem}
    x^* = \arg \min_x \|Ax-b\|_2^2 + \lambda_1 \|x\|_2^2 \,.
\end{align}
We study various distributed algorithms for solving problems of this form based on model averaging. Some of these algorithms are tailored for the unregularized case $\lambda_1=0$.

\subsection{Closed Form Expressions for Gaussian Sketch} \label{subsec:closed_form_gaus}

We begin with a non-iterative model averaging based algorithm, which we refer to as distributed randomized regression. This algorithm is for the unregularized case $\lambda_1=0$. Each of the worker nodes computes an approximate solution $\hat{x}_k$ and these are averaged at the master node to compute the final solution $\bar{x}$. This method is outlined in Algorithm \ref{distrib_avg_alg}. The theoretical analysis in this section assumes that the sketch matrix is Gaussian, which is generalized to other sketching matrices in Section \ref{sec:other_sketching_matrices}. Although computing the Gaussian sketch is not as efficient as other fast sketches such as the randomized Hadamard sketch, it has several significant advantages over other sketches: (1) exact relative error expressions can be derived as we show in this section, (2) the solution is unbiased, (3) a differential privacy bound can be provided as we show in Section \ref{sec:privacy}, and (4) computation of the sketch can be trivially parallelized.

\begin{figure}[!t]
 \removelatexerror
  \begin{algorithm}[H]
  \caption{Distributed Randomized Regression}
  \label{distrib_avg_alg}
\begin{algorithmic}
\STATE{\textbf{Input:} Data matrix $A \in \mathbb{R}^{n\times d}$, target vector $b \in \mathbb{R}^{n}$.}
 \FOR{worker $k=1,\dots,q$ in parallel}
 \STATE{Obtain the sketched data and sketched output: $S_kA$ and $S_kb$.}
 \STATE{Solve $\hat x_k = \arg\min_{x} \|S_kAx-S_kb\|^2_2$ and send $\hat x_k$ to the master node.}
 \ENDFOR
 \STATE{\textbf{Master node:} return $\bar x = \frac{1}{q} \sum_{k=1}^q \hat x_k$.}
\end{algorithmic}
\end{algorithm}
\end{figure}

We first obtain a characterization of the expected error for the single sketch estimator in Lemma \ref{gaussian_one_sketch}.

\begin{lemma} \label{gaussian_one_sketch}
For the Gaussian sketch with $m>d+1$, the estimator $\hat{x}$ satisfies
\begin{align}
    \Exs [\| A (\hat{x} - x^*) \|_2^2] = \Exs[ f(\hat{x}) ] - f(x^*) =  f(x^*) \frac{d}{m-d-1}\,,
\end{align}
where $f(x) := \|Ax-b\|_2^2$ and $x^*=\arg\min_x f(x)$.
\end{lemma}
\begin{IEEEproof}[Proof of Lemma \ref{gaussian_one_sketch}]
Suppose that the matrix $A$ is full column rank. Then, for $m\ge d$, the matrix $A^T S^T S A$ follows a Wishart distribution, and is invertible with probability one. Conditioned on the invertibility of $A^T S^T S A$, we have
\begin{align*}
    \hat{x} &= (A^T S^T S A)^{-1} A^T S^T S b \\
    &= (A^T S^T S A)^{-1} A^T S^T S(Ax^* + b^\perp)\\
             &= x^* + (A^T S^T S A)^{-1} A^T S^T S b^\perp\,,
\end{align*}
where we have defined $b^\perp := b-Ax^*$\,. Note that $SA$ and $S b^\perp$ are independent since they are Gaussian and uncorrelated as a result of the normal equations $A^T b^\perp = 0$.
Conditioned on the realization of the matrix $SA$ and the event $A^T S^TS A \succ 0$, a simple covariance calculation shows that
\begin{align}
    \hat{x} \sim \mathcal{N}\big( x^*, \frac{1}{m} f(x^*) (A^T S^T S A)^{-1}\big)\,.
\end{align}
Multiplying with the data matrix $A$ on the left yields the distribution of the prediction error, conditioned on $SA$, as
\begin{align}
    A(\hat{x}-x^*) \sim \mathcal{N}\big( 0,  \frac{1}{m} f(x^*) A(A^T S^T S A)^{-1}A^T\big)\,.\label{eq:cond_Gaussian}
\end{align}
Then we can compute the conditional expectation of the squared norm of the error
\begin{align}
    \Exs [ \| A(\hat{x} - x^*) \|_2^2 ~\, \big \vert SA ] = \frac{f(x^*)}{m} \Exs [\tr( A(A^T S^T S A)^{-1}A^T)] \,.\nonumber
\end{align}
Next we recall that the expected inverse of the Wishart matrix $A^T S^T S A$ satisfies (see, e.g., \cite{letac2004all})
\begin{align}
    \Exs [(A^T S^T S A)^{-1}] = (A^T A)^{-1} \frac{m}{m-d-1}\,.\nonumber
\end{align}
Plugging in the previous result and using the tower property of expectations and then noting that $\tr(A (A^TA)^{-1} A^T) = d$ give us the claimed result.
\end{IEEEproof}

To the best of our knowledge, this result of exact error characterization is novel in the theory of sketching. Similar tools to those that we used in this proof were used in \cite{richardson2020sketching}, however, the exact context in which they are used and the end result are different from this work. Furthermore, existing results (see e.g. \cite{Mahoney11,mahoney2018averaging,woodruff2014sketching}) characterize a high probability upper bound on the error, whereas the above is a sharp and exact formula for the expected squared norm of the error. Theorem \ref{Thm:AvgGaussian} builds on Lemma \ref{gaussian_one_sketch} to characterize the expected error for the averaged solution $\bar{x}$.

\begin{theorem}[Expected error of the averaging estimator]
\label{Thm:AvgGaussian}
Let $S_k$, $k=1,\dots,q$ be Gaussian sketching matrices, then Algorithm \ref{distrib_avg_alg} runs in time $O(md^2)$, and the error of the averaged solution $\bar{x}$ satisfies
% Let each sketch dimension be chosen as $m=\frac{2rank(A)}{q\epsilon}$ using Gaussian sketches then the parallel scheme runs in time $O(\frac{rank(A)^2 d}{q \epsilon })$ and
\begin{align}
    \frac{\Exs[f(\bar{x})]-f(x^*)}{f(x^*)} = \frac{1}{q} \frac{d}{m-d-1}\,. 
\end{align}
%
%given that $A^TS_k^TS_kA \succ 0$ for all $k=1,\dots,q$.
Consequently, Markov's inequality implies that $
    \frac{f(\bar{x})-f(x^*)}{f(x^*)} \leq  \epsilon $
holds with probability at least $\left(1 -  \frac{1}{q\epsilon}\frac{d}{m-d-1} \right)$ for any $\epsilon>0$.
% \burak{Buradaki Markov's inequality ile yaptigimiz eski bound'i cikaralim mi?}
\end{theorem}

Theorem \ref{Thm:AvgGaussian} illustrates that the expected error scales as $1/q$, and converges to zero as $q\rightarrow 0$ as long as $m\ge d+2$. This is due to the unbiasedness of the Gaussian sketch. Other sketching methods such as uniform sampling or randomized Hadamard sketch do not have this property, as we further investigate in Section \ref{sec:other_sketching_matrices}.

\begin{remark}
In Algorithm \ref{distrib_avg_alg}, the worker nodes are tasked with obtaining the sketched data $S_kA$ and the sketched output $S_kb$. We identify two options for this step:
\begin{itemize}
    \item Option 1: The master node computes the sketched data $S_kA,S_kb$, $k=1,\dots,q$ and transmits to worker nodes. This option preserves data privacy, which we discuss in Section \ref{sec:privacy}.
    \item Option 2: The worker nodes have access to the data $A,b$ and compute the sketched data $S_kA$ and $S_kb$. This option does not preserve privacy, however, parallelizes the computation of the sketch via the workers. %In this case, it becomes very critical that we employ a low complexity sketching technique. This is discussed in Section \ref{sec:other_sketching_matrices}.
\end{itemize}
\end{remark}

%%%%%
\subsection{Exponential Concentration of the Gaussian Sketch Estimator} \label{subsec:highprob_bound}

In this subsection, we show high probability concentration bounds for the error of the Gaussian sketch. Importantly, our result provides both an upper and lower bounds on the relative error with exponentially small failure probability. We begin with the single sketch estimator $\hat{x}$ and extend the results to the averaged estimator $\bar{x}=\frac{1}{q}\sum_{k=1}^q \hat{x}_k$.

The next theorem states the main concentration result for the single sketch estimator. Both upper and lower bounds are given for the ratio of the cost of the sketched solution to the optimal solution. The full proof is given in the Appendix.

\begin{theorem} [Concentration bound for the single sketch estimator $\hat{x}$] \label{thm:concent_single_sketch}
Suppose that the sketch size is such that $m \gtrsim d$. Then, the optimality ratio of $\hat{x}$ with respect to the optimal solution $x^*$ is concentrated around its mean as
\begin{align} 
    P&\bigg(\Big| \frac{f(\hat{x})}{f(x^*)} - 1 - \frac{d}{m-d-1}\Big| < \epsilon \bigg) \geq 1 - C_1 e^{-C_2\epsilon^4 m}
\end{align}
% and
% \begin{align} 
%     P&\big(\frac{f(\hat{x})}{f(x^*)} < 1+\frac{d}{m-d-1} - \epsilon \big) \leq C_1 e^{-C_2\epsilon^4 m}
% \end{align}
%
where $C_1,C_2$ are positive constants. 
\end{theorem}

\noindent \emph{Proof sketch:}
Our proof technique involves the concentration of Gaussian quadratic forms for $Sb$ conditioned on $SA$. We leverage that the error $\|A(\hat x-x^*)\|_2^2$ is a Gaussian quadratic when conditioned on $SA$ as shown in \eqref{eq:cond_Gaussian}.  In particular, we use the concentration result given in Lemma \ref{lem:gausquad_concent} for quadratic form of Gaussian random variables.

\begin{lemma} [Concentration of Gaussian quadratic forms, \cite{boucheron2013concentration}] \label{lem:gausquad_concent}
Let the entries of $z \in \mathbb{R}^m$ be distributed as i.i.d. $\mathcal{N}(0,1)$. For any $G \in \mathbb{R}^{m\times m}$ and $\epsilon > 0$, 
\begin{align}
    P\big(z^TGz - \Exs[z^TGz] > 2\|G\|_F \sqrt{\epsilon} + 2\|G\|_2\epsilon \big) \leq e^{-\epsilon} \,.
\end{align}
\end{lemma}

Then we note that after applying Lemma \ref{lem:gausquad_concent}, the conditional expectation of the error is given by $\Exs [\|A(\hat x-x^*)\|_2^2\,|SA]=\frac{f(x^*)}{m} \tr((U^TS^TSU)^{-1})$. Next, we focus on the trace term $\tr((U^TS^TSU)^{-1})$ and relate it to the Stieltjes transform.

\begin{definition}[Stieltjes Transform]
We define the Stieltjes Transform for a random rectangular matrix $S \in \mathbb{R}^{m \times d}$ such that $d\le m$ as
\begin{align}
m_S(z):&=\frac{1}{d} \tr (S^T S-z I)^{-1}\\
& = \frac{1}{d} \sum_{i=1}^{d} \frac{1}{\lambda_{i}(S^TS)-z}\,.
\end{align}
Here $\lambda_i(S^TS)$ denotes the $i$'th eigenvalue of the symmetric matrix $S^TS$.
\end{definition}

We derive a high probability bound for the trace term around its expectation by leveraging the concentration of the empirical Stieltjes transform to its expectation, which is given in Lemma \ref{lem:trace_concent} below.

\begin{lemma} \label{lem:trace_concent}
The trace $\tr ( (U^TS^TSU)^{-1})$ is concentrated around its mean with high probability as follows
\begin{align} \label{eq:trace_concent_final}
    &P\Big(| \tr ( (U^TS^TSU)^{-1}) -\Exs[ \tr ( (U^TS^TSU)^{-1}) ] | \leq \epsilon \Big) \geq 1 - 4 e^{-\frac{\epsilon^4(1-\sqrt{d/m}-\delta)^8}{2^{10}md^2}} - e^{-m\delta^2/2}\,,
\end{align}
for any $\epsilon,\delta>0$.
\end{lemma}
The proof of Lemma \ref{lem:trace_concent} is provided in the Appendix. Combining the concentration of the Stieltjes transform with the concentration of Gaussian quadratic forms completes the proof of Theorem \ref{thm:concent_single_sketch}.

%%%%
Next, we show that the relative error of the averaged estimator is also concentrated with exponentially small failure probability. 

\begin{theorem} [Concentration bound for the averaged estimator $\bar{x}$] \label{thm:concent_averaged_sketch}
Let the sketch size satisfy $m \gtrsim d$. Then, the ratio of the cost for the averaged estimator $\bar{x}=\frac{1}{q}\sum_{k=1}^q \hat{x}_k$ to the optimal cost is concentrated around its mean as follows
\begin{align}
    &P\bigg(\Big| \frac{f(\bar{x})}{f(x^*)} - 1 - \frac{1}{q}\frac{d}{m-d-1} \Big| < \epsilon \bigg) \geq 1 - qC_1e^{-C_2(q\epsilon)^4m}  
\end{align}
where $C_1,C_2$ are positive constants.
\end{theorem}

The above result shows that the relative error of the averaged estimator using distributed sketching 
concentrates considerably faster compared to a single sketch. Moreover, the above bound offers a method to choose the values of the sketch size $m$ and number of workers $q$ to achieve a desired relative error with exponentially small error probability.

\begin{figure}[!t]
\begin{minipage}[b]{\linewidth}
  \centering
  \centerline{\includegraphics[width=0.42\columnwidth]{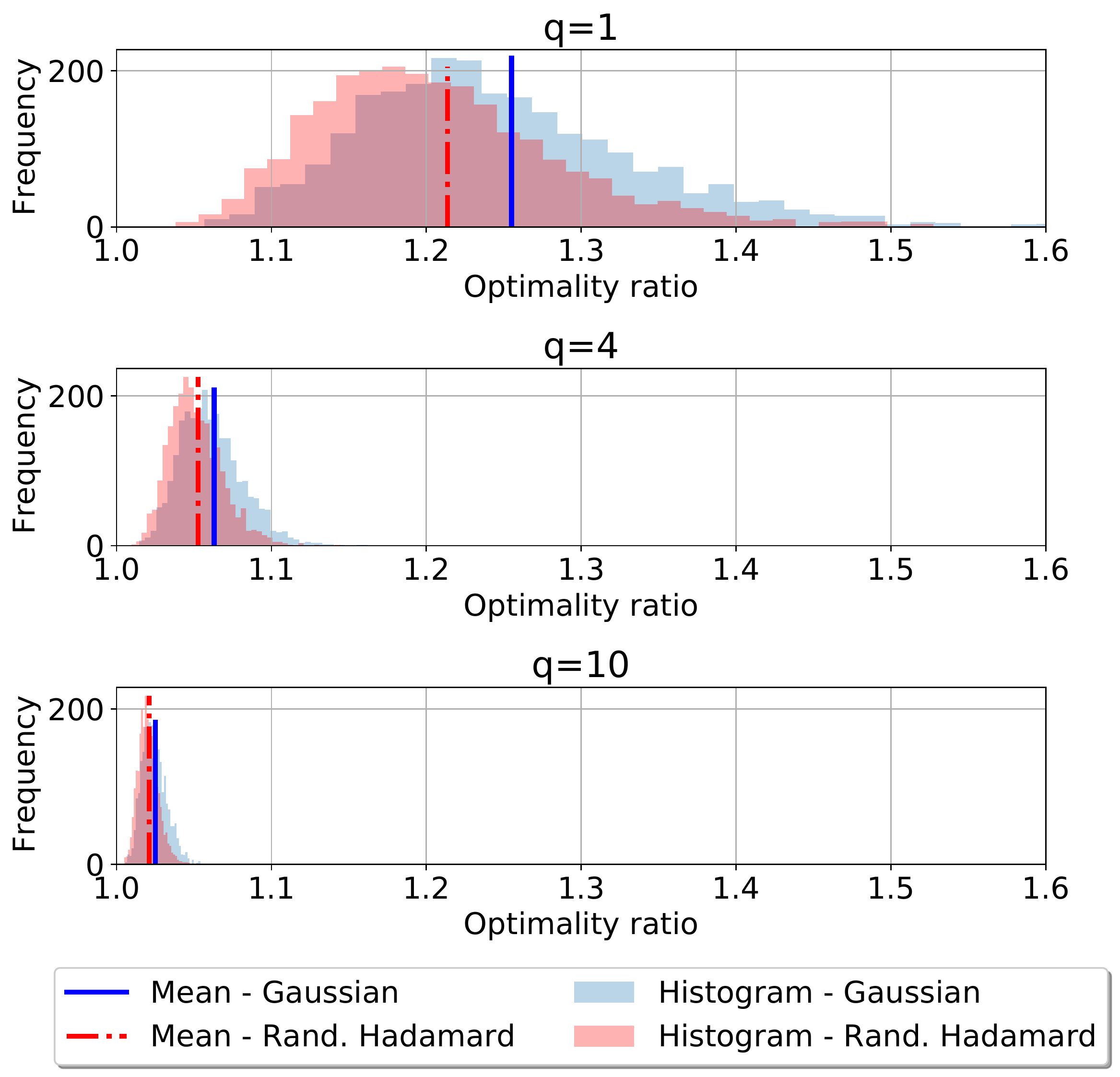}}
\end{minipage}
\caption{Histograms of the optimality ratio $f(\bar{x})/f(x^*)$ for both Gaussian sketch and randomized Hadamard sketch. The vertical lines indicate the empirical mean of the error. The dataset is synthetically generated with $n=512,d=20$ and the sketch size is $m=100$. The histogram is calculated over $2500$ independent trials. The plot at the top shows the performance of the single sketch estimator, and the middle and bottom plots are for the averaged estimator with $q=4$ and $q=10$ workers, respectively.%\mert{ Figure cok iyi, tesekkurler!}
}
\label{fig:error_histogram}
\end{figure}

Figure \ref{fig:error_histogram} shows the histograms for the single sketch ($q=1$) and averaged estimator ($q=4$ and $q=10)$ for Gaussian and randomized Hadamard sketches. This experiment demonstrates that the optimality ratio $f(\bar{x})/f(x^*)$ is concentrated around its mean and  its mean and concentration improve as we increase the number of workers $q$.

% \burak{Proof tamamlandi, appendix'e tasidim: "Proof of Theorem \ref{thm:concent_averaged_sketch}".}

% \burak{Boundlardaki $q$ carpani, union bound'dan dolayi geliyordu. Simdi union bound'siz yaptim. Bu sefer Bernoulli's inequality $(1-e^{-ms^2/2})^q \geq 1-qe^{-ms^2/2}$ yapinca ortaya cikti.}

% eski yazdiklarimizi gormek icin asagidaki satiri uncomment edebiliriz:
% \input{highprob_bound_OLD.tex}
\subsection{Error Lower Bounds via Fisher Information} \label{subsec:error_lowerbound}
We now present an error lower bound result for the Gaussian sketch. We first consider single sketch estimators and then discuss the distributed sketching case. Two different lower bounds can be obtained depending on whether the estimator is restricted to be unbiased or not. Lemma \ref{lem:error_lb_singlesketch} and \ref{lem:error_lb_singlesketch_biased} provide error lower bounds for all unbiased and general (i.e., possibly biased) estimators, respectively. The proof of Lemma \ref{lem:error_lb_singlesketch} is based on Fisher information and Cram\'er-Rao lower bound while Lemma \ref{lem:error_lb_singlesketch_biased} additionally employs the van Trees inequality. We refer the reader to \cite{sridhar20lowerbounds} for details on these single sketch lower bounds. 

\begin{lemma} [Unbiased estimators, \cite{sridhar20lowerbounds}] \label{lem:error_lb_singlesketch}
For any single sketch \textit{unbiased} estimator $\hat{x}$ obtained from the Gaussian sketched data $SA$ and $Sb$, the expected error is lower bounded as follows
\begin{align}
    \Exs [f(\hat{x})] - f(x^*) \geq f(x^*) \frac{d}{m-d-1}.
\end{align}
\end{lemma}

We note that the expected error of the single sketch estimator that we provide in Lemma \ref{gaussian_one_sketch} exactly matches the error lower bound in Lemma \ref{lem:error_lb_singlesketch}. Therefore, we conclude that no other unbiased estimator based on a single sketch $SA$, $Sb$ can achieve a better expected relative error. 
\begin{lemma} [General estimators, \cite{sridhar20lowerbounds}] \label{lem:error_lb_singlesketch_biased}
For any single sketch estimator $\hat{x}$, which is possibly biased, obtained from the Gaussian sketched data $SA$ and $Sb$, the expected error is lower bounded as follows
\begin{align}
    \Exs [f(\hat{x})] - f(x^*) \geq f(x^*) \frac{d}{m}.
\end{align}
\end{lemma}

%%%%% averaged estimator %%%%%
We now provide a novel generalized lower bound that applies to averaged estimators $\bar{x}$ obtained via distributed sketching. Theorem \ref{thm:fisher_lower_bounds} states the main lower bound result for the averaged estimator.

\begin{theorem} \label{thm:fisher_lower_bounds}
For any averaged estimator $\bar{x} = \frac{1}{q} \sum_{k=1}^q \hat{x}_k$, where each $\hat{x}_k$ is based on Gaussian sketched data $S_kA, S_kb$, $k=1,\dots,q$, the expected error is lower bounded as follows
\begin{enumerate}
\item[(i)] for unbiased estimators satisfying $\Exs[\hat{x}_k] = x^*$
\begin{align}
    \Exs [f(\bar{x})] - f(x^*) \geq \frac{f(x^*)}{q}\frac{d}{m-d-1}\,.
\end{align}

\item[(ii)] for general (i.e., possibly biased) estimators
\begin{align}
    \Exs [f(\bar{x})] - f(x^*) \geq \frac{f(x^*)}{q}\frac{d}{m}\,.
\end{align}
\end{enumerate}
\end{theorem}

We note that the expected error of the averaged estimator given in Theorem \ref{Thm:AvgGaussian} matches exactly the lower bound in Theorem \ref{thm:fisher_lower_bounds} for unbiased estimators. For general estimators, we observe that the lower bound matches the upper bound for large $m$, e.g., when $m-d-1 = O(m)$.

% proof of the theorem is moved to the appendix.

% Next, in the case of the averaged solution $\bar{x}$, an error lower bound can be obtained using Lemma \ref{lem:error_lb_singlesketch}. 

% \begin{lemma} \label{lem:error_lb_avgsol}
% For any set of $q$ unbiased estimators $\hat{x}_1,\dots,\hat{x}_q$ obtained from the Gaussian sketched data $S_1A,S_1y,\dots,S_qA,S_qy$ where the sketch matrices $S_1,\dots,S_q$ are independent, the error is lower bounded by
% \begin{align}
%     \Exs [\|A(\bar{x} - x^*) \|_2^2] \geq f(x^*) \frac{d}{mq-d-1} 
% \end{align}
% \end{lemma}

%%%%%

\subsection{Distributed Iterative Hessian Sketch} \label{subsec:distributed_ihs}
%In this section, we consider the well-known problem of unconstrained linear least squares which is stated as
% \begin{align} \label{eq:standard_least_sq}
%     x^* = \arg \min_x \frac{1}{2}\|Ax-b \|_2^2,
% \end{align}
% where $A \in \mathbb{R}^{n \times d}$ and $b \in \mathbb{R}^n$ are the problem data.
In this section, we consider an iterative algorithm for solving the unregularized least squares problem in \eqref{eq:least_squares_problem}, where $\lambda_1=0$ with higher accuracy. Note that, Newton's method terminates in one step when applied to this problem since the Hessian is $A^TA$ and
\begin{align}
    x_{t+1} = x_t - \mu (A^TA)^{-1}A^T(Ax_t-b)\,,
\end{align}
reduces to directly solving the normal equations $(A^TA)x=A^Tb$.
However, the computational cost of this direct solution is often prohibitive for large scale problems. To remedy this, the method of iterative Hessian sketch introduced in \cite{PilWai14b} employs a randomly sketched Hessian $A^TS_t^T S_t^T A$ as follows
\begin{align*}
    x_{t+1} = x_t - \mu (A^TS_t^TS_tA)^{-1}A^T(Ax_t-b),
\end{align*}
where $S_t$ corresponds to the sketching matrix at iteration $t$. Sketching reduces the row dimension of the data from $n$ to $m$ and hence computing an approximate Hessian $A^TS_t^TS_tA$ is computationally cheaper than the exact Hessian $A^TA$. Moreover, for regularized problems one can choose $m$ to be smaller than $d$ as we investigate in Section \ref{subsec:gaussian_bias_correction}.

In a distributed computing setting, one can obtain more accurate update directions by averaging multiple trials, where each worker node computes an independent estimate of the update direction. These approximate update directions can be averaged at the master node and the following update takes place
\begin{align} \label{averaged_ihs_update}
    x_{t+1} = x_t - \mu \frac{1}{q} \sum_{k=1}^q (A^TS_{t,k}^TS_{t,k}A)^{-1}A^T(Ax_t-b).
\end{align}
Here $S_{t,k}$ is the sketching matrix for the $k$'th node at iteration $t$. The details of the distributed IHS algorithm are given in Algorithm \ref{alg:dist_ihs}. %Note that in the algorithm we use $\arg\min$ instead of $(A^TS_{t,k}^TS_{t,k}A)^{-1}A^T(Ax_t-b)$. 
We note that although the update equation involves a matrix inverse term, in practice, this can be replaced with an approximate linear system solver. In particular, it might be computationally more efficient for worker nodes to compute their approximate update directions using indirect methods such as conjugate gradient.

We note that in Algorithm \ref{alg:dist_ihs}, worker nodes communicate their approximate update directions and not the approximate Hessian matrix, which reduces the communication complexity from $O(d^2)$ to $O(d)$ for each worker per iteration.

\begin{figure}[!t]
 \removelatexerror
  \begin{algorithm}[H]
  \caption{Distributed Iterative Hessian Sketch}
  \label{alg:dist_ihs}
  \begin{algorithmic}
  \STATE{\textbf{Input:} Number of iterations $T$, step size $\mu$.}
  \FOR{$t=1$ to $T$}
    \FOR{worker $k=1,\dots,q$ in parallel}
        \STATE{Obtain the sketched data $S_{t,k}A$.}
        \STATE{Compute gradient $g_t = A^T(Ax_t-b)$.}
        \STATE{Solve $\hat{\Delta}_{t,k}=\arg\min_{\Delta} \frac{1}{2} \|S_{t,k} A \Delta \|_2^2 + g_t^T\Delta $ and send to master node.}
    \ENDFOR
    \STATE{\textbf{Master node:} Update $x_{t+1} = x_t + \mu \frac{1}{q}\sum_{k=1}^q \hat{\Delta}_{t,k}$ and send $x_{t+1}$ to worker nodes.}
  \ENDFOR
  \RETURN{$x_T$}
  \end{algorithmic}
\end{algorithm}
\end{figure}

We establish the convergence rate for Gaussian sketch in Theorem \ref{thm:IHS_error_decay}, which provides an exact characterization of the expected error. First, we give the definition of the error in Definition \ref{error_vec_def}. 

\begin{definition} \label{error_vec_def}
To quantify the approximation quality of the iterate $x_t \in \mathbb{R}^d$ with respect to the optimal solution $x^* \in \mathbb{R}^d$, we define the error as $e^A_t \coloneqq A(x_t-x^*)$.
\end{definition}

To state the result, we first introduce the following moments of the inverse Wishart distribution (see Appendix).
\begin{align} \label{eq:theta_1_2_definitions}
        \theta_1 &\coloneqq \frac{m}{m-d-1}, \nonumber \\
     \theta_2 &\coloneqq \frac{m^2(m-1)}{(m-d)(m-d-1)(m-d-3)}.
\end{align}

\begin{theorem}[Expected error decay for Gaussian sketch] \label{thm:IHS_error_decay}
In Algorithm \ref{alg:dist_ihs}, let us set $\mu = 1/\theta_1$ and assume $S_{t,k}$'s are i.i.d. Gaussian sketches, then the expected squared norm of the error $e^A_{t}$ evolves according to the following relation:
\begin{align*}
\Exs [\| e^A_{t+1} \|_2^2] = \frac{1}{q} \left( \frac{\theta_2}{\theta_1^2} - 1 \right) \|e^A_t\|_2^2.
\end{align*}
\end{theorem}

Next, Corollary \ref{ihs_main_corollary} builds on Theorem \ref{thm:IHS_error_decay} to characterize the number of iterations for Algorithm \ref{alg:dist_ihs} to achieve an error of $\epsilon$. The number of iterations required for error $\epsilon$ scales with $\log(1/\epsilon)/\log(q)$.

\begin{corollary} \label{ihs_main_corollary}
Let $S_{t,k} \in \mathbb{R}^{m\times n}$ ($t=1,\dots,T$, $k=1,\dots,q$) be Gaussian sketching matrices. Then, Algorithm \ref{alg:dist_ihs} outputs $x_T$ that is $\epsilon$-accurate with respect to the initial error in expectation, that is, $\frac{\Exs [\|e^A_T\|_2^2}{\|Ax^*\|_2^2} = \epsilon$ where $T$ is given by
\begin{align*}
    T &= \frac{\log(1/\epsilon)}{\log(q)-\log\left(\frac{\theta_2}{\theta_1^2}-1\right)},
\end{align*}
where the overall required communication is $Tqd$ real numbers, and the computational complexity per worker is 
\begin{align*}
    O(Tmnd + Tmd^2 + Td^3).
\end{align*}
\end{corollary}
\begin{remark} 
Provided that $m$ is at least $2d$, the term $\log\left(\frac{\theta_2}{\theta_1^2}-1\right)$ is negative. Hence, $T$ is upper-bounded by $\frac{\log(1/\epsilon)}{\log(q)}$.
\end{remark}

\subsection{Distributed Sketching for Least-Norm Problems} \label{subsec:leastnorm}

In this section, we consider the underdetermined case where $n<d$ and applying the sketching matrix from the right, i.e., on the features. We will refer to this method as right sketch as well. Let us define the minimum norm solution
\begin{align}
     x^* = &\arg\min_x \|x\|_2^2 \quad \mbox{s.t.} ~ Ax = b. 
 \end{align}
The above problem has a closed-form solution given by $x^* = A^T (AA^T)^{-1}b$ when the matrix $A$ is full row rank. We will assume that the full row rank condition holds in the sequel. Let us denote the optimal value of the minimum norm objective as $f(x^*)=\|x^*\|_2^2 = b^T (AA^T)^{-1} b$. The $k$'th worker node will compute the approximate solution via $\hat x_k = S_k^T \hat{z}_k$ where
\begin{align}
    \hat{z}_k = &\arg \min_z \|z\|_2^2 \quad \mbox{s.t.} ~ AS_k^T z = b,
\end{align}
where $S_k \in \mathbb{R}^{m\times d}$ and $z \in \mathbb{R}^m$. The averaged solution is then computed at the master node as $\bar{x}=\frac{1}{q}\sum_{k=1}^q \hat{x}_k$.
%%% removed algorithm 2 to save space %%%
%The algorithm is listed in Alg. \ref{dist_rightsketch_alg}.
% \DontPrintSemicolon
% \begin{algorithm}
%  \KwIn{Data matrix $A \in \mathbb{R}^{n\times d}$, target vector $b \in \mathbb{R}^{n}$.}
%  \textbf{Workers $k=1,...,q$ in parallel: } \;
%  Sample $S_k \in \mathbb{R}^{m\times d}$. \;
%  Compute sketched data $AS_k^T$. \;
%  Solve $\hat{z}_k = \arg\min \|z\|_2^2$ subject to $ AS_k^Tz=b$. \;
%  Compute $\hat{x}_k = S_k^T \hat{z}_k$ and send to master. \;
%  \textbf{Master:} \;
%  return $\bar x \defn \frac{1}{q} \sum_{k=1}^q \hat x_k$. \;
%  \caption{Distributed right sketch for $n<d$ case.}
%  \label{dist_rightsketch_alg}
% \end{algorithm}
%\subsection{Gaussian Sketch}
We will assume that the sketch matrices $S_k$ are i.i.d. Gaussian in deriving the error expressions. Lemma \ref{gaussian_one_sketch_rightsketch} establishes the approximation error for a single right sketch estimator.

\begin{lemma} \label{gaussian_one_sketch_rightsketch}
For the Gaussian sketch with sketch size $m>n+1$, the estimator $\hat{x}$ satisfies
\begin{align*}
    \Exs [\| \hat{x} - x^* \|_2^2] = \frac{d-n}{m-n-1} f(x^*).
\end{align*}
\end{lemma}

\begin{IEEEproof}[Proof of Lemma \ref{gaussian_one_sketch_rightsketch}]
Conditioned on $AS^T$, we have 
\begin{align*}
    \hat{x} \sim \mathcal{N}\Big(x^*, P_{\mbox{Null}(A)}\| AS^T (AS^TSA^T)^{-1} b \|_2^2 \Big)\,.
\end{align*}
Noting that $\Exs [(AS^TSA^T)^{-1}] = AA^T \frac{m}{m-n-1} $, taking the expectation and noting that $\tr( P_{\mbox{Null}(A)})=d-n$, we obtain
\begin{align}
    \Exs [\| \hat{x} - x^* \|_2^2] &= \frac{d-n}{m-n-1} b^T (AA^T)^{-1} b = \frac{d-n}{m-n-1} f(x^*). 
\end{align} 
\end{IEEEproof}

An exact formula for averaging multiple outputs in right sketch that parallels Theorem \ref{Thm:AvgGaussian} can be obtained in a similar fashion: %We defer the details to the appendix. 
%
%\textbf{Error for the averaged solution:} Note that Lemma \ref{gaussian_one_sketch_rightsketch} establishes the expected error for the single sketch estimator $\hat{x}_k$. We now present an exact formula for the averaged solution $\bar{x}$. 
% We first note the unbiasedness of the estimates $\Exs[\hat{x}_k] = x^*$. Next, following similar steps as those used in proving Lemma \ref{expected_obj_val_diff}, we obtain
\begin{align*}
    \Exs[\|\bar{x} - x^* \|_2^2] &= \frac{1}{q} \Exs[\|\hat{x}_k - x^* \|_2^2] = \frac{1}{q} \frac{d-n}{m-n-1} f(x^*).
\end{align*}
Hence, for the distributed Gaussian right sketch, we establish the approximation error as stated in Theorem \ref{Thm:AvgGaussian_leastnorm}.
% \begin{align}
%     \frac{\Exs[f(\bar{x})] - f(x^*)}{f(x^*)} = \frac{1}{q} \frac{d-n}{m-n-1}.
% \end{align}

\begin{theorem}[Cost approximation for least-norm problems]
\label{Thm:AvgGaussian_leastnorm}
Let $S_k$, $k=1,\dots,q$ be Gaussian sketching matrices, then the error of the averaged solution $\bar{x}$ satisfies 
\begin{align}
    \frac{\Exs[f(\bar{x})]-f(x^*)}{f(x^*)} =  \frac{1}{q} \frac{d-n}{m-n-1}\,.
\end{align}
Consequently, Markov's inequality implies that $
    \frac{f(\bar{x})-f(x^*)}{f(x^*)} \leq  \epsilon $
holds with probability at least $\left(1 - \frac{1}{q\epsilon} \frac{d-n}{m-n-1} \right)$ for any $\epsilon > 0$.
\end{theorem}
\subsection{Bias Correction for Regularized Least Squares} \label{subsec:gaussian_bias_correction}

We have previously studied distributed randomized regression for unregularized problems and showed that Gaussian sketch leads to unbiased estimators. In this section, we focus on the regularized case and show that using the original regularization coefficient for the sketched problems causes the estimators to be biased. In addition, we provide a bias correction procedure. More precisely,
the method described in this section is based on non-iterative averaging for solving the linear least squares problem with $\ell_2$ regularization, i.e., ridge regression.%, and is fully asynchronous.

% We consider the problem given by
% \begin{align}
%     x^* = \arg\min_x \|Ax-b\|_2^2 + \lambda_1 \|x\|_2^2,
% \end{align}
% where $A \in \mathbb{R}^{n\times d}$, $b \in \mathbb{R}^n$ denote input data, and $\lambda_1 > 0$ is the regularization parameter.
% %
% Each worker applies sketching on $A$ and $b$ and obtains the estimate $\hat{x}_k$ given by
% \begin{align}
%     \hat{x}_k = \arg\min_x \|S_kAx-S_kb\|_2^2 + \lambda_2 \|x\|_2^2
% \end{align}
% for $k=1,\dots,q$,
% and the averaged solution is computed by the master node as
% \begin{align}
%     \bar{x} = \frac{1}{q} \sum_{k=1}^q \hat{x}_k.
% \end{align}
Note that we have $\lambda_1$ as the regularization coefficient of the original problem and $\lambda_2$ for the sketched sub-problems. If $\lambda_2$ is chosen to be equal to $\lambda_1$, then this scheme reduces to the framework given in the work of \cite{mahoney2018averaging} and we show in Theorem \ref{thm:opt_lambda_2_newton_LS} that $\lambda_2 = \lambda_1$ leads to a biased estimator, which does not converge to the optimal solution.

We first introduce the following results on traces involving random Gaussian matrices which are instrumental in our result.

\begin{lemma} [\cite{liu2019ridge}] \label{expectation_inverse_regularization}
For a Gaussian sketching matrix $S$, the following asymptotic formula holds
\begin{align*}
    \lim_{n \rightarrow \infty} \Exs[\tr( (U^TS^TSU + \lambda_2 I)^{-1} )] = d \times \theta_3(d/m, \lambda_2),
\end{align*}
where $\theta_3(d/m, \lambda_2)$ is defined as 
\begin{align*}
    \theta_3(d/m, \lambda_2) :=\frac{-\lambda_2+d/m-1 + \sqrt{(-\lambda_2+d/m-1)^2+4\lambda_2d/m}}{2\lambda_2 d/m}\,.
\end{align*}
\end{lemma}

\begin{lemma} \label{expectation_inverse_regularization_2}
For a Gaussian sketching matrix $S$, the following asymptotic formula holds
\begin{align*}
    \lim_{n \rightarrow \infty} \Exs[(U^TS^TSU + \lambda_2 I)^{-1}] = \theta_3(d/m, \lambda_2) I_d,
\end{align*}
where $\theta_3(d/m, \lambda_2)$ is as defined in Lemma \ref{expectation_inverse_regularization}.
\end{lemma}

We list the distributed randomized ridge regression method in Algorithm \ref{alg:dist_regularized_LS}. This algorithm assumes that our goal is to solve a large scale regression problem for a given value of regularization coefficient $\lambda_1$. Theorem \ref{thm:opt_lambda_2_newton_LS} states the main result of this subsection. In short, for the averaged result to converge to the optimal solution $x^*$, the regularization coefficient needs to be modified to $\lambda_2^*$.

\begin{theorem} \label{thm:opt_lambda_2_newton_LS}
Given the thin SVD decomposition $A=U\Sigma V^T \in \mathbb{R}^{n\times d}$ and $n \geq d$, and assuming $A$ has full rank and has identical singular values (i.e., $\Sigma = \sigma I_d$ for some $\sigma>0$), there is a value of $\lambda_2$ that yields a zero bias of the single sketch estimator $\Exs [A(\hat{x}_k - x^*)]$ as $n$ tends to infinity if 
\begin{itemize}
    \item[(i)] $m > d$\,\, or
    \item[(ii)] $m \leq d$ and $\lambda_1 \geq \sigma^2 \left( \frac{d}{m} - 1\right)$
\end{itemize}
and the value of $\lambda_2$ that achieves zero bias is given by
\begin{align} \label{opt_lambda_2_formula_LS}
    \lambda_2^* = \lambda_1 - \frac{d}{m}\frac{\lambda_1}{1+\lambda_1/\sigma^2},
\end{align}
where the matrix $S_k$ in $\hat{x}_k = \arg\min_x \|S_kAx-S_kb\|_2^2 + \lambda_2 \|x\|_2^2$ is the Gaussian sketch.
\end{theorem}

\begin{figure}[!t]
 \removelatexerror
  \begin{algorithm}[H]
  \caption{Distributed Randomized Ridge Regression}
  \label{alg:dist_regularized_LS}
  \begin{algorithmic}
  \STATE{Set $\sigma$ to the mean of singular values of $A$.}
  \STATE{Calculate $\lambda_2^* = \lambda_1 - \frac{d}{m}\frac{\lambda_1}{1+\lambda_1/\sigma^2}$.}
  \FOR{worker $k=1,\dots,q$ in parallel}
  \STATE{Obtain the sketched data and sketched output: $S_kA$ and $S_kb$.}
  \STATE{Solve $\hat{x}_k = \arg\min_x \|S_kAx-S_kb\|_2^2 + \lambda_2^* \|x\|_2^2$ and send $\hat x_k$ to master node.}
  \ENDFOR
  \STATE{\textbf{Master node:} return $\bar x = \frac{1}{q} \sum_{k=1}^q \hat x_k$.}
  \end{algorithmic}
\end{algorithm}
\end{figure}

\begin{figure}%[!t]
\hspace{2cm}
\begin{minipage}[b]{0.3\linewidth}
  \centering
  \centerline{\includegraphics[width=\columnwidth]{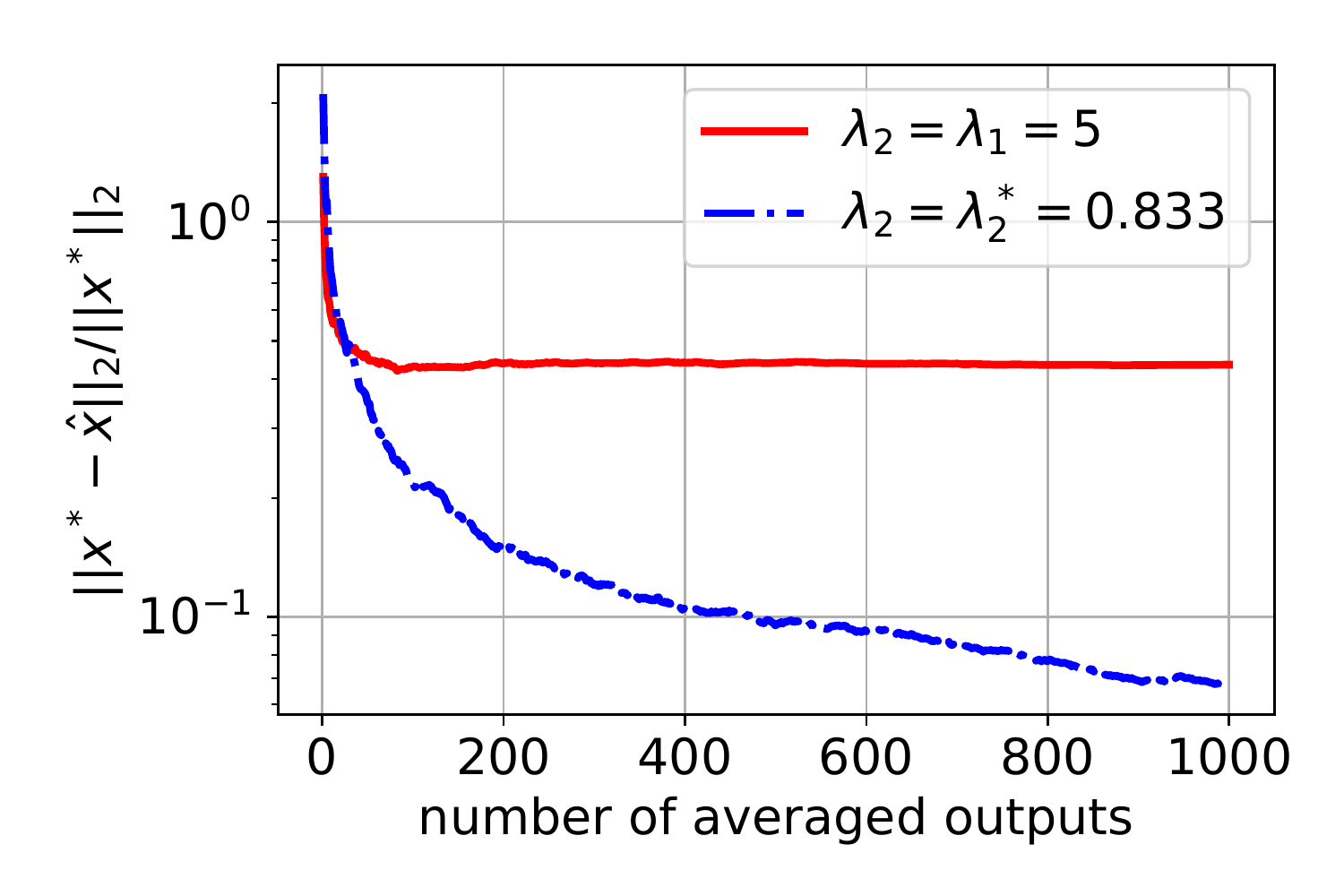}}
  \centerline{(a) Identical singular values }\medskip
\end{minipage}
% \hfill
\hspace{2cm}
\begin{minipage}[b]{0.3\linewidth}
  \centering
 \centerline{\includegraphics[width=\columnwidth]{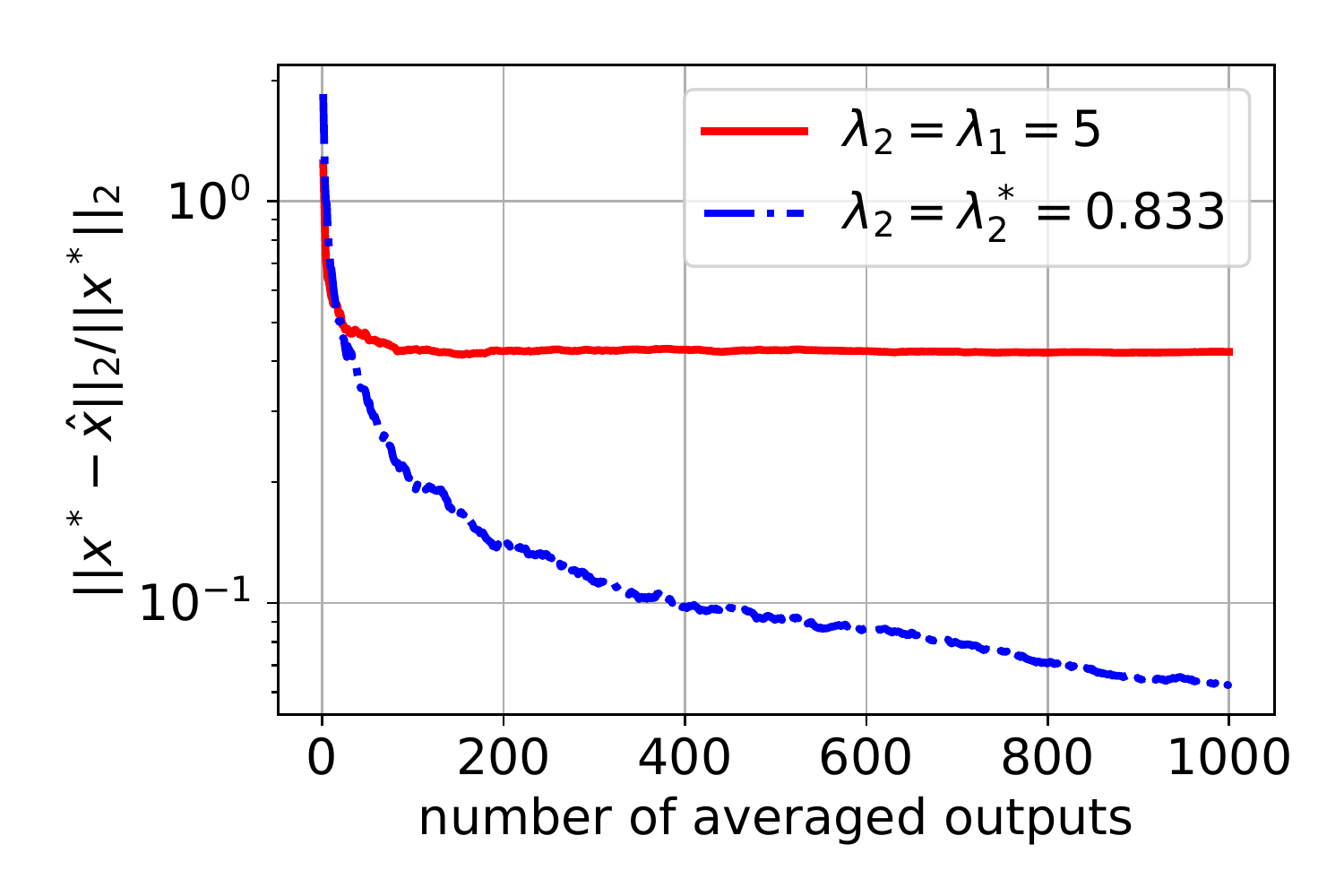}}
  \centerline{(b) Non-identical singular values}\medskip
\end{minipage}
\caption{Plots of $\|\bar{x}-x^*\|_2 / \|x^*\|_2$ against the number of averaged worker outputs for an unconstrained least squares problem with regularization using Algorithm \ref{alg:dist_regularized_LS}. The dashed blue line corresponds to the case where $\lambda_2$ is determined according to the formula \eqref{opt_lambda_2_formula_LS}, and the solid red line corresponds to the case where $\lambda_2$ is the same as $\lambda_1$. The parameters are as follows: $n=1000$, $d=100$, $\lambda_1=5$, $m=20$, sketch type is Gaussian. (a) All singular values of $A$ are $1$. (b) Singular values of $A$ are not identical and their mean is $1$.}
\label{fig:reg_ls_gaussketch}
\end{figure}

Figure \ref{fig:reg_ls_gaussketch} illustrates the practical implications of Theorem \ref{thm:opt_lambda_2_newton_LS}. If $\lambda_2$ is chosen according to the formula in \eqref{opt_lambda_2_formula_LS}, then the averaged solution $\bar{x}$ gives a significantly better approximation to $x^*$ than if we had used $\lambda_2=\lambda_1$. The data matrix $A$ in Figure \ref{fig:reg_ls_gaussketch}(a) has identical singular values, and \ref{fig:reg_ls_gaussketch}(b) shows the case where the singular values of $A$ are not identical. When the singular values of $A$ are not all equal to each other, we set $\sigma$ to the mean of the singular values of $A$ as a heuristic, which works extremely well as shown in Figure \ref{fig:reg_ls_gaussketch}(b).
%and in this case it has been set to $1$. 
According to the formula in \eqref{opt_lambda_2_formula_LS}, the value of $\lambda_2$ that we need to use to achieve zero bias is found to be $\lambda_2^*=0.833$ whereas $\lambda_1 = 5$. The plot in Figure \ref{fig:reg_ls_gaussketch}(b) illustrates that even if the assumption that $\Sigma=\sigma I_d$ in Theorem \ref{thm:opt_lambda_2_newton_LS} is violated, the proposed bias corrected averaging method outperforms vanilla averaging in \cite{mahoney2018averaging} where $\lambda_2=\lambda_1$.

% \subsubsection{Varying Sketch Sizes}
\begin{remark}[Varying sketch sizes]
Let us now consider the scenario where we have different sketch sizes for each worker node. This situation frequently arises in heterogeneous computing environments. Specifically, let us assume that the sketch size for worker $k$ is $m_k$, $k=1,\dots,q$. It follows from Theorem \ref{thm:opt_lambda_2_newton_LS} that by choosing the regularization parameter for worker node $k$ as
\begin{align*}
    \lambda_2^*(k) = \lambda_1 - \frac{d}{m_k}\frac{\lambda_1}{1+\lambda_1/\sigma^2},
\end{align*}
it is possible to obtain unbiased estimators $\hat{x}_k$, $k=1,\dots,q$ and hence an unbiased averaged result $\bar{x}$. Note that here we assume that the sketch size for each worker satisfies the condition in Theorem \ref{thm:opt_lambda_2_newton_LS}, that is, either $m_k > d$ or $m_k \leq d$ and $\lambda_1 \geq \sigma^2 (d/m_k-1)$.
\end{remark}

\section{Privacy Preserving Property} \label{sec:privacy}

We now digress from the error and convergence properties of distributed sketching methods to consider the privacy preserving properties of distributed sketching. We use the notion of differential privacy for our privacy result given in Definition \ref{def:diff_privacy}. Differential privacy is a worst case type of privacy notion which does not require distribution assumptions on the data. It has been the privacy framework adopted by many works in the literature. %Theorem \ref{thm:diff_priv_cited} from \cite{sheffet15privacy} indicates that publishing the product of a data matrix and a Gaussian random projection matrix preserves $(\varepsilon, \delta)$-differential privacy if the minimum singular value of the data matrix is at least $w$, which is defined in \eqref{eq:def_w}.

\begin{definition}[$(\varepsilon, \delta)$-Differential Privacy, \cite{dwork06privacy}, \cite{sheffet15privacy}] \label{def:diff_privacy}
An algorithm ALG which maps $(n\times d)$-matrices into some range
$\mathcal{R}$ satisfies $(\varepsilon, \delta)$-differential privacy if for all pairs of neighboring inputs $A$ and $A^\prime$ (i.e. they differ only in a single row) and all subsets $\mathcal{S} \subset \mathcal{R}$, it holds that $P(\text{ALG}(A) \in \mathcal{S}) \leq e^\varepsilon P(\text{ALG}(A^\prime) \in \mathcal{S}) + \delta$. %When $\delta=0$, we say the algorithm is $\varepsilon$-differentially private.
\end{definition}

Theorem \ref{thm:diff_priv_result} characterizes the required conditions and the sketch size to guarantee $(\varepsilon, \delta)$-differential privacy for a given dataset. The proof of Theorem \ref{thm:diff_priv_result} is based on the $(\varepsilon, \delta)$-differential privacy result for i.i.d. Gaussian random projections from \cite{sheffet15privacy}.
%We note that the values $m$ given in Theorem \ref{thm:diff_priv_result} are upper bounds for $m$, and ...

\begin{theorem} \label{thm:diff_priv_result}
Given the data matrix $A \in \mathbb{R}^{n\times d}$ and the output vector $b \in \mathbb{R}^{n}$, let $A_c$ denote the concatenated version of $A$ and $b$, i.e., $A_c = [A, b] \in \mathbb{R}^{n\times {(d+1)}}$. Define $\beta := \ln(4/\delta)$. Define $\sigma_0:=\sigma_{\min}(A_c)/ \sqrt{n}$, and $B_0:=\max_{i,j} |A_{c,ij}|$. Suppose that the condition
\begin{align} \label{eq:diff_priv_condition}
    \frac{n}{d+1} \geq \left( 3 + \frac{2\beta}{\varepsilon} \right) \frac{B_0^2}{\sigma_0^2}
    %\sigma_{\min} \geq B \sqrt{3 + \frac{2\beta}{\varepsilon}}
\end{align}
is satisfied and the sketch size satisfies
\begin{align} \label{eq:sketch_size_expressions}
    m &= O\left(\beta \frac{n^2}{(d+1)^2} \frac{\varepsilon^2}{(\varepsilon + \beta)^2}\right) \,\text{for privacy w.r.t. one worker,} \\
    %\left\lfloor \frac{\beta}{8} \left(\left(\frac{\sigma_{\min}^2}{B^2} - 1\right) \frac{\varepsilon }{\varepsilon + \beta} - 2 \right)^2 \right\rfloor \quad\text{for privacy w.r.t. a single worker,} \\
    m &= O\left(\frac{\beta}{q} \frac{n^2}{(d+1)^2} \frac{\varepsilon^2}{(\varepsilon + \beta)^2}\right) \,\text{for privacy w.r.t. $q$ workers,}
    %\left\lfloor \frac{1}{q}\frac{\beta}{8} \left(\left(\frac{\sigma_{\min}^2}{B^2} - 1\right) \frac{\varepsilon }{\varepsilon + \beta} - 2 \right)^2 \right\rfloor \quad\text{for privacy w.r.t. $q$ workers},
\end{align}
and publishing the sketched matrices $S_kA_c \in \mathbb{R}^{m\times {(d+1)}}$, $k=1,\dots,q$, provided that $d+1 < m$, is $(\varepsilon, \delta)$-differentially private, where $S_k \in \mathbb{R}^{m\times n}$ is the Gaussian sketch, and $\varepsilon > 0$, $\beta > 1+\ln(4)$.
\end{theorem}

%According to Theorem \ref{thm:diff_priv_result}, provided that the conditions are satisfied, we can pick the sketch size values as given in \eqref{eq:sketch_size_expressions} for $(\varepsilon, \delta)$-differential privacy. 
% Burayi cikarabiliriz
%In the cases where one is not able to pick $m$ that would satisfy the privacy requirements, the authors in \cite{sheffet15privacy} utilize the method of regularization as a way of ensuring privacy, which corresponds to appending a scaled identity matrix to the data matrix. The scaling factor of the identity is chosen to satisfy the differential privacy requirements.

% \bremark
% Theorem \ref{thm:diff_priv_result} in \eqref{eq:diff_priv_condition} implies that if $\frac{n}{d}$ is on the order of $O(\frac{\beta}{\varepsilon})$, then we can pick $m = O\left(\frac{\beta}{q} \frac{n^2}{d^2} \frac{\varepsilon^2}{(\varepsilon + \beta)^2}\right)$ for $(\varepsilon, \delta)$-differential privacy with respect to all $q$ workers.
% \eremark

\begin{remark}
For fixed values of $\beta$, $\sigma_{\min}$, $B_0$, $n$ and $d$, the approximation error is on the order of $O\left(\frac{1}{\varepsilon^2} \right)$ for $(\varepsilon, \delta)$-differential privacy with respect to all $q$ workers.
\end{remark}

The work \cite{huang15distprivacy} considers convex optimization under privacy constraints and shows that  $\varepsilon$-differential privacy (i.e. equivalent to Definition \ref{def:diff_privacy} with $\delta=0$), the approximation error of their distributed iterative algorithm is on the order of $O(\frac{1}{\varepsilon^2})$, which is on the same order as Algorithm \ref{distrib_avg_alg}. The two algorithms however have different dependencies on parameters, which are hidden in the $O$-notation. %For instance, the error of the algorithm in \cite{huang15distprivacy} has dependence on the order of $O(\sigma_{\max}^4)$ while Algorithm \ref{distrib_avg_alg} has dependence on the order of $O(B^4 / \sigma_{\min}^4)$. 
We note that the approximation error of the algorithm in \cite{huang15distprivacy} depends on parameters that Algorithm \ref{distrib_avg_alg} does not have such as the initial step size, the step size decay rate, and noise decay rate. The reason for this is that the algorithm in \cite{huang15distprivacy} is a synchronous iterative algorithm designed to solve a more general class of optimization problems. Algorithm \ref{distrib_avg_alg}, on the other hand, is designed to solve linear regression problems and offers a significant advantage due to its single-round communication requirement.
%and this makes it possible to design Algorithm \ref{distrib_avg_alg} as a single-round communication algorithm. This also establishes Algorithm \ref{distrib_avg_alg} as more robust against slower nodes.

\begin{remark}
A similar privacy statement holds for the right sketch method discussed in Section \ref{subsec:leastnorm}. In right sketch, we only sketch the data matrix $A$ and not the output vector $b$. For publishing $AS_k^T$ to be $(\varepsilon, \delta)$-differentially private, Theorem \ref{thm:diff_priv_result} still holds with the modification that we replace $A_c$ with $A^T$.
\end{remark}

% S is scaled R
% used \varepsilon in two different places
% consider all workers ++
% A should be: concatenated [A, y]: n x (d+1) ++
% m is an integer ++

%%%% papers on differential privacy:
% https://ieeexplore.ieee.org/stamp/stamp.jsp?arnumber=8437722
% https://arxiv.org/pdf/0901.1365.pdf
% https://arxiv.org/pdf/1204.2136.pdf
% https://arxiv.org/pdf/1507.02482.pdf
% Private Approximations of the 2nd-Moment Matrix Using Existing Techniques in Linear Regression: https://arxiv.org/abs/1507.00056 

%%%% papers 06/02 %%%%
% Differentially Private Distributed Optimization
% Randomness Efficient Fast-Johnson-Lindenstrauss Transform with Applications in Differential Privacy
% A Private and Finite-Time Algorithm for Solving a Distributed System of Linear Equations
% Privacy-Preserving Distributed Optimization via Subspace Perturbation: A General Framework
% Acceleration and Privacy Protection for Distributed Computing

\section{Bounds for Other Sketching Matrices}
\label{sec:other_sketching_matrices}

In this section, we consider randomized Hadamard sketch, uniform sampling, and leverage score sampling for solving the unregularized linear least squares problem $\min_x f(x)=\|Ax-b\|_2^2$. For each of these sketching methods, we present upper bounds on the bias of the corresponding single sketch estimator. Then, we combine these in Theorem \ref{thm:errorbounds_other_sketches}, which provides high probability upper bounds for the error of the averaged estimator.
%In particular, the results of this section focus on the bias bounds for a single output $\hat{x}$, and the way these results are related to the averaged solution $\bar{x}$ is through the decomposition given in Lemma \ref{expected_obj_val_diff}.

% We make use of the following identity in the sequel: The expected difference between the costs of the averaged solution $f(\bar{x})$ and the optimal solution $f(x^*)$ is equal to:
% \begin{align} \label{eq:expected_diff}
%     \Exs[f(\bar{x})]-f(x^*)% &= \Exs[\|A\bar{x}-b\|_2^2]-f(x^*) \\
%     &= \Exs[\|A(\bar{x}-x^*)+Ax^*-b\|_2^2]-f(x^*) \nonumber \\
%     & = \Exs [\|A(\bar{x}-x^*)\|_2^2 + \|Ax^*-b\|_2^2]-f(x^*) \nonumber \\
%     & = \Exs [\|A(\bar{x}-x^*)\|_2^2]\,,
% \end{align}
% where we have used the orthogonality property of the optimal least squares solution $x^*$ given by the normal equations $A^T (Ax^*-b) = 0$.%, which implies $(\bar x - x^*)^T A^T (Ax^*-b) = 0$.

Lemma \ref{expected_obj_val_diff} expresses the expected objective value difference in terms of the bias and variance of the single sketch estimator. 

%For the sketch types studied in this section, it is possible to obtain high probability bounds for the error $(f(\bar{x})-f(x^*))/f(x^*)$ based on the bias bounds given in this section, using an argument similar to the one given in the proof of Theorem \ref{Thm:AvgGaussian}. This approach would involve defining an event that bounds the error for the single sketch estimator. 

\begin{lemma} \label{expected_obj_val_diff}
For any i.i.d. sketching matrices $S_k$, $k=1,\dots,q$, the expected objective value difference for the averaged estimator $\bar{x}$ can be decomposed as
\begin{align} \label{approx_guarantee}
    &\Exs[f(\bar{x})]-f(x^*) = \frac{1}{q} \Exs \left[ \|A\hat{x}-Ax^* \|_2^2 \right] + \frac{q-1}{q} \| \Exs[A\hat{x}] - Ax^* \|_2^2,
\end{align}
where $\hat{x}$ is the single sketch estimator.
\end{lemma}

Lemma \ref{norm_of_bias} establishes an upper bound on the norm of the bias for any i.i.d. sketching matrix. The results in the remainder of this section will build on this lemma.

\begin{lemma} \label{norm_of_bias}
Let the SVD of $A$ be denoted as $A = U\Sigma V^T$. Let $z := U^TS^TSb^\perp$ and $Q := (U^TS^TSU)^{-1} $ where $b^\perp := b - Ax^*$. Define the event $E$ as $(1-\epsilon)I_d \preceq Q \preceq (1+\epsilon)I_d$. Then, for any sketch matrix $S$ with $\Exs[S^TS]=I_n$, the bias of the single sketch estimator, conditioned on $E$, is upper bounded as
\begin{align}
    \|\Exs[A\hat{x} | E] - Ax^*\|_2 \leq \sqrt{4\epsilon \Exs[\| z\|_2^2|E]} \,.
\end{align}
\end{lemma}

The event $E$ is a high probability event and it is equivalent to the subspace embedding property (e.g. \cite{mahoney2018averaging}). We will analyze the unconditioned error later in Theorem \ref{thm:errorbounds_other_sketches}.

We note that Lemma \ref{expected_obj_val_diff} and \ref{norm_of_bias} apply to all of the sketching matrices considered in this work. We now give specific bounds for the bias of the single sketch estimator separately for each of randomized Hadamard sketch, uniform sampling, and leverage score sampling.
%In terms of the computational complexity required to form a sketch of the data, uniform sampling is the fastest method while its bias is the largest as we show in the sequel. 
%This section lists all the derived bounds and the proofs are deferred to the Appendix.

% \subsection{Randomized Hadamard Sketch}
%Randomized Hadamard sketch can be represented as $S=PHD$ where $P$ is for sampling $m$ rows out of $n$, $H$ is the ($n\times n$)-dimensional Hadamard matrix, and $D$ is a diagonal matrix with diagonal entries sampled from the Rademacher distribution. 
\textbf{Randomized Hadamard sketch:} This method has the advantage of low computational complexity due to the Fast Hadamard Transform, which can be computed in $O(n\log n)$, whereas applying the Gaussian sketch takes $O(mn)$ time for length $n$ vectors. Lemma \ref{bound_z_norm_ROS} states the upper bound for the bias for randomized Hadamard sketch.

\begin{lemma} \label{bound_z_norm_ROS}
For randomized Hadamard sketch, the bias is upper bounded as
\begin{align}
    \|\Exs[A\hat{x}|E] - Ax^*\|_2 \leq \sqrt{4\epsilon \frac{d}{m} f(x^*)}.
\end{align}
\end{lemma}

% \subsection{Uniform Sampling}
%The sketching matrix $S$ for uniform sampling is equal to $P$. 
% In uniform sampling, each row of $S$ consists of a single $1$ and $(n-1)$ $0$'s and the position of $1$ in every row is distributed according to uniform distribution. Then, the $1$'s in $S$ are scaled so that $\Exs[S^TS]=I_n$.

\textbf{Uniform sampling:} We note that the bias of the uniform sampling estimator is different when it is computed with or without replacement. The reason for this is that the rows of the sketching matrix $S$ for uniform sampling are independent in the case of sampling with replacement, which breaks down in the case of sampling without replacement. 
Lemma \ref{bound_z_norm_uniform} provides bounds for both cases.

\begin{lemma} \label{bound_z_norm_uniform}
For uniform sampling, the bias can be upper bounded as
\begin{align}
    &\|\Exs[A\hat{x}|E] - Ax^*\|_2 \leq \sqrt{4\epsilon \frac{\mu d}{m}\, f(x^*) } \\
    &\|\Exs[A\hat{x}|E] - Ax^*\|_2 \leq \sqrt{4\epsilon \frac{\mu d}{m} \frac{n-m}{n-1}\, f(x^*) },
\end{align}
for sampling with and without replacement, respectively, where $\tilde{u}_i \in \mathbb{R}^d$ denotes the $i$'th row of $U$.
\end{lemma}

\textbf{Leverage score sampling:} Recall that the row leverage scores of a matrix $A=U\Sigma V^T$ are computed via $\ell_i = \|\tilde{u}_i\|_2^2$ for $i=1,\dots,n$ where $\tilde{u}_i \in \mathbb{R}^d$ denotes the $i$'th row of $U$. Since leverage score sampling looks at the data to help guide its sampling strategy, it achieves a lower error compared to uniform sampling which treats all the samples equally.
% The probability that the $j$'th entry of $s_i$ is nonzero is proportional to the leverage score $\|\tilde{u}_j\|_2^2$, i.e., $\mprob[s_{ij}\neq 0] = \frac{\|\tilde{u}_j\|_2^2}{\sum_{j=1}^n \|\tilde{u}_j\|_2^2}$. The denominator is equal to $d$ since it is the Frobenius norm of $U$ and the columns of $U$ are normalized. Hence, more concisely, we have $\mprob[s_{ij}\neq 0] = \frac{1}{d}\|\tilde{u}_j\|_2^2=\frac{1}{d}\ell_i$. 
Lemma \ref{bound_z_norm_leverage} gives the upper bound for the bias of the leverage score sampling estimator.

\begin{lemma} \label{bound_z_norm_leverage}
For leverage score sampling, the bias can be upper bounded as
\begin{align}
    \|\Exs[A\hat{x}|E] - Ax^*\|_2 \leq \sqrt{4\epsilon \frac{d}{m} f(x^*)}.
\end{align}
\end{lemma}

The results given in Lemma \ref{bound_z_norm_ROS}, \ref{bound_z_norm_uniform}, and \ref{bound_z_norm_leverage} can be combined with Lemma \ref{expected_obj_val_diff}. In particular, the error of the averaged estimator contains the squared bias term which is scaled with $(q-1)/q$. For a distributed computing system with a large number of worker nodes $q$, the contribution of the bias term is very close to $1$ while the variance term will vanish as it is scaled with $1/q$.

%%%%%%%%%%%%%%%%%%%%%%%%%%
We now leverage our analysis of the bias upper bounds to establish upper bounds on the error of the averaged estimator. Theorem \ref{thm:errorbounds_other_sketches} gives high probability error bounds for different sketch types when the sketch size is set accordingly.
\begin{theorem} \label{thm:errorbounds_other_sketches}
Suppose that the sketch size is selected as
\begin{align}
    m &\gtrsim \frac{d+\log(n)}{\epsilon^2} \log(d/\delta)  \mbox{  for randomized Hadamard sketch,} \nonumber \\
    m &\gtrsim \frac{\mu d}{\epsilon^2} \log(d/\delta)  \mbox{  for uniform sampling,} \nonumber \\
    m &\gtrsim \frac{d}{\epsilon^2} \log(d/\delta)  \mbox{  for leverage score sampling},
\end{align}
for some $\epsilon\in(0,\frac{1}{4}]$ and $\delta>0$.
Here $\mu$ is the row coherence defined in Section \ref{subsec:prior_work}, that is, $\mu=n/d \max_i \|\tilde{u}_i\|_2^2$. Then, the relative optimality gap of the averaged estimator obeys the following upper bounds:
\begin{align}
    &P\Big(\frac{f(\bar{x})}{f(x^*)} \leq 1 + \gamma \Big) \geq \nonumber \\
    &1 - q\delta - \frac{d}{q\gamma m} \big((1+\epsilon)^2 + 4\epsilon(q-1)\big) \mbox{  for randomized Hadamard sketch,} \\
    & 1 - q\delta - \frac{\mu d}{q\gamma m} \big((1+\epsilon)^2 + 4\epsilon(q-1)\big)  \mbox{  for uniform sampling with replacement,} \\
    & 1 - q\delta - \frac{\mu d}{q\gamma m} \frac{n-m}{n-1} \big((1+\epsilon)^2 + 4\epsilon(q-1)\big)  \mbox{  for uniform sampling without replacement,} \\
    & 1 - q\delta - \frac{d}{q\gamma m} \big((1+\epsilon)^2 + 4\epsilon(q-1)\big) \mbox{  for leverage score sampling,}
\end{align}
for any $\gamma>0$.
\end{theorem}

Remarkably, the optimality ratio of the distributed sketching estimator converges to $1$ as $q$ or $m$ gets large. This is different from the results of \cite{mahoney2018averaging}, which do not imply convergence to $1$ as $q\rightarrow \infty$ due to an additional bias term. We note that the relative error of the uniform sampling sketch has a dependence on the row coherence unlike the randomized Hadamard and leverage score sampling sketch.

The main idea for the proof of Theorem \ref{thm:errorbounds_other_sketches} involves using the decomposition given in Lemma \ref{expected_obj_val_diff} along with the upper bounds for $\Exs[\|U^TS^TSb^\perp\|_2^2|E]$. Then, we use Markov's inequality to find a high probability bound for the approximation quality of the averaged estimator, conditioned on the events $E_k$, $k=1,\dots,q$ defined as $(1-\epsilon)I_d \preceq (U^TS_k^TS_kU)^{-1} \preceq (1+\epsilon)I_d$. This inequality is a direct result of the subspace embedding property and holds with high probability. The last component of the proof deals with removing the conditioning on the events $E_k$ to find the upper bound for the relative error $f(\bar{x})/f(x^*)$.

\section{Distributed Sketching for Nonlinear Convex Optimization Problems} \label{sec:nonlinear_problems}

In this section, we consider randomized second order methods for solving a broader range of problems by leveraging a distributed version of the Newton Sketch algorithm described in \cite{pilanci2017newton}.

\subsection{Distributed Newton Sketch} \label{subsec:newton_sketch}
The update equation for Newton's method is of the form $x_{t+1} = x_t - \alpha_1 H_t^{-1} g_t$, where $H_t \in \mathbb{R}^{d\times d}$ and $g_t \in \mathbb{R}^{d}$ denote the Hessian matrix and the gradient vector at iteration $t$ respectively, and $\alpha_1$ is the step size. In contrast, Newton Sketch performs the approximate updates
\begin{align}
    x_{t+1} = x_t + \alpha_1 \arg\min_{\Delta}(\frac{1}{2}\|S_t H_t^{1/2}\Delta \|_2^2 + g_t^T\Delta),
\end{align}
where the sketching matrices $S_t \in \mathbb{R}^{m\times n}$ are refreshed every iteration. There can be a multitude of options for devising a \textit{distributed} Newton's method or a \textit{distributed} Newton sketch algorithm. Here we consider a scheme that is similar in spirit to the GIANT algorithm \cite{mahoney2018giant} where worker nodes communicate length-$d$ approximate update directions to be averaged at the master node. Another alternative scheme would be to communicate the approximate Hessian matrices, which would require an increased communication load of $d^2$ numbers.

We consider Hessian matrices of the form $H_t=(H_t^{1/2})^T H_t^{1/2}$, where we assume that $H_t^{1/2} \in \mathbb{R}^{n\times d}$ is a full rank matrix and $n \geq d$. Note that this factorization is already available in terms of scaled data matrices in many problems as we illustrate in the sequel. This enables the fast construction of an approximation of $H_t$ by applying sketching $S_tH_t^{1/2}$ which leads to the approximation $\hat{H}_t = (S_tH_t^{1/2})^T S_tH_t^{1/2}$. Averaging for regularized problems with Hessian matrices of the form $H_t=(H_t^{1/2})^T H_t^{1/2} + \lambda_1 I_d$ will be considered in the next subsection.

The update equation for \textit{distributed} Newton sketch for a system with $q$ worker nodes can be written as
\begin{align}
    x_{t+1} = x_t + \alpha_2 \frac{1}{q} \sum_{k=1}^q \arg\min_{\Delta}\,\left(\frac{1}{2}\|S_{t,k} H_t^{1/2}\Delta \|_2^2 + g_t^T\Delta \right).
\end{align}
Note that the above update requires access to the full gradient $g_t$. If worker nodes do not have access to the entire dataset, then this requires an additional communication round per iteration where worker nodes communicate their local gradients with the master node, which computes the full gradient and broadcasts to worker nodes. The details of the distributed Newton sketch method are given in Algorithm \ref{alg:dist_newton_sketch}. 

\begin{figure}[!t]
 \removelatexerror
  \begin{algorithm}[H]
  \caption{Distributed Newton Sketch}
  \label{alg:dist_newton_sketch}
  \begin{algorithmic}
  \STATE{\textbf{Input:} Tolerance $\epsilon$}
  %\STATE{\bfseries repeat}
  \WHILE{$g_t^T \left(\sum_{k=1}^q \hat{\Delta}_{t,k}\right)/2 \leq \epsilon$}
  \FOR{worker $k=1,\dots,q$ in parallel}
  \STATE{Obtain $S_{t,k}H_t^{1/2}$.}
  \STATE{Obtain the gradient $g_t$.}
  \STATE{Compute approximate Newton direction $\hat{\Delta}_{t,k} = \arg\min_{\Delta}(\frac{1}{2}\|S_{t,k} H_t^{1/2}\Delta \|_2^2 + g_t^T\Delta)$ and send to master node.}
  \ENDFOR
  \STATE{\textbf{Master node:} Determine step size $\alpha_2$ and update $x_{t+1} = x_t + \alpha_2 \frac{1}{q} \sum_{k=1}^q \hat{\Delta}_{t,k}$.}
  \ENDWHILE
  %\STATE{\bfseries until} {$g_t^T \left(\sum_{k=1}^q \hat{\Delta}_{t,k}\right)/2 \geq \epsilon$ is satisfied}
  \end{algorithmic}
\end{algorithm}
\end{figure}

% \begin{algorithm}
% \caption{Distributed Randomized Regression}
% \label{distrib_avg_alg}
% \begin{algorithmic}
%  \STATE{\textbf{Input:} Data matrix $A \in \mathbb{R}^{n\times d}$, target vector $b \in \mathbb{R}^{n}$.}
%  \FOR{worker $k=1,\dots,q$ in parallel}
%  \STATE{Obtain the sketched data and sketched output: $S_kA$ and $S_kb$.}
%  \STATE{Solve $\hat x_k = \arg\min_{x} \|S_kAx-S_kb\|^2_2$ and send $\hat x_k$ to the master node.}
%  \ENDFOR
%  \STATE{\textbf{Master node:} return $\bar x = \frac{1}{q} \sum_{k=1}^q \hat x_k$.}
% \end{algorithmic}
% \end{algorithm}

% \subsubsection{Gaussian Sketch}
We analyze the bias and the variance of the update directions for distributed Newton sketch, and derive exact expressions for Gaussian sketching matrices. First, we establish the notation for update directions.
We will let $\Delta_t^*$ denote the optimal Newton update direction at iteration $t$
\begin{align}
    \Delta_t^* = ((H_t^{1/2})^T H_t^{1/2})^{-1}g_t
\end{align}
and let $\hat{\Delta}_{t,k}$ denote the approximate update direction returned by worker node $k$ at iteration $t$, which has the closed form expression
\begin{align}
    \hat{\Delta}_{t,k} = \alpha_s \left((H_t^{1/2})^TS_{t,k}^TS_{t,k}H_t^{1/2}\right)^{-1}g_t.
\end{align}
Note that the step size for the averaged update direction will be calculated as $\alpha_2 = \alpha_1 \alpha_s$. Theorem \ref{thm:dist_newton_sketch} characterizes how the update directions need to be scaled to obtain an unbiased update direction, and also a minimum variance estimator.

\begin{theorem} \label{thm:dist_newton_sketch}
For Gaussian sketches $S_{t,k}$, assuming $H_t^{1/2}$ is full column rank, the variance $\Exs[\|H_t^{1/2}(\hat{\Delta}_{t,k} - \Delta_t^*)\|_2^2]$ is minimized when $\alpha_s$ is chosen as $\alpha_s =  \frac{\theta_1}{\theta_2}$ whereas the bias $\Exs[H_t^{1/2}(\hat{\Delta}_{t,k} - \Delta_t^*)]$ is zero when $\alpha_s =  \frac{1}{\theta_1}$, where $\theta_1$ and $\theta_2$ are as defined in \eqref{eq:theta_1_2_definitions}.
\end{theorem}

\begin{figure}%[!t]
\hspace{2cm}
\begin{minipage}[b]{0.3\linewidth}
  \centering
  \centerline{\includegraphics[width=\columnwidth]{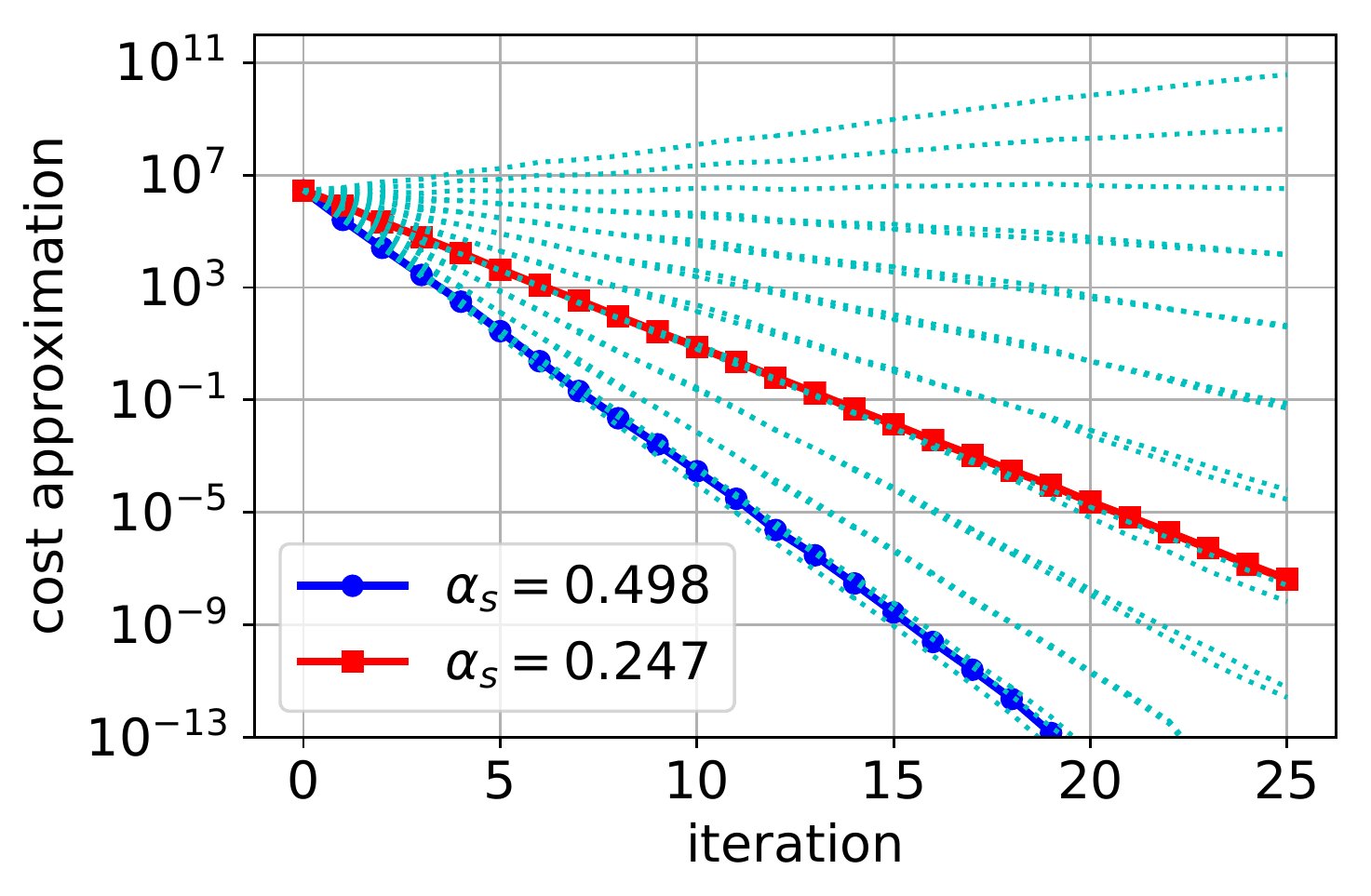}}
  \centerline{(a) $q=10$ workers}\medskip
\end{minipage}
\hspace{2cm}
\begin{minipage}[b]{0.3\linewidth}
  \centering
  \centerline{\includegraphics[width=\columnwidth]{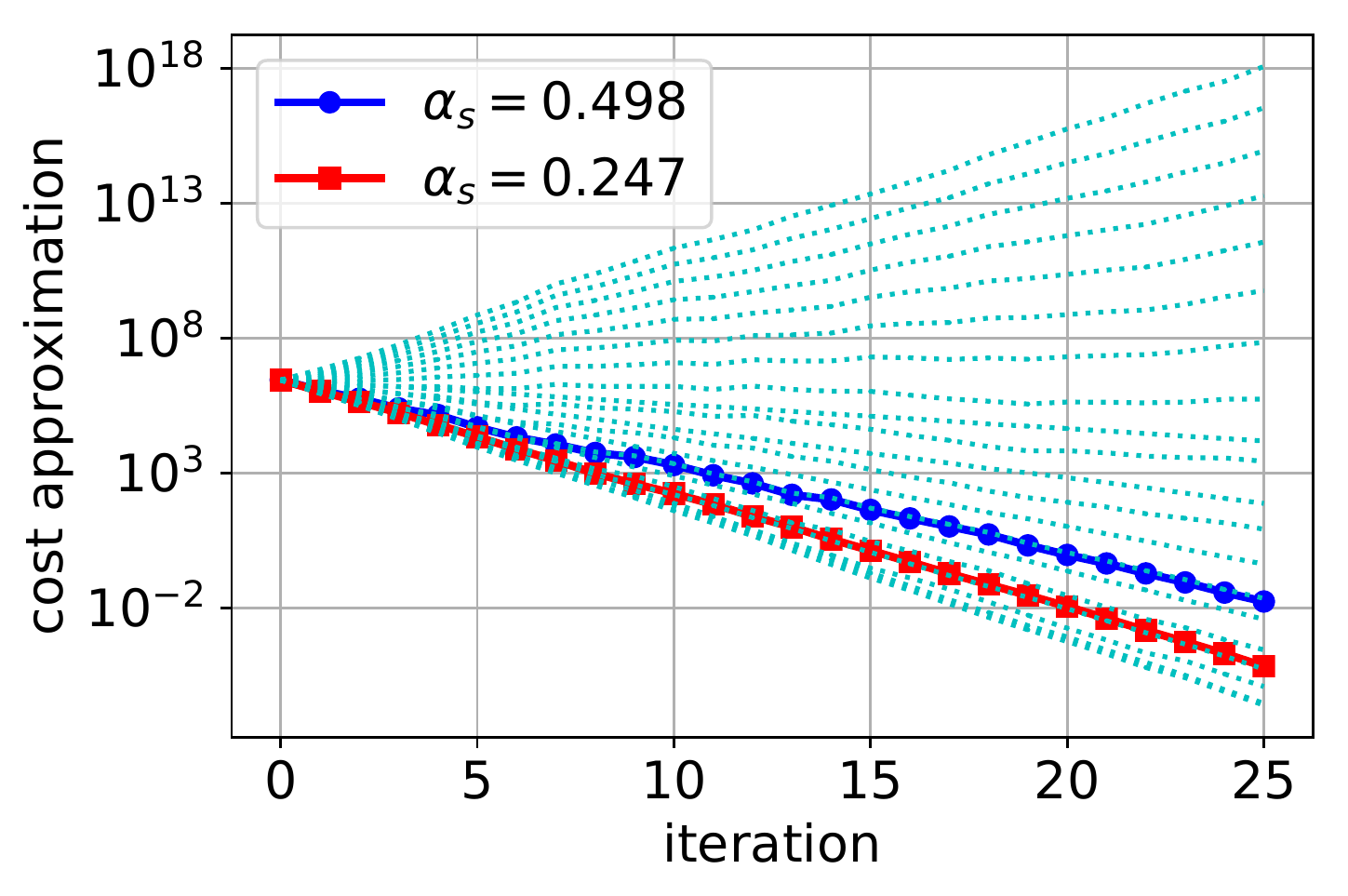}}
  \centerline{(b) $q=2$ workers}\medskip
\end{minipage}
\caption{Cost approximation $(f(x_t)-f(x^*))/f(x^*)$ for Algorithm \ref{alg:dist_newton_sketch} against iteration number $t$ for various step sizes in solving a linear least squares problem on randomly generated data. The cyan colored dotted lines show cost approximation when we make a search for the learning rate $\alpha_s$ between $0.05$ and $1$. The blue line with circle markers corresponds to $\alpha_s=1/\theta_1$ that leads to the unbiased estimator and the red line with square markers corresponds to $\alpha_s=\theta_1/\theta_2$ that gives the minimum variance. The resulting step size scaling factors $\alpha_s$ are shown in the legends of the plots. The parameters used in these experiments are $n=1000$, $d=200$, $m=400$.}
\label{fig:newton_sketch_data_gauss4th}
\end{figure}

Figure \ref{fig:newton_sketch_data_gauss4th} demonstrates that choosing $\alpha_2=\alpha_1 \alpha_s$ when $\alpha_s$ is calculated using the unbiased estimator formula $\alpha_s=1/\theta_1$ leads to faster decrease of the objective value when the number of workers is high. If the number of workers is small, one should choose the step size that minimizes variance using $\alpha_s=\theta_1/\theta_2$ instead. Furthermore, Figure \ref{fig:newton_sketch_data_gauss4th}(a) illustrates that the blue curve with square markers is in fact the best one could hope to achieve as it is very close to the best cyan dotted line.

% \textit{Non-identical sketch sizes:} Theorem \ref{thm:dist_newton_sketch} establishes that whenever the sketch dimension varies among workers, it is possible to obtain an unbiased update direction by computing $\alpha_s$ for every worker individually.

%%%%%%%%%%%%%%%%%%%%%%%%%%%%%%%

\subsection{Distributed Newton Sketch for Regularized Problems} \label{subsec:regularized_newton}
We now consider problems with squared $\ell_2$-norm regularization. In particular, we study problems with Hessian matrices of the form $H_t=(H_t^{1/2})^TH_t^{1/2} + \lambda_1 I_d$. Sketching can be applied to obtain approximate Hessian matrices as $H_t=(S_tH_t^{1/2})^TS_tH_t^{1/2} + \lambda_2 I_d$. Note that the case $\lambda_2=\lambda_1$ corresponds to the setting in the GIANT algorithm described in \cite{mahoney2018giant}. 

Theorem \ref{opt_lambda_2_newton} establishes that $\lambda_2$ should be chosen according to the formula \eqref{opt_lambda_2_formula} under the assumption that the singular values of $H_t^{1/2}$ are identical. We later verify empirically that when the singular values are not identical, plugging the mean of the singular values into the formula still leads to improvements over the case of $\lambda_2 = \lambda_1$.

\begin{theorem} \label{opt_lambda_2_newton}
Given the thin SVD decomposition $H_t^{1/2}=U\Sigma V^T \in \mathbb{R}^{n\times d}$ and $n \geq d$ where $H_t^{1/2}$ is assumed to have full rank and satisfy $\Sigma = \sigma I_d$, the bias of the single sketch Newton step estimator $\Exs[ H_t^{1/2}(\hat{\Delta}_{t,k} - \Delta_t^*)]$ is equal to zero as $n$ goes to infinity when $\lambda_2$ is chosen as
\begin{align} \label{opt_lambda_2_formula}
    \lambda_2^* = \left(\lambda_1 + \frac{d}{m}\sigma^2 \right) \left(1 - \frac{d/m}{1 + \lambda_1 \sigma^{-2} + d/m} \right),
\end{align}
where $\Delta_t^* = ((H_t^{1/2})^TH_t^{1/2} + \lambda_1 I_d)^{-1} g_t$ and $\hat{\Delta}_{t,k} = ((S_{t,k}H_t^{1/2})^TS_{t,k}H_t^{1/2} + \lambda_2 I_d)^{-1} g_t$, and $S_{t,k}$ is the Gaussian sketch.
\end{theorem}

\begin{remark}
The proof of Theorem \ref{opt_lambda_2_newton} builds on the proof of Theorem \ref{thm:opt_lambda_2_newton_LS} (see Appendix). The main difference between these results is that the setting in Theorem \ref{opt_lambda_2_newton} is such that we assume access to the full gradient, i.e. only the Hessian matrix is sketched. In contrast, the setting in Theorem \ref{thm:opt_lambda_2_newton_LS} does not use the exact gradient.
\end{remark}

\section{Applications of Distributed Sketching} \label{sec:example_problems}
This section describes some example problems where our methodology can be applied. In particular, the problems in this section are convex problems that are efficiently addressed by our methods in distributed systems. We present numerical results on these problems in Section \ref{sec:numerical_results}.

\subsection{Logistic Regression} 
We begin by considering the well-known logistic regression model with squared $\ell_2$-norm penalty. The optimization problem for this model can be formulated as $\textrm{minimize}_{x} f(x)$ where
\begin{align}
     f(x) = -\sum_{i=1}^n \left( y_i \log(p_i) + (1-y_i)\log(1-p_i)  \right) + \frac{\lambda_1}{2} \|x\|_2^2,
\end{align}
and $p \in \mathbb{R}^{n}$ is defined such that $p_i = 1/(1+\exp(-\tilde{a}_i^Tx))$. $\tilde{a}_i \in \mathbb{R}^d$ denotes the $i$'th row of the data matrix $A \in \mathbb{R}^{n\times d}$. The output vector is denoted by $y \in \mathbb{R}^n$.

We can use the distributed Newton sketch algorithm to solve this optimization problem. The gradient and Hessian for $f(x)$ are as follows
\begin{align*}
    g &= A^T(p-y) + \lambda_1 x, \\
    H &= A^TDA + \lambda_1 I_d \,.
\end{align*}
$D$ is a diagonal matrix with the elements of the vector $p(1-p)$ as its diagonal entries. The sketched Hessian matrix in this case can be formed as
$(SD^{1/2}A)^T(SD^{1/2}A) + \lambda_2^* I_d$ where $\lambda_2^*$ can be calculated using \eqref{opt_lambda_2_formula}, and we can set $\sigma$ to the mean of the singular values of $D^{1/2}A$. Since the entries of the matrix $D$ are a function of the variable $x$, the matrix $D$ gets updated every iteration. It might be computationally expensive to re-compute the mean of the singular values of $D^{1/2}A$ in every iteration. However, we have found through experiments that it is not required to compute the exact value of the mean of the singular values for bias reduction. For instance, setting $\sigma$ to the mean of the diagonals of the matrix $D^{1/2}$ as a heuristic works sufficiently well.

\subsection{Inequality Constrained Optimization} 
The second application that we consider is the inequality constrained optimization problem of the form
\begin{align}
    \textrm{minimize}_{x} \quad & \|x-c\|_2^2 \nonumber \\
    \textrm{subject to} \quad & \|Ax\|_{\infty} \leq \lambda %\mathbf{1}_{n\times 1}
\end{align}
where $A \in \mathbb{R}^{n \times d}$, and $c \in \mathbb{R}^{d}$ are the problem data, and $\lambda \in \mathbb{R}$ is a positive scalar. Note that this problem is the dual of the Lasso problem given by $\text{min}_x$ $\lambda \|x\|_1 + \frac{1}{2}\|Ax-c\|_2^2$.

The above problem can be tackled by the standard log-barrier method \cite{boyd2004convex}, by solving sequences of unconstrained barrier penalized problems as follows
\begin{align} \label{log_barrier_opt_prob_2}
    \textrm{minimize}_x \, - \sum_{i=1}^n \log(&-\tilde{a}_i^Tx+\lambda) -\sum_{i=1}^n \log(\tilde{a}_i^Tx+\lambda) + \lambda_1 \|x\|_2^2 - 2\lambda_1 c^Tx + \lambda_1 \|c\|_2^2
\end{align}
where $\tilde{a}_i$ represents the $i$'th row of $A$. The distributed Newton sketch algorithm could be used to solve this problem.
The gradient and Hessian of the objective are given by
\begin{align*}
    g &= -A_c^T D \mathbf{1}_{2n\times 1} + 2\lambda_1 x - 2\lambda_1 c, \\
    H &= (DA_c)^T(DA_c) + 2\lambda_1 I_d.
\end{align*}
Here $A_c = [A^T, -A^T]^T$ and $D$ is a diagonal matrix with the element-wise inverse of the vector $(A_cx-\mathbf{1}_{2n\times 1})$ as its diagonal entries. $\mathbf{1}_{2n\times 1}$ is the length-$2n$ vector of all $1$'s. The sketched Hessian can be written in the form of $(SDA_c)^T(SDA_c) + \lambda_2I_d$.

\begin{remark}
Since the contribution of the regularization term in the Hessian matrix is $2\lambda_1 I_d$, we need to plug $2\lambda_1$ instead of $\lambda_1$ in the formula for computing $\lambda_2^*$.
\end{remark}

\subsection{Fine Tuning of Pre-Trained Neural Networks}
The third application is neural network fine tuning. Let $f_{NN}(\cdot): \mathbb{R}^d \rightarrow \mathbb{R}^r$ represent the output of the $(L-1)$'st layer of a pre-trained network which consists of $L$ layers. We are interested in re-using a pre-trained neural network for, say, a task on a different dataset by learning a different final layer. More precisely, assuming squared loss, one way that we can formulate this problem is as follows
\begin{align} \label{eq:tuning_NN_least_sq}
    W_L^* = \arg\min_{W_L} \Bigg\Vert  \underbrace{\begin{bmatrix}f_{NN}(\tilde{a}_1)^T \\ \vdots \\ f_{NN}(\tilde{a}_n)^T \end{bmatrix}}_{f_{NN}(A)} W_L - b \Bigg\Vert _F^2 + \lambda_1 \|W_L\|_F^2 \,,
\end{align}
where $W_L \in \mathbb{R}^{r\times C}$ is the weight matrix of the final layer. The vectors $\tilde{a}_i \in \mathbb{R}^d$ $i=1,\dots,n$ denote the rows of the data matrix $A \in \mathbb{R}^{n\times d}$ and $b \in \mathbb{R}^{n\times C}$ is the output matrix. In a classification problem, $C$ would correspond to the number of classes.

For large scale datasets, it is important to develop efficient methods for solving the problem in \eqref{eq:tuning_NN_least_sq}. Distributed randomized ridge regression algorithm (listed in Algorithm \ref{alg:dist_regularized_LS}) could be applied to this problem where the sketched sub-problems are of the form
\begin{align}
    \hat{W}_{L,k} = \arg\min_{W_L} \Vert S_kf_{NN}(A) W_L - S_kb \Vert _F^2 + \lambda_2^* \|W_L\|_F^2
\end{align}
and the averaged solution can be computed at the master node as $\bar{W}_L = \frac{1}{q} \sum_{k=1}^q \hat{W}_{L,k}$. This leads to an efficient non-iterative asynchronous method for fine-tuning of a neural network. 
\section{Numerical Results} \label{sec:numerical_results}

We have implemented the distributed sketching methods for AWS Lambda in Python using the Pywren package \cite{jonas2017pywren}. The setting for the experiments is a centralized computing model where a single master node collects and averages the outputs of the $q$ worker nodes. The majority of the figures in this section plots the approximation error which we define as $(f(\bar{x})-f(x^*)) / f(x^*)$.

\subsection{Hybrid Sketch}
In a distributed computing setting, the amount of data that can be fit into the memory of nodes and the size of the largest problem that can be solved by the nodes often do not match. The hybrid sketch idea is motivated by this mismatch and it is basically a sequentially concatenated sketching scheme where we first perform uniform sampling with dimension $m^\prime$ and then sketch the sampled data using another sketching method, preferably with better convergence properties (say, Gaussian) with dimension $m$. Worker nodes load $m^\prime$ rows of the data matrix $A$ into their memory and then perform sketching, reducing the number of rows from $m^\prime$ to $m$. In addition, we note that if $m^\prime = m$, hybrid sketch reduces to sampling and if $m^\prime = n$, then it reduces to Gaussian sketch. 
The hybrid sketching scheme does not take privacy into account as worker nodes are assumed to have access to the data matrix.

% We present experiment results in the numerical results section showing the practicality of the hybrid sketch idea. 
For the experiments involving very large scale datasets, we have used Sparse Johnson-Lindenstrauss Transform (SJLT) \cite{nelson2013osnap} as the second sketching method in the hybrid sketch due to its low computational complexity.

\subsection{Airline Dataset}
We have conducted experiments with the publicly available Airline dataset \cite{airline_dataset}. This dataset contains information on domestic USA flights between the years 1987-2008. We are interested in predicting whether there is going to be a departure delay or not, based on information about the flights. %The dataset contains information for around $120$ million flights with each flight having $29$ attributes.
More precisely, we are interested in predicting whether \texttt{DepDelay} $>15$ minutes using the attributes \texttt{Month}, \texttt{DayofMonth}, \texttt{DayofWeek}, \texttt{CRSDepTime}, \texttt{CRSElapsedTime}, \texttt{Dest}, \texttt{Origin}, and \texttt{Distance}. Most of these attributes are categorical and we have used dummy coding to convert these categorical attributes into binary representations. The size of the input matrix $A$, after converting categorical features into binary representations, becomes $(1.21 \times 10^8) \times 774$.%, and the output vector $b$ is a $(1.21 \times 10^8)$-dimensional vector.

%The aim of the experiments presented in this subsection is not to address the predictive power of linear least squares in predicting flight delays. Instead, 

%The goal of the experiments presented in this subsection is to illustrate the ability of distributed sketching to solve large-scale linear least squares problems. Even though solving a linear least squares problem by itself does not lead to high accuracy predictions, many other approaches to this problem requires solving least squares problems as a step. 
%Hence, it is of significance to have efficient tools for solving large-scale linear least squares problems with convergence guarantees.

\begin{figure}[!t]
\begin{minipage}[b]{0.24\linewidth}
  \centering
  \centerline{\includegraphics[width=\columnwidth]{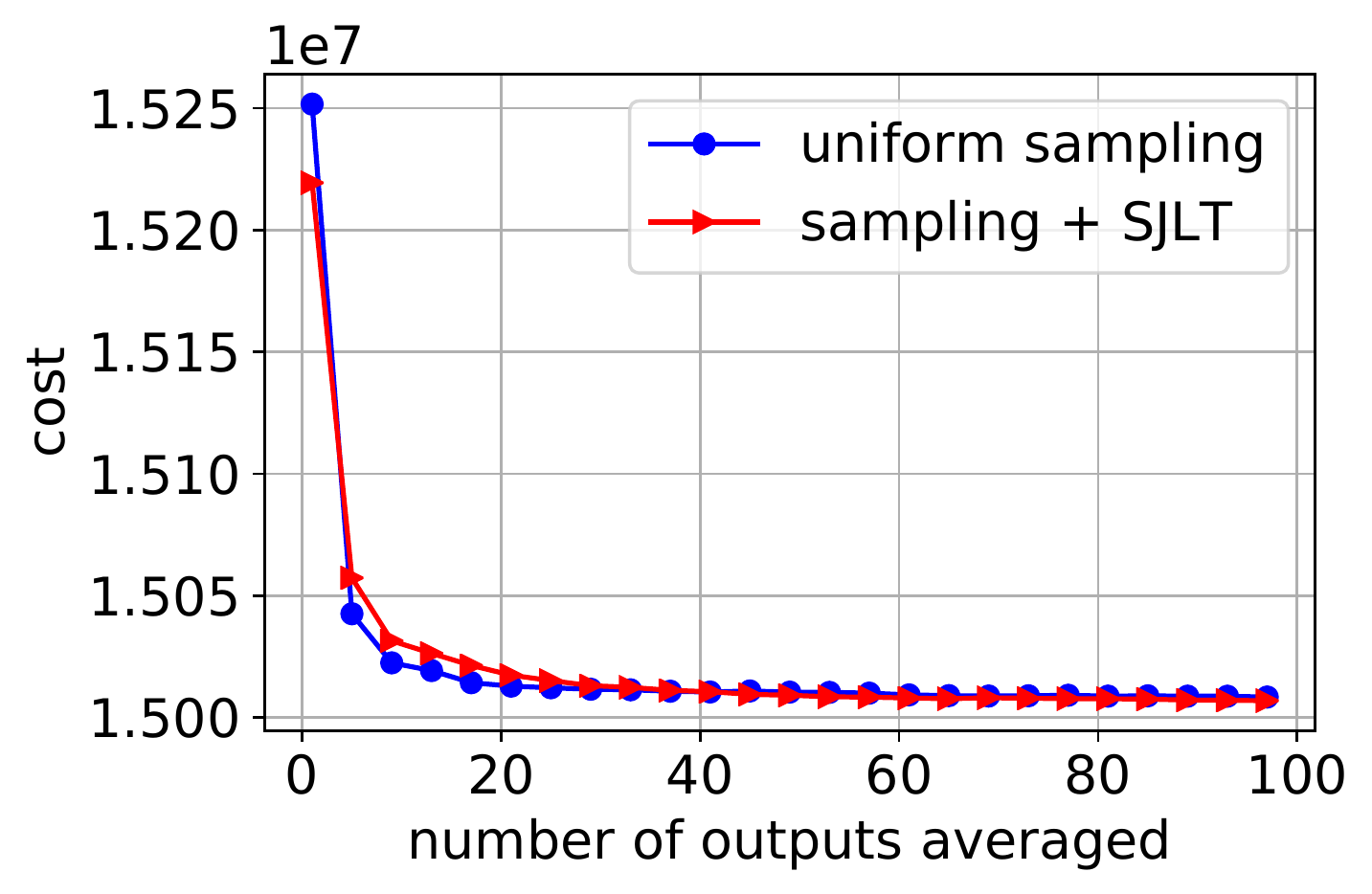}}
  \centerline{a) $m^\prime = 5m=5\times 10^5$}\medskip
\end{minipage}
\hfill
\begin{minipage}[b]{0.24\linewidth}
  \centering
  \centerline{\includegraphics[width=\columnwidth]{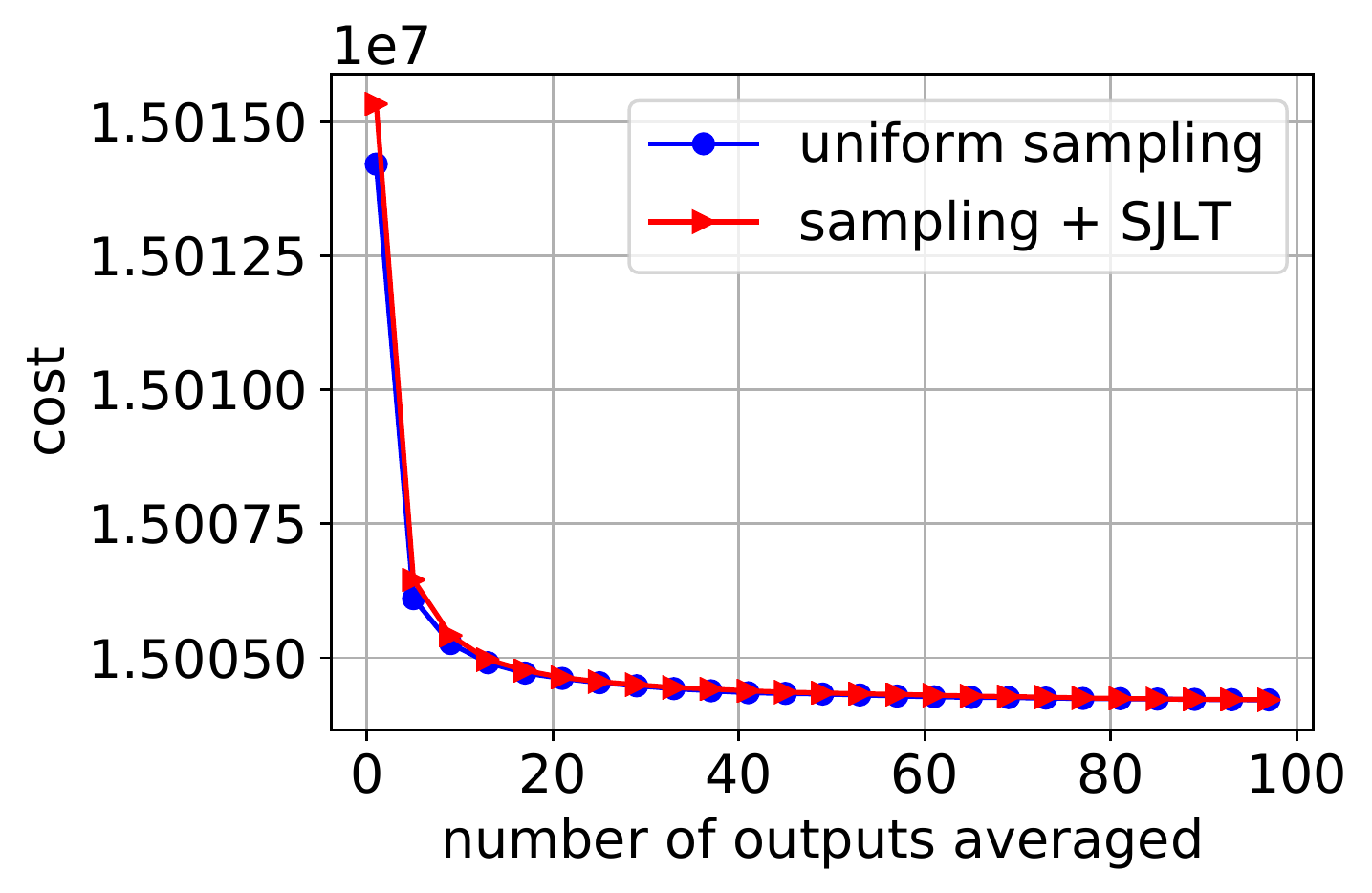}}
  \centerline{b) $m^\prime=2m=2\times 10^6$}\medskip
\end{minipage}
\hfill
\begin{minipage}[b]{0.24\linewidth}
  \centering
  \centerline{\includegraphics[width=\columnwidth]{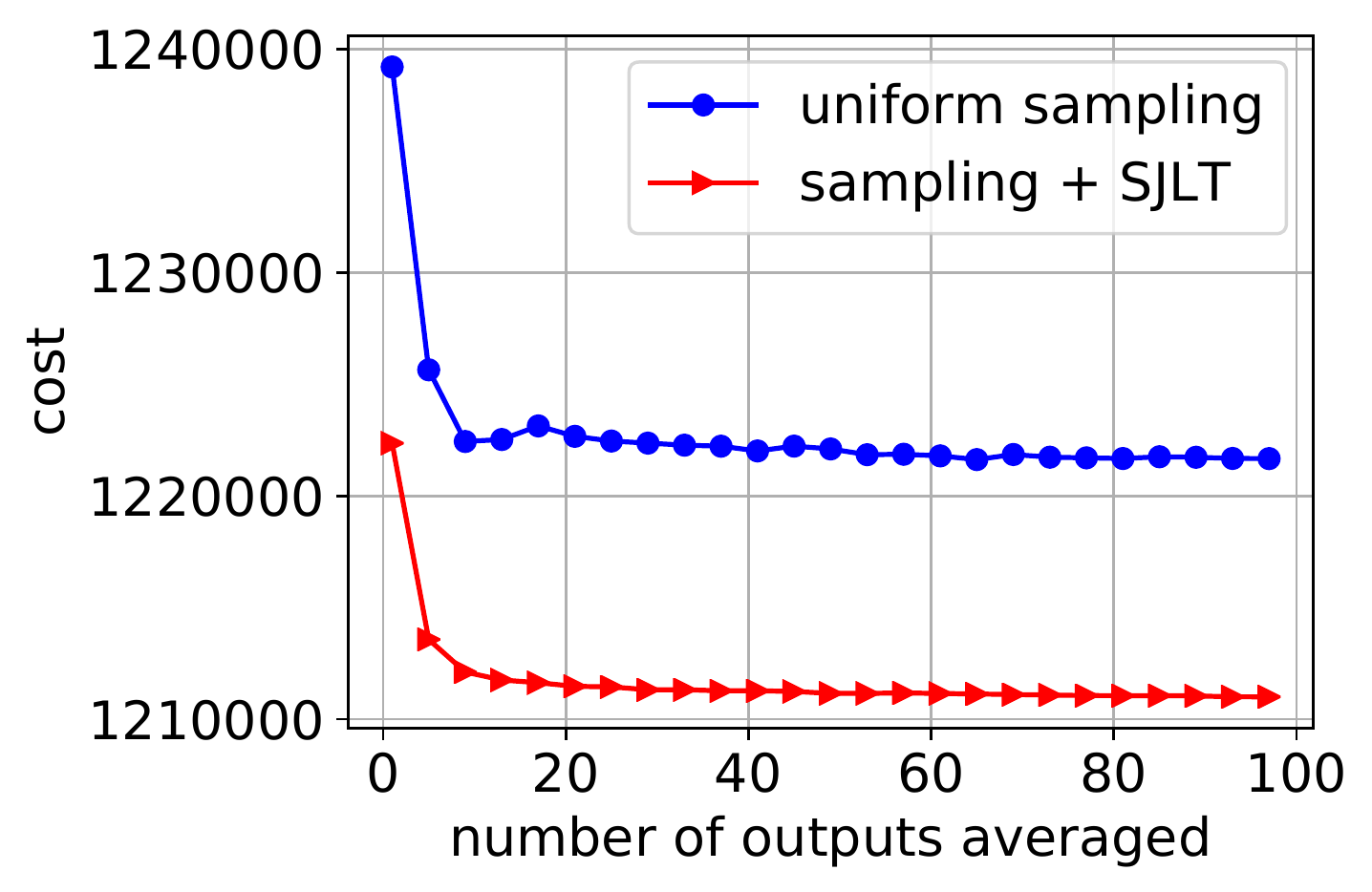}}
  \centerline{c) $m^\prime = 5m=5\times 10^5$}\medskip
\end{minipage}
\hfill
\begin{minipage}[b]{0.24\linewidth}
  \centering
  \centerline{\includegraphics[width=\columnwidth]{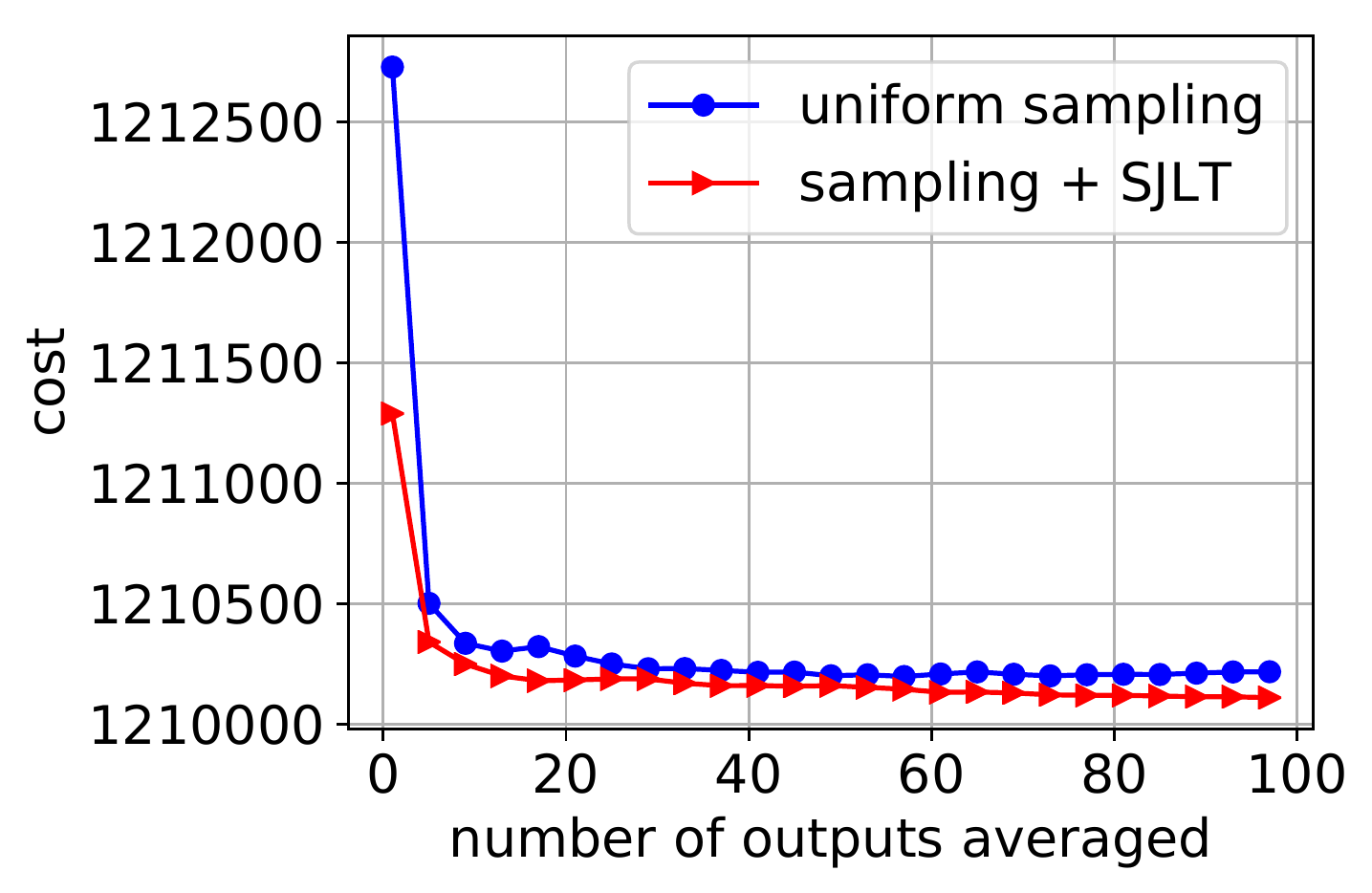}}
  \centerline{d) $m^\prime = 2m=2\times 10^6$}\medskip
\end{minipage}
\caption{AWS Lambda experiments on the entire airline dataset ($n=1.21\times 10^8$) with $q=100$ workers. %The vertical axis represents the approximation error, and the horizontal axis shows the number of outputs used in obtaining the averaged estimate. 
The run times for each plot are as follows (given in this order: sampling, sampling followed by SJLT): a: 37.5, 43.9 seconds, b: 48.3, 60.1 seconds, c: 39.8, 52.9 seconds, d: 78.8, 107.6 seconds.}
\label{airline_cloud}
\end{figure}

We have solved the linear least squares problem on the entire airline dataset: $\text{minimize}_x \|Ax-b\|_2^2$ using $q$ workers on AWS Lambda. The output $b$ for the plots a and b in Figure \ref{airline_cloud} is a vector of binary variables indicating delay. The output $b$ for the plots c and d in Figure \ref{airline_cloud} is artificially generated via $b=Ax_{truth}+ \epsilon$ where $x_{truth}$ is the underlying solution and $\epsilon$ is random Gaussian noise distributed as $\mathcal{N}(0, 0.01I)$. Figure \ref{airline_cloud} shows that sampling followed by SJLT leads to a lower error.

%Figure \ref{airline_cloud} shows the results for the entire dataset where we have used AWS Lambda functions in parallel as the computing platform. 
We note that it increases the convergence rate to choose $m$ and $m^\prime$ as large as the available resources allow. Based on the run times given in the caption of Figure \ref{airline_cloud}, we see that the run times are slightly worse if SJLT is involved. Decreasing $m^\prime$ will help reduce this processing time at the expense of error performance.

\subsection{Image Dataset: Extended MNIST}
The experiments of this subsection are performed on the image dataset EMNIST (extended MNIST) \cite{emnist_dataset}. We have used the "bymerge" split of EMNIST, which has 700K training and 115K test images. The dimensions of the images are $28\times 28$ and there are 47 classes in total (letters and digits). Some capital letter classes have been merged with small letter classes (like C and c), and thus there are 47 classes and not 62.

\begin{figure}%[!t]
\hspace{2cm}
\begin{minipage}[b]{0.3\linewidth}
  \centering
  \centerline{\includegraphics[width=\columnwidth]{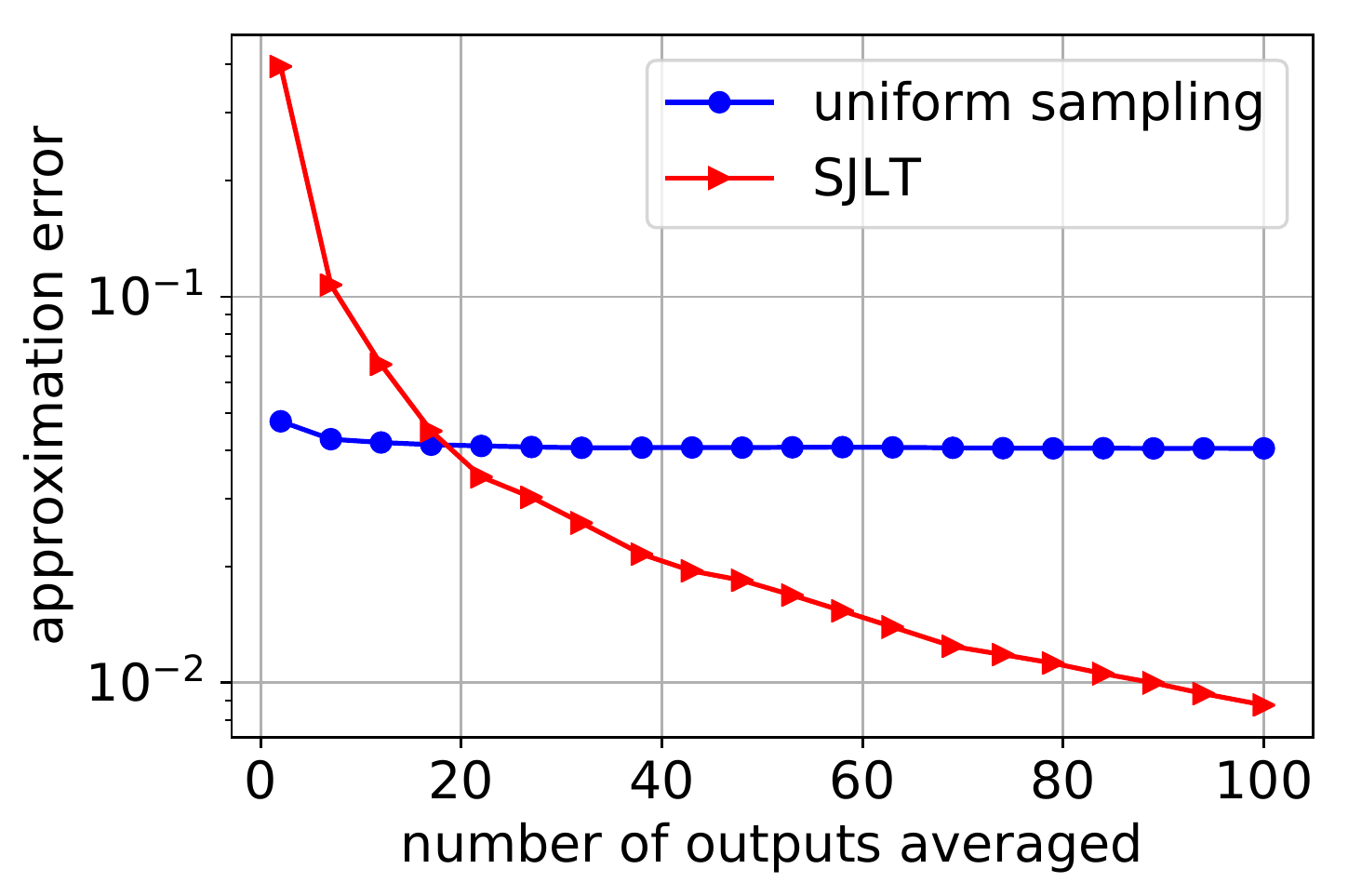}}
  \centerline{a) Approximation error}\medskip
\end{minipage}
\hspace{2cm}
\begin{minipage}[b]{0.3\linewidth}
  \centering
  \centerline{\includegraphics[width=\columnwidth]{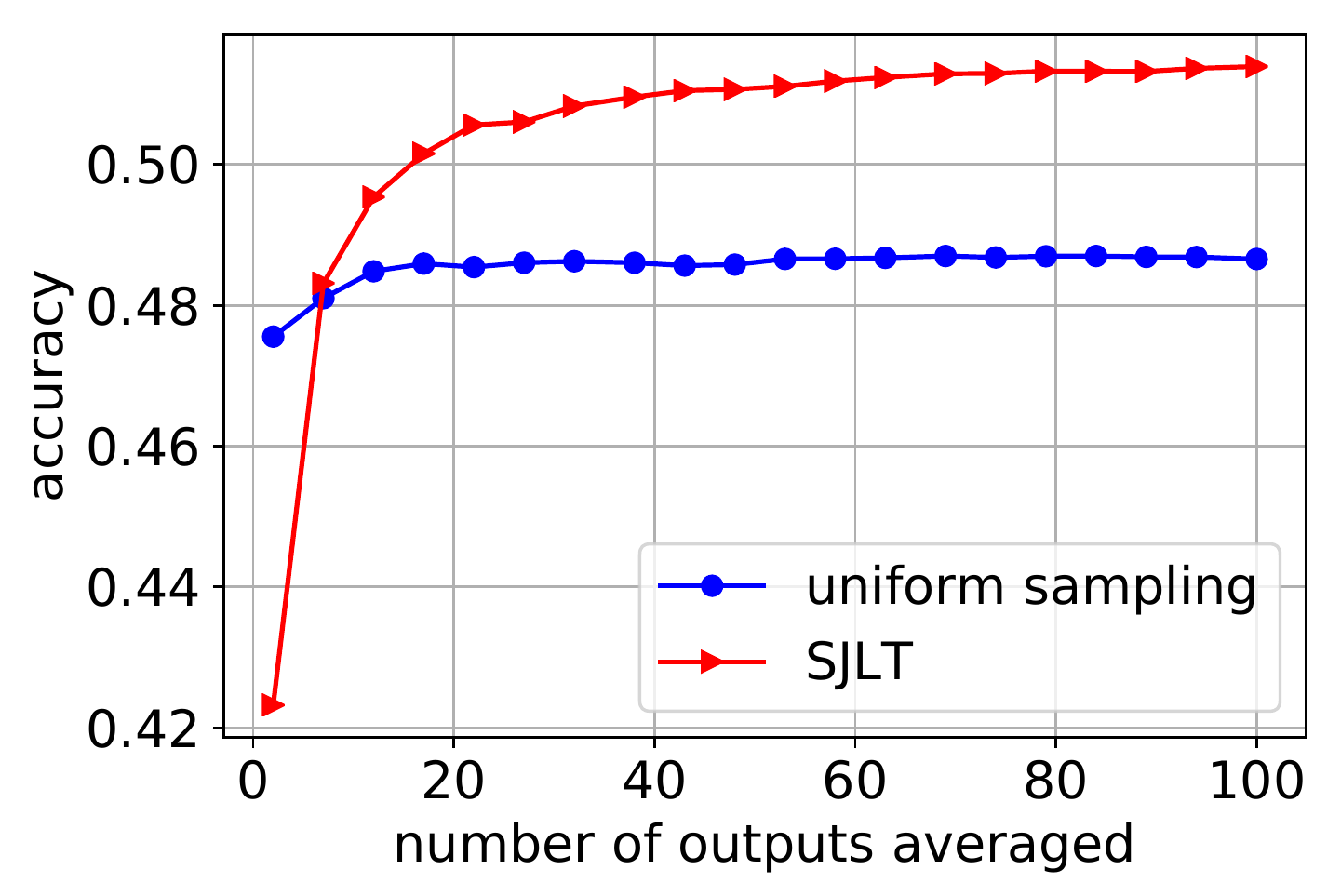}}
  \centerline{b) Test accuracy}\medskip
\end{minipage}
\caption{Approximation error and test set classification accuracy plots for the EMNIST-bymerge dataset where $q=100$, $m=2000$, $s=20$. The run times on AWS Lambda for uniform sampling and SJLT are $41.5$ and $66.9$ seconds, respectively.}
\label{emnist_experiments}
\end{figure}

Figure \ref{emnist_experiments} shows the approximation error and test accuracy plots when we fit a least squares model on the EMNIST-bymerge dataset using the distributed randomized regression algorithm. Because this is a multi-class classification problem, we have one-hot encoded the labels. %For these experiments, we have used $q=100$ workers in AWS Lambda and the sketch dimension is $m=2000$. 
Figure \ref{emnist_experiments} demonstrates that SJLT is able to drive the cost down more and leads to a better classification accuracy than uniform sampling.

\subsection{Performance on Large Scale Synthetic Datasets}
\begin{figure}%[!t]
\hspace{2cm}
\begin{minipage}[b]{0.3\linewidth}
  \centering
  \centerline{\includegraphics[width=\columnwidth]{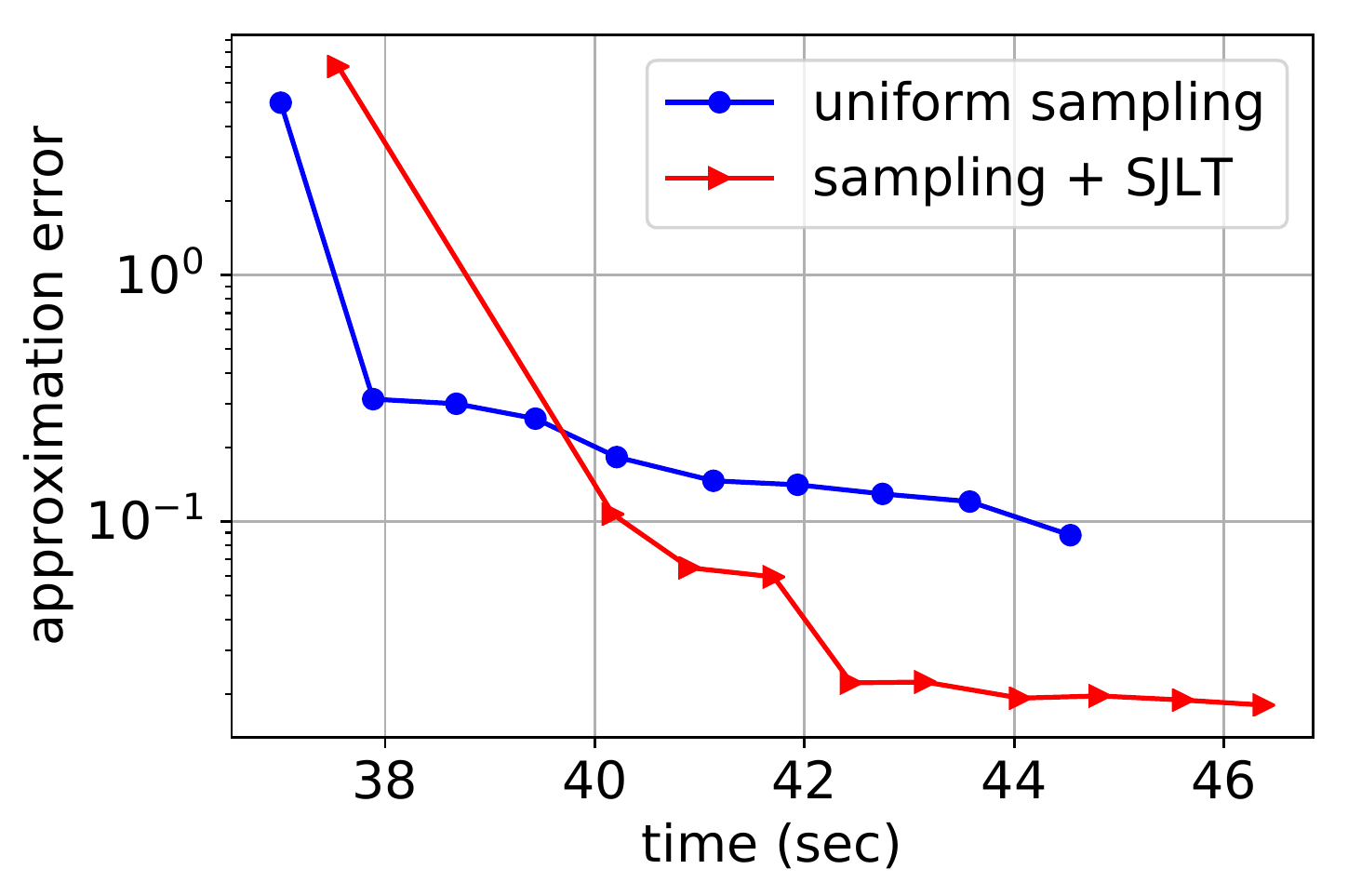}}
  \centerline{a) $A\in \mathbb{R}^{10^7 \times 10^3}$}\medskip
\end{minipage}
\hspace{2cm}
\begin{minipage}[b]{0.3\linewidth}
  \centering
  \centerline{\includegraphics[width=\columnwidth]{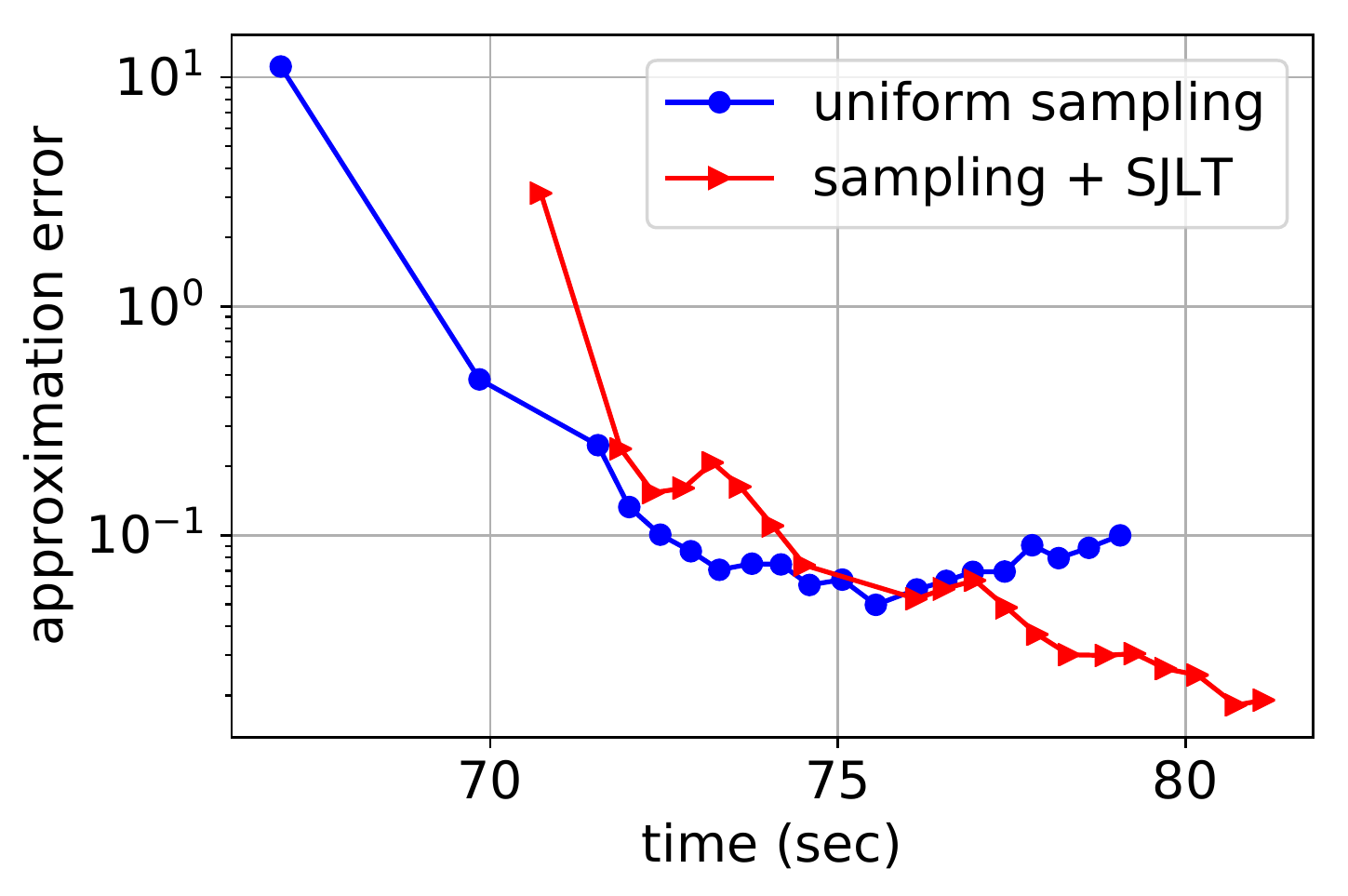}}
  \centerline{b) $A\in \mathbb{R}^{(2\times 10^7) \times (2\times 10^3)}$}\medskip
\end{minipage}
% \hfill
% \begin{minipage}[b]{0.48\linewidth}
%   \centering
%   \centerline{\includegraphics[width=4.5cm]{src/test_loss_1.pdf}}
% %  \vspace{2.0cm}
%   \centerline{(c) Test set 1}\medskip
% \end{minipage}
% \hfill
% \begin{minipage}[b]{0.48\linewidth}
%   \centering
%   \centerline{\includegraphics[width=4.5cm]{src/test_loss_2.pdf}}
% %  \vspace{2.0cm}
%   \centerline{(d) Test set 2}\medskip
% \end{minipage}
\caption{Approximation error vs time for AWS Lambda experiments on randomly generated large scale datasets ($q=200$ AWS Lambda functions have been used). %Plots a, b show the approximation errors on the training sets and plots c, d show the test set errors.
}
\label{random_data_experiments}
\end{figure}

We now present the results of the experiments carried out on randomly generated large scale data to illustrate scalability of the methods. Plots in Figure \ref{random_data_experiments} show the approximation error as a function of time, where the problem dimensions are as follows: $A\in \mathbb{R}^{10^7 \times 10^3}$ for plot a and $A\in \mathbb{R}^{(2\times 10^7) \times (2\times 10^3)}$ for plot b. These data matrices take up $75$ GB and $298$ GB, respectively. The data used in these experiments were randomly generated from the student's t-distribution with degrees of freedom of $1.5$ for plot a and $1.7$ for plot b. The output vector $b$ was computed according to $b = A x_{truth} + \epsilon$ where $\epsilon \in \mathbb{R}^n$ is i.i.d. noise distributed as $\mathcal{N}(0, 0.1 I_n)$. Other parameters used in the experiments are $m=10^4, m^\prime = 10^5$ for plot a, and $m=8\times 10^3, m^\prime = 8\times 10^4$ for plot b. %Furthermore, plots c and d in Figure \ref{random_data_experiments} show the approximation error on the test set. The test sets have also been randomly sampled from the same distribution as the training set and we have generated $10^6$ samples in both cases.
We observe that both plots in Figure \ref{random_data_experiments} reveal similar trends where the hybrid approach leads to a lower approximation error but takes longer due to the additional processing required for SJLT. %We also observe that the curve for the hybrid approach is below the curve for sampling, that is, the average due to the hybrid approach is a better approximate than the one due to sampling most of the time. 

\subsection{Numerical Results for the Right Sketch Method}
\begin{figure}%[!t]
\hspace{2cm}
\begin{minipage}[b]{0.3\linewidth}
  \centering
  \centerline{\includegraphics[width=\columnwidth]{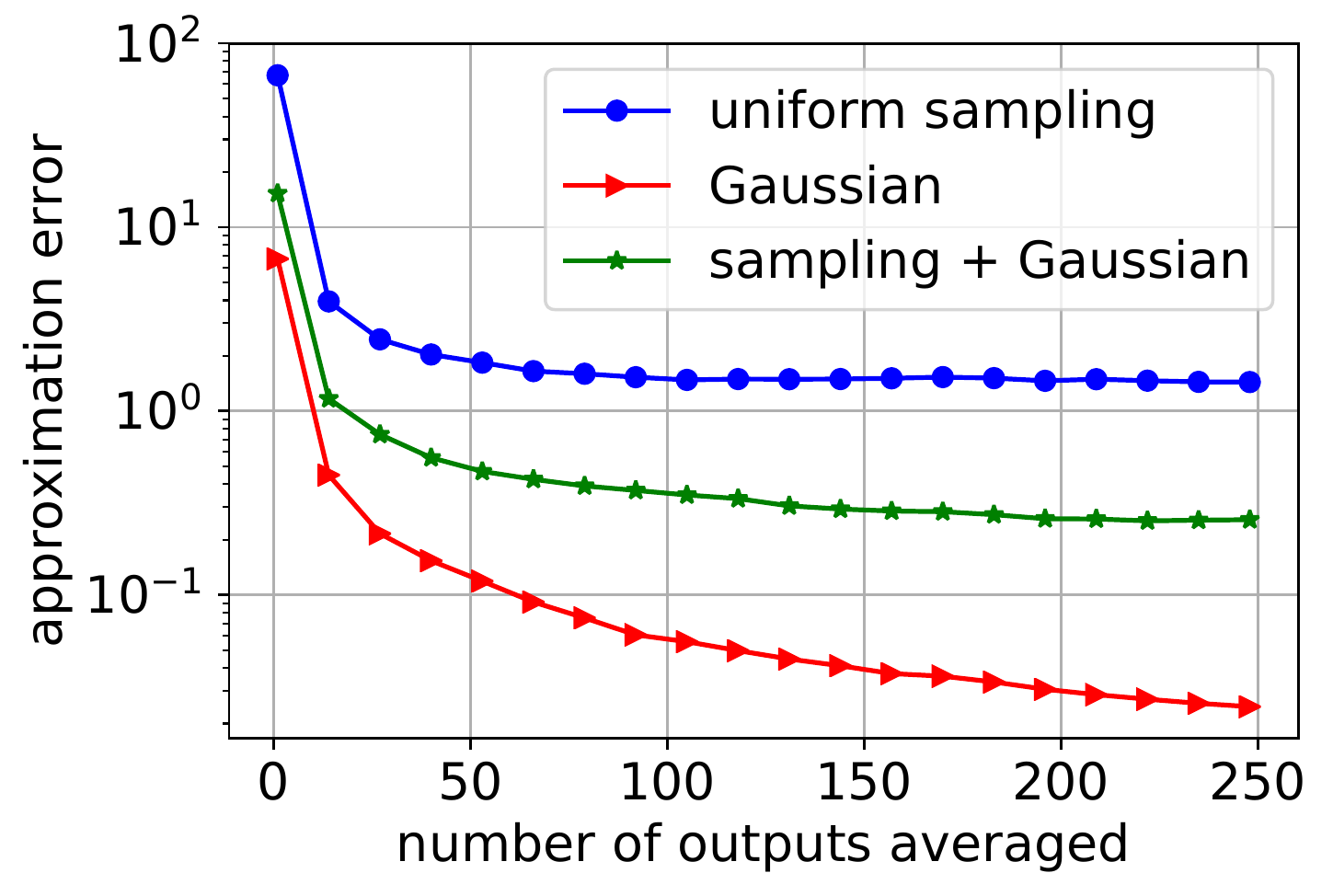}}
  \centerline{a) Random data}\medskip
\end{minipage}
\hspace{2cm}
\begin{minipage}[b]{0.3\linewidth}
  \centering
  \centerline{\includegraphics[width=\columnwidth]{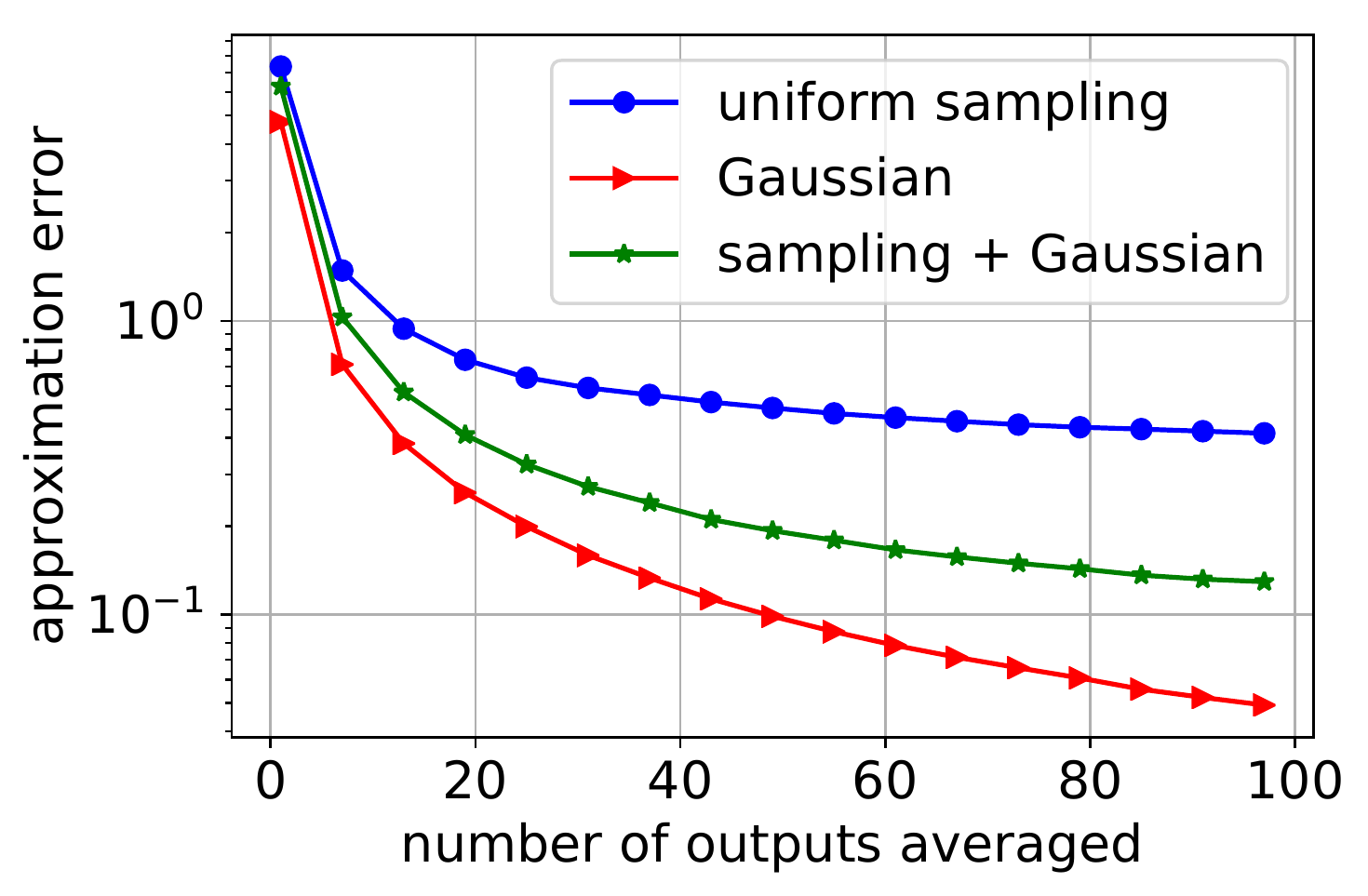}}
  \centerline{b) Airline data}\medskip
\end{minipage}
\caption{Averaging for least norm problems. \textit{Plot (a):} The parameters are $n=50$, $d=1000$, $m=200$, $m^\prime = 500$. \textit{Plot (b):} Least norm averaging applied to a subset of the airline dataset. The parameters are $n=2000$, $d=11588$, $m=4000$, $m^\prime =8000$. For this plot, the features include the pairwise interactions in addition to the original features.}
\label{right_sketch_plots}
\end{figure}

Figure \ref{right_sketch_plots} shows the approximation error as a function of the number of averaged outputs in solving the least norm problem for two different datasets. %The approximation error is the ratio of the $\ell_2$-norms of the averaged solution $\bar{x}$ and the optimal solution $x^*$, that is $\|\bar{x}\|_2^2 / \|x^*\|_2^2 - 1$.
The dataset for Figure \ref{right_sketch_plots}(a) is randomly generated with dimensions $A \in \mathbb{R}^{50\times 1000}$. We observe that Gaussian sketch outperforms uniform sampling in terms of the approximation error. Furthermore, Figure \ref{right_sketch_plots}(a) verifies that if we apply the hybrid approach of sampling first and then applying Gaussian sketch, then its performance falls between the extreme ends of only sampling and only using Gaussian sketch. Moreover, Figure \ref{right_sketch_plots}(b) shows the results for the same experiment on a subset of the airline dataset where we have included the pairwise interactions as features which makes this an underdetermined linear system. Originally, we had $774$ features for this dataset and if we include all $x_ix_j$ terms as features, we would have a total of $299925$ features, most of which are zero for all samples. We have excluded the all-zero columns from this extended matrix to obtain the final dimensions $2000\times 11588$. %This is to make sure that the sketched sub-problems are underdetermined as well. 

%%%%%%%%%%%%%%%%%%%%%%%%%%%%%%%%%%%%%%%%%%
%%%%%%%%%%%%%%%%%%%%%%%%%%%%%%%%%%%%%%%%%%
%%%%%%%%%%%%%%%%%%%%%%%%%%%%%%%%%%%%%%%%%%
%%%%%%%%%%%%%%%%%%%%%%%%%%%%%%%%%%%%%%%%%%

\subsection{Distributed Iterative Hessian Sketch}
We have evaluated the performance of the distributed IHS algorithm on AWS Lambda. 
% The way we have set up the implementation is as follows. 
In the implementation, each serverless function is responsible for solving one sketched problem per iteration. Worker nodes wait idly once they finish their computation for that iteration until the next iterate $x_{t+1}$ becomes available. The master node is implemented as another AWS Lambda function and is responsible for collecting and averaging the worker outputs and broadcasting the next iterate $x_{t+1}$.

\begin{figure}%[!t]
\begin{center}
\centerline{\includegraphics[width=0.5\columnwidth]{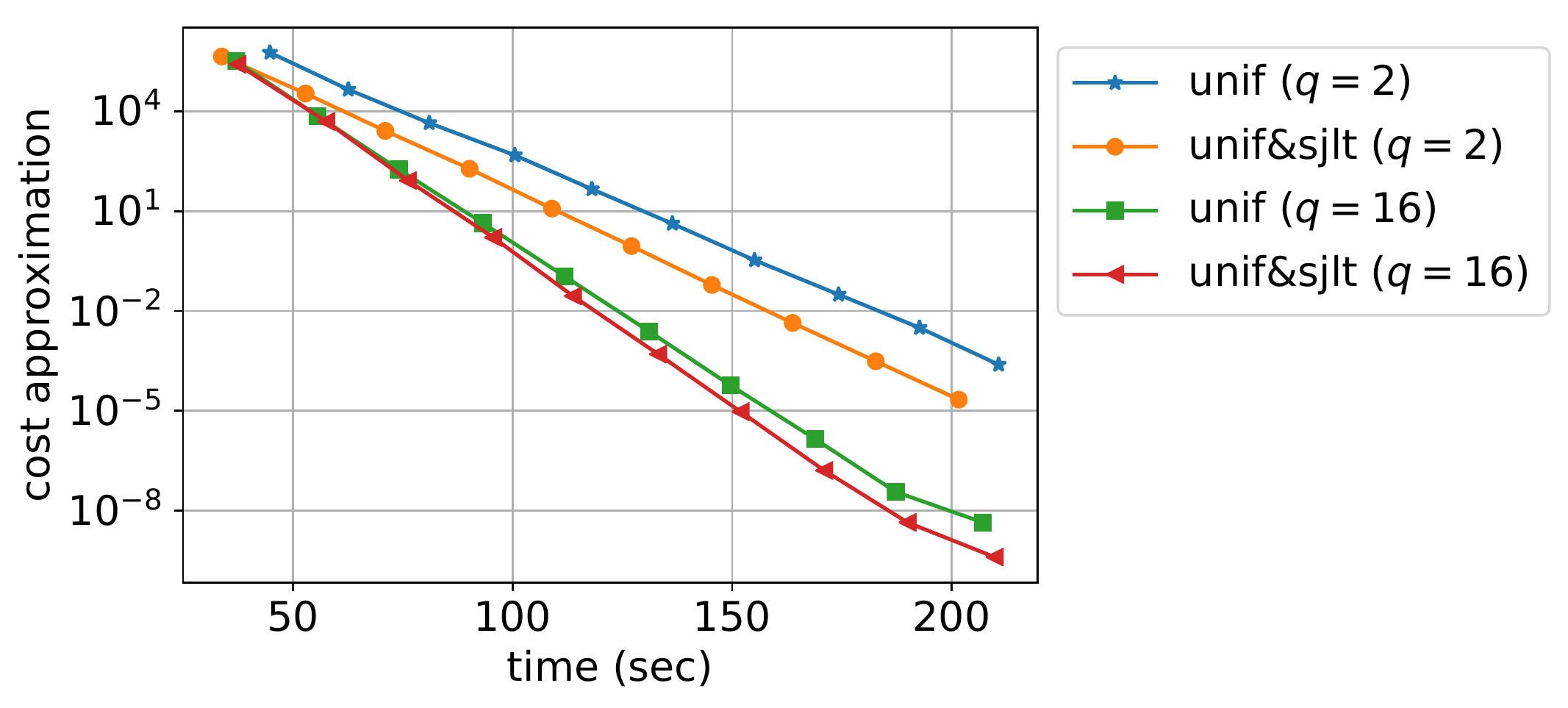}}
\caption{Cost approximation $(f(x_t)-f(x^*))/f(x^*)$ vs time for the distributed IHS algorithm running on AWS Lambda to fit the linear least squares model on randomly generated data. Unif is short for uniform sampling and unif\&sjlt is short for hybrid sketch where uniform sampling is followed by SJLT. Problem parameters are as follows: $n=250000$, $d=500$, $m=6000$, $m_2=20000$, and $q$ as specified in the legend.}
\label{fig:IHS_plot_aws_lambda}
\end{center}
\end{figure}

Figure \ref{fig:IHS_plot_aws_lambda} shows the scaled difference between the cost for the $t$'th iterate and the optimal cost (i.e. $(f(x_t)-f(x^*))/f(x^*)$) against wall-clock time for the distributed IHS algorithm given in Algorithm \ref{alg:dist_ihs}. Due to the relatively small size of the problem, we have each worker compute the exact gradient without requiring an additional communication round per iteration to form the full gradient. We note that, for problems where it is not feasible for workers to form the full gradient due to limited memory, one can include an additional communication round where each worker sends their local gradient to the master node, and the master node forms the full gradient and distributes it to the worker nodes. 

\subsection{Inequality Constrained Optimization}
\begin{figure}%[!t]
\begin{center}
\centerline{\includegraphics[width=0.45\columnwidth]{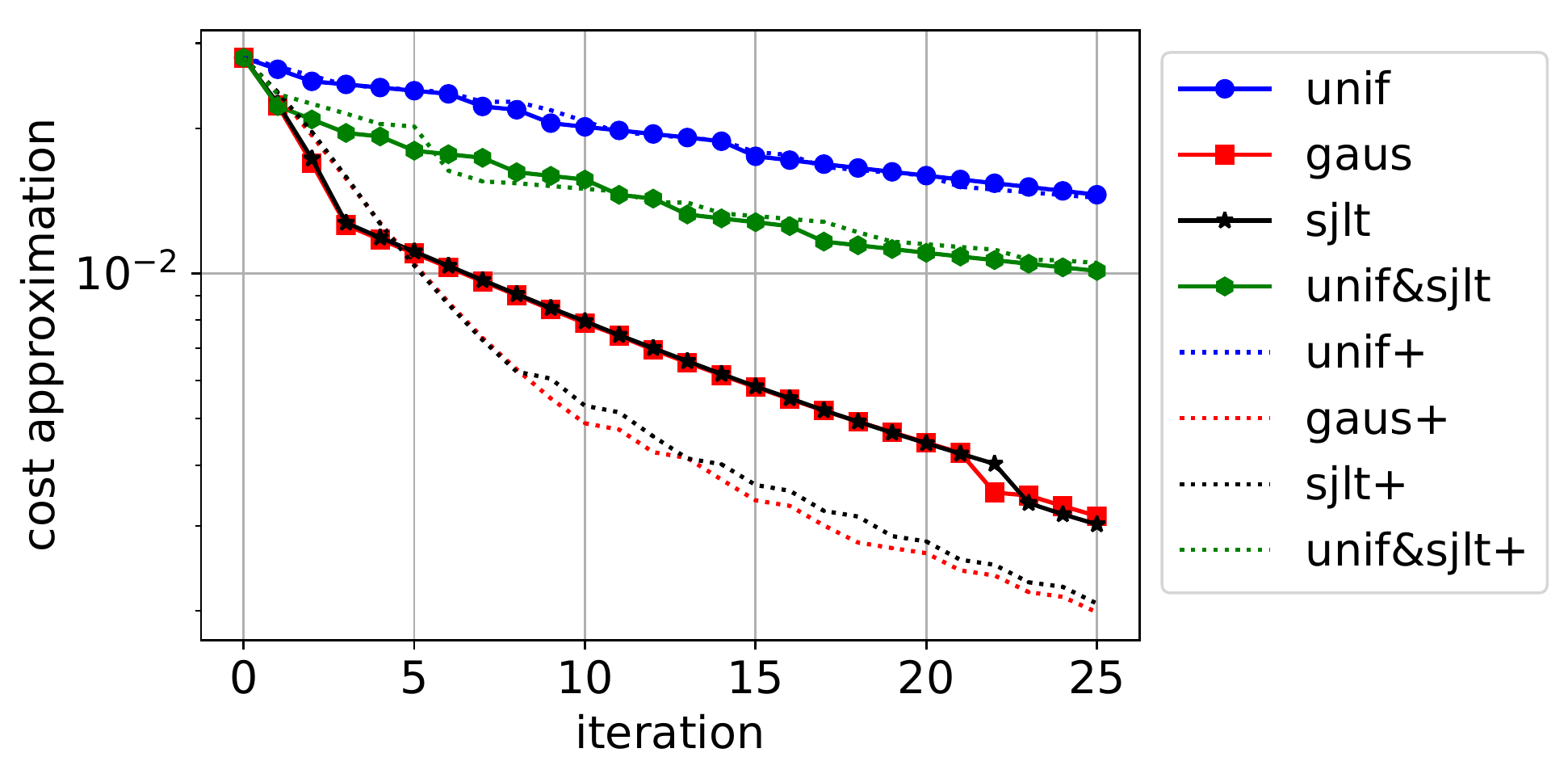}}
\caption{Plot shows cost approximation of the iterate $x_t$ (i.e., $(f(x_t)-f(x^*))/f(x^*)$) against iteration number $t$ for various sketches in solving the log-barrier version of an inequality constrained optimization problem given in \eqref{log_barrier_opt_prob_2}. Abbreviations used in the plot are as follows. Unif: Uniform sampling, gaus: Gaussian sketch, unif\&sjlt: Hybrid sketch where uniform sampling is followed by SJLT. The abbreviations followed by $+$'s in the legend refer to the bias corrected versions. Problem parameters are as follows: $n=500$, $d=200$, $\lambda_1=1000$, $m=50$, $m_2=8m=400$, $q=10$, $\lambda=0.01$, $s=10$. }
\label{fig:log_barrier_experiment}
\end{center}
\end{figure}

Figure \ref{fig:log_barrier_experiment} compares the error performance of sketches with and without bias correction for the distributed Newton sketch algorithm when it is used to solve the log-barrier penalized problem given in \eqref{log_barrier_opt_prob_2}. For each sketch, we have plotted the performance for $\lambda_2=\lambda_1$ and the bias corrected version $\lambda_2=\lambda_2^*$ (see Theorem \ref{opt_lambda_2_newton}). The bias corrected versions are shown as the dotted lines.
In these experiments, we have set $\sigma$ to the minimum of the singular values of $DA$ as we have observed that setting $\sigma$ to the minimum of the singular values of $DA$ performed better than setting it to their mean. 

Even though the bias correction formula is derived for Gaussian sketch, we observe that it improves the performance of SJLT as well. We see that Gaussian sketch and SJLT perform the best out of the 4 sketches we have experimented with. We note that computational complexity of sketching for SJLT is lower than it is for Gaussian sketch, and hence the natural choice would be to use SJLT in this case. %We note that we set the sparsity parameter to $s=10$.

\section{Discussion} \label{sec:discussion}

In this work, we study averaging sketched solutions for linear least squares problems and averaging for randomized second order optimization algorithms. We discuss distributed sketching methods from the perspectives of convergence, bias, privacy, and distributed computing. Our results and numerical experiments suggest that distributed sketching methods offer a competitive straggler-resilient solution for solving large scale problems.

We have shown that for problems involving regularization, averaging requires a more detailed analysis compared to problems without regularization. When the regularization term is not scaled properly, the resulting estimators are biased, and averaging a large number of independent sketched solutions does not converge to the true solution. We have provided closed-form formulas for scaling the regularization coefficient to obtain unbiased estimators that guarantee convergence to the optimum. %We also extended our analysis to non-identical sketch dimensions for heterogeneous computing environments. 

\appendices

% \section{}
\section{Proofs} \label{sec:appendix_proofs}
This section contains the proofs for the lemmas and theorems stated in the paper.

\subsection{Proofs of Theorems and Lemmas in Section \ref{sec:quad_problems}}

\begin{IEEEproof} [Proof of Theorem \ref{Thm:AvgGaussian}] 
% We first derive a result on the expected error of the averaged solution that relates the error of a single estimator $\hat{x}_k$ to that of the averaged solution $\bar{x}$.
% \begin{align*}
% \Exs &[\|A\bar x - Ax^*\|_2^2] = \Exs \left[ \left \|\frac{1}{q} \sum_{k=1}^q (A \hat x_k - \Exs[A \hat x_k]) \right\|_2^2 \right] \\
% & = \frac{1}{q^2}  \Exs \sum_{k=1}^q \sum_{l=1}^q \inprod{A \hat x_k - \Exs[A \hat x_k]}{A \hat x_l - \Exs[A  \hat x_l]}\\
% & = \frac{1}{q^2} \sum_{k=1}^q \Exs[\|A \hat x_k - \Exs [ A \hat x_k]\|^2]\\
% &= \frac{1}{q} \Exs[\|A \hat x_1 -  A \xstar]\|_2^2,
% \end{align*}
Since the Gaussian sketch estimator is unbiased (i.e., $\Exs[\hat{x}_k]=x^*$), Lemma \ref{expected_obj_val_diff} reduces to $
    \Exs[f(\bar{x})] - f(x^*) = \frac{1}{q} \Exs[\|A \hat x_1 -  A x^*]\|_2^2.
$
By Lemma \ref{gaussian_one_sketch}, the error of the averaged solution conditioned on the events that $E_k = A^TS_k^TS_kA \succ 0$, $\forall i=1,\dots,q$ can exactly be written as
\begin{align*}
    \Exs[ \| A(\bar{x} - x^*) \|_2^2 | E_1 \cap \dots\cap E_q ] = \frac{1}{q} \frac{d}{m-d-1} f(x^*).
\end{align*}
Using Markov's inequality, it follows that
\begin{align*}
    P(\| A(\bar{x} - x^*) \|_2^2 \geq a |  E_1 \cap \dots\cap E_q) \leq \frac{1}{qa} \frac{d}{m-d-1} f(x^*).
\end{align*}
The LHS can be lower bounded as 
\begin{align}
    \frac{P(\| A(\bar{x} - x^*) \|_2^2 \geq a \cap (\bigcap_{k=1}^q E_k))}{P(\bigcap_{k=1}^q E_k)} &\geq \frac{P(\| A(\bar{x} - x^*) \|_2^2 \geq a) + P(\bigcap_{k=1}^q E_k) - 1}{P(\bigcap_{k=1}^q E_k)} \label{temp_1}\\
    & = \frac{P(\| A(\bar{x} - x^*) \|_2^2 \geq a) + P(E_1)^q - 1}{P(E_1)^q}, \label{temp_2}
\end{align}
where we have used the identity $P(A\cap B)\geq P(A)+P(B)-1$ in \eqref{temp_1} and the independence of the events $E_k$ in \eqref{temp_2}. It follows
\begin{align*}
    P(\| A(\bar{x} - x^*) \|_2^2 \leq a) \geq P(E_1)^q \left(1 - \frac{1}{qa} \frac{d}{m-d-1} f(x^*) \right).
\end{align*}
Setting $a=f(x^*)\frac{\epsilon}{q}$ and plugging in $P(E_1)=1$, which holds for $m\ge d$, we obtain
\begin{align*}
    P\left(\frac{\| A(\bar{x} - x^*) \|_2^2}{f(x^*)} \leq \frac{\epsilon}{q}\right)
    \geq \left(1 - \frac{d/\epsilon}{m-d-1} \right).
\end{align*}
\end{IEEEproof}

\begin{lemma} \label{lem:smallest_sing_val_bound}
For the Gaussian sketching matrix $S$ and the matrix $U$ whose columns are orthonormal, the smallest singular value of the product $SU$ is bounded as:
\begin{align}
P(\sigma_{\min}(SU) \leq 1-\sqrt{d/m} - \delta) \leq \exp(-m\delta^2/2).
\end{align}
This result follows from the concentration of Lipshitz functions of Gaussian random variables.
%http://www.stat.cmu.edu/~larry/=sml/random_matrix_theory
\end{lemma}

\begin{IEEEproof} [Proof of Lemma \ref{lem:trace_concent}]
% Let us define $m(z) = \frac{1}{d} \tr \left((U^TS^TSU - \epsilon i I)^{-1}\right)$. 
We start by invoking a result on the concentration of Stieltjes transforms, namely Lemma 6 of \cite{karoui09concentration}. This implies that
\begin{align} \label{eq:stieltjes_concent}
    P&\big(| \tr (U^TS^TSU - \epsilon i I)^{-1} -\Exs[ \tr ( U^TS^TSU - \epsilon i I)^{-1}] | > t \big) \leq 4\exp \left(-t^2\epsilon^2 / (16m)\right)
\end{align}
where $i = \sqrt{-1}$.
The matrix $U^TS^TSU$ can be written as a sum of rank-1 matrices as follows:
\begin{align}
    U^TS^TSU = \sum_{i=1}^m (\tilde{s}_i^TU)^T (\tilde{s}_i^TU)
\end{align}
where $\tilde{s}_i^T \in \mathbb{R}^{1\times n}$ is the $i$'th row of $S$. The vectors $(\tilde{s}_i^TU)^T$ are independent as required by Lemma 6 of \cite{karoui09concentration}.
Next, we will bound the difference between the trace terms $\tr ( (U^TS^TSU)^{-1})$ and  $\tr ( (U^TS^TSU - \epsilon i I)^{-1})$. First, we let the SVD decomposition of $SU$ be $\tilde{U} \tilde{\Sigma} \tilde{V}^T$. Then, we have the following relations:
\begin{align}
    SU &= \tilde{U} \tilde{\Sigma} \tilde{V}^T \nonumber \\
    U^TS^TSU &= \tilde{V} \tilde{\Sigma}^2 {\tilde{V}}^T \nonumber \\
    (U^TS^TSU)^{-1} &= \tilde{V} \tilde{\Sigma}^{-2} {\tilde{V}}^T \nonumber \\
    U^TS^TSU - \epsilon i I &= \tilde{V} \diag(\tilde{\sigma}_j^2-\epsilon i) {\tilde{V}}^T \nonumber \\
    (U^TS^TSU - \epsilon i I)^{-1} &= \tilde{V} \diag\left(\frac{1}{\tilde{\sigma}_j^2-\epsilon i}\right) {\tilde{V}}^T
\end{align}
where the notation $\diag(\tilde{\sigma}_j^2)$ refers to a diagonal matrix with diagonal entries equal to $\tilde{\sigma}_1^2, \tilde{\sigma}_2^2, \dots$. 
The difference between the trace terms can be bounded as follows:
\begin{align}
    &\left|\tr ( (U^TS^TSU)^{-1}) - \tr ( (U^TS^TSU - \epsilon i I)^{-1})\right| = \left| \tr\left(\tilde{V} \diag\left(\frac{1}{\tilde{\sigma}_j^2} - \frac{1}{\tilde{\sigma}_j^2-\epsilon i}\right) {\tilde{V}}^T \right) \right| \nonumber \\
    &= \left| \sum_{j=1}^d \left(\frac{1}{\tilde{\sigma}_j^2} - \frac{1}{\tilde{\sigma}_j^2-\epsilon i}\right) \right| = \left| \sum_{j=1}^d \left( \frac{-\epsilon i}{\tilde{\sigma}_j^2 (\tilde{\sigma}_j^2-\epsilon i)}\right) \right| \nonumber \\
    &\leq \sum_{j=1}^d \left| \frac{-\epsilon i}{\tilde{\sigma}_j^2 (\tilde{\sigma}_j^2-\epsilon i)} \right| = \sum_{j=1}^d \frac{\epsilon }{\tilde{\sigma}_j^2 |\tilde{\sigma}_j^2-\epsilon i|}  \nonumber \\
    &\leq \sum_{j=1}^d \frac{\epsilon}{\tilde{\sigma}_j^4} \leq  \sum_{j=1}^d \frac{\epsilon}{\tilde{\sigma}_j^2 \sigma^2_{\min}(SU)} \nonumber \\
    &= \frac{\epsilon}{\sigma^2_{\min}(SU)} \sum_{j=1}^d \frac{1}{\tilde{\sigma}_j^2} = \frac{\epsilon}{\sigma^2_{\min}(SU)} \tr\left((U^TS^TSU)^{-1} \right) \,.
\end{align}
Equivalently, we have
\begin{align} \label{eq:trace_sandwich}
    &\big(1-\frac{\epsilon}{\sigma^2_{\min}(SU)}\big) \tr((U^TS^TSU)^{-1}) \leq \tr ( (U^TS^TSU - \epsilon i I)^{-1}) \leq \big(1+\frac{\epsilon}{\sigma^2_{\min}(SU)}\big) \tr((U^TS^TSU)^{-1}) \,.
\end{align}

\noindent The relations in \eqref{eq:stieltjes_concent} and \eqref{eq:trace_sandwich} together imply the following concentration result:
\begin{align}
    &P\big(| \tr ( (U^TS^TSU)^{-1}) -\Exs[ \tr ( (U^TS^TSU)^{-1}) ] | \leq t + \frac{\epsilon}{\sigma^2_{\min}(SU)}\left(\tr((U^TS^TSU)^{-1}) + \Exs[\tr((U^TS^TSU)^{-1})] \right) \big) \nonumber \\
    &\geq 1 - 4\exp \left(-t^2\epsilon^2 / (16m)\right)\,.
\end{align}

\noindent We know that  $\Exs[\tr((U^TS^TSU)^{-1})]=\frac{md}{m-d-1}$. Using the upper bound $\tr ((U^TS^TSU)^{-1}) \leq d\sigma_{\min}^{-2}(SU)$ yields
\begin{align} \label{eq:trace_concent_1}
    &P\Big(| \tr ( (U^TS^TSU)^{-1}) -\Exs[ \tr ( (U^TS^TSU)^{-1}) ] | \leq t + \frac{\epsilon}{\sigma^2_{\min}(SU)}\big( \frac{d}{\sigma^2_{\min}(SU)} + \frac{md}{m-d-1} \big) \Big) \nonumber \\
    &\geq 1 - 4\exp \left(-t^2\epsilon^2 / (16m)\right) .
\end{align}

The concentration inequality in \eqref{eq:trace_concent_1} and the concentration of the minimum singular value $P(\frac{1}{\sigma_{\min}(SU)} \leq \frac{1}{1-\sqrt{d/m} - \delta}) \geq 1-\exp(-m\delta^2/2)$ from Lemma \ref{lem:smallest_sing_val_bound} together imply that
\begin{align} \label{eq:trace_concent_2}
    &P\Big(| \tr ( (U^TS^TSU)^{-1}) -\Exs[ \tr ( (U^TS^TSU)^{-1}) ] |  \leq t + \frac{\epsilon}{(1-\sqrt{d/m}-\delta)^2}\big( \frac{d}{(1-\sqrt{d/m}-\delta)^2} + \frac{md}{m-d-1} \big) \Big) \nonumber \\
    &\geq 1 - 4\exp \left(-t^2\epsilon^2 / (16m)\right) - \exp(-m\delta^2/2).
\end{align}
We now simplify this bound by observing that
\begin{align}
     &\frac{d}{(1-\sqrt{d/m}-\delta)^2} = \frac{d}{(1+\frac{d}{m}-2\sqrt{\frac{d}{m}}) + \delta^2 - 2\delta (1-\sqrt{\frac{d}{m}})} \nonumber \\
     &\quad= \frac{md}{m+d-2\sqrt{md}+\delta^2m-2\delta m (1-\sqrt{\frac{d}{m}})} \nonumber \\
     &\quad \geq \frac{md}{m-d(2\sqrt{\alpha}-1)}  \nonumber \\
     &\quad \geq \frac{md}{m-d-1}
\end{align}
where we used $m \geq \alpha d$ in the third line. We have also used the following simple inequality
\begin{align}
    2\delta m (1-\sqrt{d/m}) &\geq 2\delta m (1-\frac{1}{\sqrt{2}}) \geq \delta m \geq \delta^2 m \,.
\end{align}
For the fourth line to be true, we need $m-2d\sqrt{\alpha}+d \leq m-d-1$. This is satisfied if $1\leq 2d(\sqrt{\alpha}-1)$ or $\alpha \geq (\frac{1}{2d}+1)^2$. We note that this is already implied by our assumption that $m \gtrsim d$.
Using the above observation, we simplify the bound as follows:
\begin{align} 
    &P\Big(| \tr ( (U^TS^TSU)^{-1}) -\Exs[ \tr ( (U^TS^TSU)^{-1}) ] |  \leq t + \frac{2d\epsilon}{(1-\sqrt{d/m}-\delta)^4} \Big) \geq 1 - 4\exp \left(-t^2\epsilon^2 / (16m)\right) - \exp(-m\delta^2/2) .
\end{align}
Scaling $\epsilon \leftarrow \epsilon \frac{(1-\sqrt{d/m}-\delta)^4}{2d}$ gives
\begin{align}
    &P\Big(| \tr ( (U^TS^TSU)^{-1}) -\Exs[ \tr ( (U^TS^TSU)^{-1}) ] | \leq t + \epsilon \Big) \geq 1 - 4 e^{-\frac{t^2\epsilon^2(1-\sqrt{d/m}-\delta)^8}{64md^2}} - e^{-m\delta^2/2}.
\end{align}
Setting $t=\epsilon$ and then $\epsilon \leftarrow \epsilon/2$ give us the final result:
\begin{align}
    &P\Big(| \tr ( (U^TS^TSU)^{-1}) -\Exs[ \tr ( (U^TS^TSU)^{-1}) ] | \leq \epsilon \Big) \geq 1 - 4 e^{-\frac{\epsilon^4(1-\sqrt{d/m}-\delta)^8}{2^{10}md^2}} - e^{-m\delta^2/2} .
\end{align}
\end{IEEEproof}

\begin{IEEEproof} [Proof of Theorem \ref{thm:concent_single_sketch}]
Let us consider the single sketch estimator $\hat{x}$. The relative error of $\hat{x}$ can be expressed as
\begin{align}
    f(\hat{x}) - f(x^*) &= \|A(\hat{x}-x^*)\|_2^2 \nonumber \\
    &= \|A(A^TS^TSA)^{-1}A^TS^TSb^{\perp}\|_2^2 \nonumber \\
    &= \|A(SA)^\dagger Sb^{\perp}\|_2^2
\end{align}
where the superscript $\hspace{-1mm}~^\dagger$ indicates the pseudo-inverse.
We will rewrite the error expression as:
\begin{align}
f(\hat{x}) - f(x^*) &= \left\Vert \frac{\|b^\perp\|_2}{\sqrt{m}} A(SA)^\dagger \frac{\sqrt{m}}{\|b^\perp\|_2} Sb^\perp \right\Vert_2^2 
\end{align}
and define $z := \frac{\sqrt{m}}{\|b^\perp\|_2}Sb^\perp$ and $M := \frac{\|b^\perp\|_2}{\sqrt{m}} A(SA)^\dagger$. This simplifies the previous expression
\begin{align}
    f&(\hat{x}) - f(x^*) = \|Mz\|_2^2 = z^TM^TMz .
\end{align}
Lemma \ref{lem:gausquad_concent} implies
\begin{align} \label{eq:quad_concent_M}
    &P\Big(f(\hat{x}) - f(x^*) - \Exs_{Sb^\perp}[f(\hat{x}) - f(x^*)] > 2\|M^TM\|_F \sqrt{\epsilon} + 2\|M^TM\|_2 \epsilon \Big| SU \Big) \leq e^{-\epsilon} .
\end{align}
The terms $\|M^TM\|_F$ and $\|M^TM\|_2$ reduce to the following expressions:
\begin{align}
    \|M^TM\|_F &= \frac{f(x^*)}{m} \sqrt{\tr((U^TS^TSU)^{-2})} \leq \frac{f(x^*)}{m} \frac{\sqrt{d}}{\sigma^2_{\min}(SU)}\,, \\
    \|M^TM\|_2 &= \frac{f(x^*)}{m} \frac{1}{\sigma^2_{\min}(SU)}\,.
\end{align}
Plugging these in \eqref{eq:quad_concent_M} and taking the expectation of both sides with respect to $SU$ give us
\begin{align} \label{eq:before_trace_concent}
    &P\Big(f(\hat{x}) - f(x^*) > \frac{f(x^*)}{m} \tr((U^TS^TSU)^{-1}) + 2 \frac{f(x^*)}{m\sigma^2_{\min}(SU)} ( \sqrt{\epsilon d} + \epsilon ) \Big) \leq e^{-\epsilon} .
\end{align}
Next, we combine \eqref{eq:before_trace_concent} with Lemma \ref{lem:trace_concent} and the concentration of the minimum singular value from Lemma \ref{lem:smallest_sing_val_bound} to obtain
\begin{align}
    &P\Big(f(\hat{x}) - f(x^*) > \frac{f(x^*)}{m} \Exs[\tr((U^TS^TSU)^{-1})] + 2 \frac{f(x^*)}{m(1-\sqrt{d/m}-s)^2} ( \sqrt{\epsilon d} + \epsilon ) + \frac{f(x^*)}{m} \gamma \Big) \nonumber \\
    &\leq e^{-\epsilon} + 4 e^{-\frac{\gamma^4(1-\sqrt{d/m}-\delta)^8}{2^{10}md^2}} + e^{-m\delta^2/2} + e^{-ms^2/2} .
\end{align}
Assuming that $\sqrt{\epsilon d} > \epsilon$, we can write:
\begin{align}
    &P\Big(f(\hat{x}) - f(x^*) > \frac{f(x^*)}{m} \Exs[\tr((U^TS^TSU)^{-1})] + \frac{4f(x^*)\sqrt{\epsilon d}}{m(1-\sqrt{d/m}-s)^2} + \frac{f(x^*)}{m} \gamma \Big) \nonumber \\
    &\leq e^{-\epsilon} + 4 e^{-\frac{\gamma^4(1-\sqrt{d/m}-\delta)^8}{2^{10}md^2}} + e^{-m\delta^2/2} + e^{-ms^2/2} .
\end{align}
Making the change of variable $\epsilon \leftarrow \epsilon^2/d$ yields:
\begin{align}
    &P\Big(f(\hat{x}) - f(x^*) > \frac{f(x^*)}{m} \Exs[\tr((U^TS^TSU)^{-1})] + \frac{4f(x^*)\epsilon}{m(1-\sqrt{d/m}-s)^2} + \frac{f(x^*)}{m} \gamma \Big) \nonumber \\
    &\leq e^{-\epsilon^2/d} + 4 e^{-\frac{\gamma^4(1-\sqrt{d/m}-\delta)^8}{2^{10}md^2}} + e^{-m\delta^2/2} + e^{-ms^2/2} .
\end{align}
Scaling $\epsilon \leftarrow \epsilon \frac{(1-\sqrt{d/m}-s)^2}{4}$ yields:
\begin{align}
    &P\Big(f(\hat{x}) - f(x^*) > \frac{f(x^*)}{m} \Exs[\tr((U^TS^TSU)^{-1})] + \frac{f(x^*)(\epsilon + \gamma)}{m}  \Big) \nonumber \\
    &\leq e^{-\frac{\epsilon^2(1-\sqrt{d/m}-s)^4}{16d}} + 4 e^{-\frac{\gamma^4(1-\sqrt{d/m}-\delta)^8}{2^{10}md^2}} + e^{-m\delta^2/2} + e^{-ms^2/2}\,.
\end{align}
Letting $\gamma = \epsilon$ and $s=\delta$, and then scaling $\epsilon \leftarrow \epsilon m/2$ lead to
\begin{align}
    &P\big(\frac{f(\hat{x})}{f(x^*)} > \frac{m-1}{m-d-1}  + \epsilon \big) \leq e^{-\frac{\epsilon^2m^2(1-\sqrt{d/m}-\delta)^4}{64d}} + 4 e^{-\frac{\epsilon^4m^4(1-\sqrt{d/m}-\delta)^8}{2^{14}md^2}} + 2e^{-m\delta^2/2}.
\end{align}
Let us also set $\delta=\epsilon$:
\begin{align}
    &P\big(\frac{f(\hat{x})}{f(x^*)} > \frac{m-1}{m-d-1}  + \epsilon \big) \leq e^{-\frac{\epsilon^2m^2(1-\sqrt{d/m}-\epsilon)^4}{64d}} + 4 e^{-\frac{\epsilon^4m^4(1-\sqrt{d/m}-\epsilon)^8}{2^{14}md^2}} + 2e^{-m\epsilon^2/2} \,.
\end{align}
Under our assumption that $m \gtrsim d$, the above probability is bounded by
\begin{align}
    &P\big(\frac{f(\hat{x})}{f(x^*)} > \frac{m-1}{m-d-1}  + \epsilon \big) \leq C_1 e^{-C_2\epsilon^4 m} \,.
\end{align} 
%

%%%
We now move on to find a lower bound. We note that Lemma \ref{lem:gausquad_concent} can be used as follows to find a lower bound
\begin{align}
    &P\big(z^T(-G)z - \Exs[z^T(-G)z] > 2\|(-G)\|_F \sqrt{\epsilon} + 2\|(-G)\|_2\epsilon \big) = P\big(z^TGz < \Exs[z^TGz] - 2\|G\|_F \sqrt{\epsilon} - 2\|G\|_2\epsilon \big) \leq e^{-\epsilon} \,.
\end{align}
Therefore, we have the following lower bound for the error:
\begin{align}
    &P\Big(f(\hat{x}) - f(x^*) < \Exs_{Sb^\perp}[f(\hat{x}) - f(x^*)] - 2\|M^TM\|_F \sqrt{\epsilon} - 2\|M^TM\|_2 \epsilon \Big| SU \Big) \leq e^{-\epsilon} .
\end{align}
Using the upper bound for $\|M^TM\|_F$ and the exact expression for $\|M^TM\|_2$, and taking the expectation with respect to $SU$, we obtain the following
\begin{align}
    &P\Big(f(\hat{x}) - f(x^*) < \frac{f(x^*)}{m} \tr((U^TS^TSU)^{-1}) - 2 \frac{f(x^*)}{m\sigma^2_{\min}(SU)} ( \sqrt{\epsilon d} + \epsilon ) \Big) \leq e^{-\epsilon} .
\end{align}
Applying the same steps as before, we arrive at the following bound:
\begin{align}
    &P\big(\frac{f(\hat{x})}{f(x^*)} < \frac{m-1}{m-d-1} - \epsilon \big) \leq C_1 e^{-C_2\epsilon^4 m}.
\end{align}
\end{IEEEproof}

\begin{IEEEproof} [Proof of Theorem \ref{thm:concent_averaged_sketch}]
Consider the averaged estimator $\bar{x}=\frac{1}{q}\sum_{k=1}^q\hat{x}_k$:
\begin{align}
    f(\bar{x}) - f(x^*) &= \|A(\bar{x} - x^*)\|_2^2 = \left\Vert \frac{1}{q} \sum_{k=1}^q A (S_kA)^\dagger S_kb^\perp \right\Vert_2^2 \nonumber \\
    &= \left\Vert \frac{1}{q} \sum_{k=1}^q M_k z_k \right\Vert_2^2
\end{align}
where we defined
\begin{align}
    M_k &:= \frac{\|b^\perp\|_2}{\sqrt{m}} A(S_kA)^\dagger, \quad z_k := \frac{\sqrt{m}}{\|b^\perp\|_2} S_k b^\perp .
\end{align}
We note that the entries of the vector $z$ are distributed as i.i.d. standard normal. We can manipulate the error expression as follows
\begin{align}
    \sum_{k=1}^q M_k z_k &= \begin{bmatrix} M_1 & \dots & M_q \end{bmatrix} \begin{bmatrix} z_1 \\ \vdots \\ z_q \end{bmatrix} = M z
\end{align}
and then
\begin{align}
    f(\bar{x}) - f(x^*) = \frac{1}{q^2} \|Mz\|_2^2 = z^T\big(\frac{1}{q^2} M^TM\big)z .
\end{align}
We will now find an equivalent expression for $\|\frac{1}{q^2} M^TM\|_F^2 = \frac{1}{q^4} \|M^TM\|_F^2$:
\begin{align}
    \|M^TM\|_F^2 &= \tr (M^TMM^TM) = \tr(MM^TMM^T) \nonumber \\
    &= \tr\left( \sum_{k=1}^q M_kM_k^T \sum_{l=1}^q M_lM_l^T \right) \nonumber \\
    % &= \frac{f^2(x^*)}{m^2} \tr\left(\sum_{k=1}^q {(S_kU)^\dagger}^T (S_kU)^\dagger \sum_{l=1}^q {(S_lU)^\dagger}^T (S_lU)^\dagger \right) \nonumber \\
    &= \frac{f^2(x^*)}{m^2} \sum_{k=1}^q \sum_{l=1}^q \tr\left({(S_kU)^\dagger}^T (S_kU)^\dagger {(S_lU)^\dagger}^T (S_lU)^\dagger \right) \nonumber \\
    &= \frac{f^2(x^*)}{m^2} \Big(\sum_{k=1}^q \tr((U^TS_k^TS_kU)^{-2}) + \sum_{k\neq l} \tr\left({(S_kU)^\dagger}^T (S_kU)^\dagger {(S_lU)^\dagger}^T (S_lU)^\dagger \right) \Big) \nonumber \\
    &\leq \frac{f^2(x^*)}{m^2} \Big(\sum_{k=1}^q \frac{d}{\sigma^4_{\min}(S_kU)} + \sum_{k \neq l} \sqrt{\tr\left({(S_kU)^\dagger}^T (S_kU)^\dagger\right)} \sqrt{\tr\left({(S_lU)^\dagger}^T (S_lU)^\dagger\right)} \Big) \nonumber \\
    &= \frac{f^2(x^*)}{m^2} \Big(\sum_{k=1}^q \frac{d}{\sigma^4_{\min}(S_kU)} + \sum_{k \neq l} \sqrt{\tr\left((U^TS_k^TS_kU)^{-1}\right)} \sqrt{\tr\left((U^TS_l^TS_lU)^{-1}\right)} \Big) \nonumber \\
    &\leq \frac{f^2(x^*)}{m^2} \Big(\sum_{k=1}^q \frac{d}{\sigma^4_{\min}(S_kU)} + \sum_{k\neq l} \frac{d}{\sigma_{\min}(S_kU) \sigma_{\min}(S_lU)} \Big)
\end{align}
where we have used the Cauchy-Schwarz inequality to bound the trace of the matrix product.

Next, we look at $\|\frac{1}{q^2}M^TM\|_2=\frac{1}{q^2}\|M^TM\|_2$:
\begin{align}
    \|M^TM\|_2 &=  \|M\|_2^2 = \|MM^T\|_2 = \left\Vert\sum_{k=1}^q M_kM_k^T \right\Vert_2 \nonumber \\
    &= \frac{\|b^\perp\|_2^2}{m}\left\Vert \sum_{k=1}^q {(S_kU)^\dagger}^T (S_kU)^\dagger \right\Vert_2 \nonumber \\
    &\leq \frac{\|b^\perp\|_2^2}{m} \sum_{k=1}^q \| {(S_kU)^\dagger}^T (S_kU)^\dagger \|_2 \nonumber \\
    &= \frac{\|b^\perp\|_2^2}{m} \sum_{k=1}^q \frac{1}{\sigma^2_{\min}(S_kU)}.
\end{align}
By Lemma \ref{lem:gausquad_concent}, we have
\begin{align}
    &P\Big(f(\bar{x}) - f(x^*) - \Exs_{S_kb^\perp}[f(\bar{x}) - f(x^*)] > \frac{2}{q^2}\|M^TM\|_F \sqrt{\epsilon} + \frac{2}{q^2} \|M^TM\|_2 \epsilon \Big| SU_k \Big) \leq e^{-\epsilon}
\end{align}
We can take expectation of both sides with respect to $SU_k$, $k=1,\dots,q$ to remove the conditioning.
The expectation term inside the probability is equal to
\begin{align}
    \Exs_{S_kb^\perp}[f(\bar{x}) - f(x^*)] &= \tr(\frac{1}{q^2}M^TM) = \frac{1}{q^2} \tr(MM^T) = \frac{f(x^*)}{mq^2} \sum_{k=1}^q \tr((U^TS_k^TS_kU)^{-1}).
\end{align}
Hence, we have
\begin{align}
    &P\Big(f(\bar{x}) - f(x^*) - \frac{f(x^*)}{mq^2} \sum_{k=1}^q \tr((U^TS_k^TS_kU)^{-1}) > \frac{2}{q^2}\|M^TM\|_F \sqrt{\epsilon} + \frac{2}{q^2} \|M^TM\|_2 \epsilon \Big) \leq e^{-\epsilon} \,.
\end{align}
We now use the concentration bound of the trace term in Lemma \ref{lem:trace_concent} to obtain
\begin{align}
    &P\Big(f(\bar{x}) - f(x^*) - \frac{f(x^*)}{mq^2} \sum_{k=1}^q \Exs[\tr((U^TS_k^TS_kU)^{-1})] > \frac{2}{q^2}\|M^TM\|_F \sqrt{\epsilon} + \frac{2}{q^2} \|M^TM\|_2 \epsilon + \frac{f(x^*)}{mq^2} q\gamma \Big) \nonumber \\
    &\leq e^{-\epsilon} + 1 - \big(1 - 4 e^{-\frac{\gamma^4(1-\sqrt{d/m}-\delta)^8}{2^{10}md^2}} - e^{-m\delta^2/2}\big)^q 
\end{align}
where we have used the independence of $S_k$'s.
Using the upper bounds for $\|M^TM\|_F$ and $\|M^TM\|_2$, we make the following observation. If the minimum singular values of all of $S_kU$, $k=1,\dots,q$ satisfy $\sigma_{\min}(S_kU) > 1-\sqrt{d/m}-s$, which occurs with probability at least $(1-e^{-ms^2/2})^q$ from Lemma \ref{lem:smallest_sing_val_bound}, then the following holds
\begin{align}
    \|M^TM\|_F\sqrt{\epsilon} + \|M^TM\|_2 \epsilon &\leq \frac{f(x^*)}{m} \bigg( \sqrt{\frac{qd}{(1-\sqrt{d/m}-s)^4}+\frac{(q^2-q)d}{(1-\sqrt{d/m}-s)^2}} \sqrt{\epsilon} + \frac{q\epsilon}{(1-\sqrt{d/m}-s)^2} \bigg) \nonumber \\
    &\leq \frac{f(x^*)}{m} \frac{q(\sqrt{\epsilon d}+\epsilon)}{(1-\sqrt{d/m}-s)^2} \nonumber \\
    &\leq \frac{f(x^*)}{m} \frac{2q\sqrt{\epsilon d}}{(1-\sqrt{d/m}-s)^2}
\end{align}
where we used $1-\sqrt{d/m}-s < 1$ and $\epsilon \leq \sqrt{\epsilon d}$ in the last two inequalities.
Combining this observation with the probability bound, we arrive at
\begin{align}
    &P\Big(f(\bar{x}) - f(x^*) - \frac{f(x^*)}{mq^2} \sum_{k=1}^q \Exs[\tr((U^TS_k^TS_kU)^{-1})] > \frac{2}{q^2}\frac{f(x^*)}{m} \frac{2q\sqrt{\epsilon d}}{(1-\sqrt{d/m}-s)^2} + \frac{f(x^*)}{mq^2} q\gamma \Big) \nonumber \\
    &\leq e^{-\epsilon} + 1 - \big(1 - 4 e^{-\frac{\gamma^4(1-\sqrt{d/m}-\delta)^8}{2^{10}md^2}} - e^{-m\delta^2/2}\big)^q + 1 - (1-e^{-ms^2/2})^q.
\end{align}
Simplifying the expressions and making the change $\epsilon \leftarrow \epsilon^2/d$ lead to
\begin{align}
    &P\Big(f(\bar{x}) - f(x^*) - \frac{f(x^*)}{q}\frac{d}{m-d-1} > \frac{4f(x^*)}{mq} \frac{\epsilon}{(1-\sqrt{d/m}-s)^2} + \frac{f(x^*)}{mq} \gamma \Big) \nonumber \\
    &\leq e^{-\epsilon^2/d} + 2 - \big(1 - 4 e^{-\frac{\gamma^4(1-\sqrt{d/m}-\delta)^8}{2^{10}md^2}} - e^{-m\delta^2/2}\big)^q - (1-e^{-ms^2/2})^q
\end{align}
Next, we let $\epsilon \leftarrow \epsilon \frac{(1-\sqrt{d/m}-s)^2}{4}$ and obtain
\begin{align}
    &P\Big(f(\bar{x}) - f(x^*) - \frac{f(x^*)}{q}\frac{d}{m-d-1} > \frac{f(x^*)}{mq}(\epsilon+\gamma) \Big) \nonumber \\
    &\leq e^{-\frac{\epsilon^2(1-\sqrt{d/m}-s)^4}{16d}} + 2  - (1-e^{-ms^2/2})^q - \big(1 - 4 e^{-\frac{\gamma^4(1-\sqrt{d/m}-\delta)^8}{2^{10}md^2}} - e^{-m\delta^2/2}\big)^q \,.
\end{align}
We now set $\gamma=\epsilon$ and scale $\epsilon \leftarrow \epsilon \frac{mq}{2}$:
\begin{align}
    &P\Big(\frac{f(\bar{x})}{f(x^*)} > 1+ \frac{1}{q}\frac{d}{m-d-1} + \epsilon \Big) \leq e^{-\frac{\epsilon^2m^2q^2(1-\sqrt{d/m}-s)^4}{64d}} + 2 - (1-e^{-ms^2/2})^q - \big(1 - 4 e^{-\frac{\epsilon^4m^4q^41-\sqrt{d/m}-\delta)^8}{2^{14}md^2}} - e^{-m\delta^2/2}\big)^q .
\end{align}
Let us use the assumption $m \gtrsim d$ to further simplify:
\begin{align}
    &P\Big(\frac{f(\bar{x})}{f(x^*)} > 1+ \frac{1}{q}\frac{d}{m-d-1} + \epsilon \Big) \leq C_1e^{-C_2(q\epsilon)^2m} + 2 - (1-e^{-ms^2/2})^q - \big(1 - 4 e^{-\frac{\epsilon^4m^4q^4(1-\sqrt{d/m}-\delta)^8}{2^{14}md^2}} - e^{-m\delta^2/2}\big)^q .
\end{align}
We can use Bernoulli's inequality $(1-e^{-ms^2/2})^q \geq 1-qe^{-ms^2/2}$ and arrive at:
\begin{align}
    &P\Big(\frac{f(\bar{x})}{f(x^*)} > 1+ \frac{1}{q}\frac{d}{m-d-1} + \epsilon \Big) \leq C_1e^{-C_2(q\epsilon)^2m}  +C_3q e^{-C_4(q\epsilon)^4m(1-\sqrt{d/m}-s)^8} + 2qe^{-ms^2/2}
\end{align}
where we also set $\delta=s$.
Picking $s = q\epsilon$ yields:
\begin{align}
    &P\Big(\frac{f(\bar{x})}{f(x^*)} > 1+ \frac{1}{q}\frac{d}{m-d-1} + \epsilon \Big) \leq C_1e^{-C_2(q\epsilon)^2m}  +C_3q e^{-C_4(q\epsilon)^4m(1-\sqrt{d/m}-q\epsilon)^8} + 2qe^{-m(q\epsilon)^2/2}\,.
\end{align}
Finally, we obtain the following simpler expression for the bound:
\begin{align}
    &P\Big(\frac{f(\bar{x})}{f(x^*)} > 1+ \frac{1}{q}\frac{d}{m-d-1} + \epsilon \Big) \leq qC_1e^{-C_2(q\epsilon)^4m}  
\end{align}
where we redefine the constants $C_1,C_2$.

Lower bound can be obtained by applying the same steps:
\begin{align}
    &P\Big(\frac{f(\bar{x})}{f(x^*)} < 1+ \frac{1}{q}\frac{d}{m-d-1} - \epsilon \Big) \leq qC_1e^{-C_2(q\epsilon)^4m} .
\end{align}
\end{IEEEproof}

\begin{IEEEproof} [Proof of Theorem \ref{thm:fisher_lower_bounds}]
\noindent \textit{i) For unbiased estimators:}  
We begin by invoking Lemma \ref{expected_obj_val_diff} which gives us
\begin{align}
    \Exs [\|A(\bar{x} - x^*) \|_2^2] &= \frac{1}{q} \Exs [\|A(\hat{x}_1 - x^*) \|_2^2]
\end{align}
since we assume that $\hat{x}_k$'s are unbiased, i.e., $\Exs[\hat{x}_k] = x^*$ for $k=1,\dots,q$. 
%Then the result follows from Lemma \ref{lem:error_lb_singlesketch}. 

This shows that the error of the averaged estimator is equal to the error of the single sketch estimator scaled by $\frac{1}{q}$. Hence, we can leverage the lower bound result for the single sketch estimator. We now give the details of how to find the lower bound for the single sketch estimator which was first shown in \cite{sridhar20lowerbounds}.

The Fisher information matrix for estimating $x^*$ from $Sb$ can be constructed as follows:
\begin{align} \label{eq:fisher_info_matrix}
    I(Sb; x^*) = \Exs_{Sb}[\nabla_{x^*} \log g(Sb; x^*) \nabla_{x^*} \log g(Sb; x^*)^T]
\end{align}
where $g(Sb;x^*)$ is the probability density function of $Sb$ conditioned on $SA$, i.e., $g(Sb;x^*)$ is the multivariate Gaussian distribution with $\mathcal{N}(SAx^*, \frac{1}{m}\|b^\perp\|_2^2I_m)$. We recall the definition that $b^\perp := b - Ax^*$.

The gradient of the logarithm of the probability density function with respect to $x^*$ is
\begin{align}
    \nabla_{x^*} \log g(Sb; x^*) = \frac{m}{\|b^\perp\|_2^2} (SA)^T (Sb-SAx^*) \,.
\end{align}
Plugging this in \eqref{eq:fisher_info_matrix} yields
\begin{align}
    I(Sb; x^*) &= \frac{m^2}{\|b^\perp\|_2^4}\Exs_{Sb}\left[(SA)^T (Sb-SAx^*) (Sb-SAx^*)^T SA \right] \nonumber \\
    % &= \frac{m^2}{\|b^\perp\|_2^4} (SA)^T \Exs_{Sb}\left[(Sb-SAx^*) (Sb-SAx^*)^T \right] SA \nonumber \\
    &= \frac{m^2}{\|b^\perp\|_2^4} (SA)^T \big(\frac{\|b^\perp\|_2^2}{m} I_m \big) SA \nonumber \\
    &= \frac{m}{\|b^\perp\|_2^2} A^TS^TSA \,.
\end{align}
We note that the expectation of the error conditioned on $SA$ is equivalent to
\begin{align}
    \Exs[\|A(\hat{x}-x^*)\|_2^2 | SA] &= \Exs[\tr(A(\hat{x}-x^*)(\hat{x}-x^*)^T A^T) | SA] \nonumber \\
    &= \tr \big(A\Exs[(\hat{x}-x^*)(\hat{x}-x^*)^T| SA] A^T \big)  \nonumber \\
    &\geq \frac{\|b^\perp\|_2^2}{m} \tr(A(A^TS^TSA)^{-1}A^T)\,,
\end{align}
where the last line follows from the Cram\'er lower bound  \cite{borovkov98mathstat} given by
\begin{align}
    \Exs[(\hat{x}-x^*)(\hat{x}-x^*)^T|SA] \succeq I^{-1}(Sb; x^*) \,.
\end{align}
Taking the expectation of both sides with respect to $SA$ yields
\begin{align}
    \Exs&[\|A(\hat{x}-x^*)\|_2^2 ] \geq f(x^*) \frac{d}{m-d-1} \, .
\end{align}
Hence, we arrive at
\begin{align}
    \Exs [\|A(\bar{x} - x^*) \|_2^2] \geq \frac{f(x^*)}{q} \frac{d}{m-d-1} \,.
\end{align}

%%%%
\noindent \textit{ii) For any estimator:} We again invoke Lemma \ref{expected_obj_val_diff}, which gives us
\begin{align}
    \Exs[\|A(\bar{x} - x^*) \|_2^2] &= \frac{1}{q} \Exs [\|A(\hat{x}_1 - x^*) \|_2^2] + \frac{q-1}{q} \|\Exs[A\hat{x}_1] - Ax^* \|_2^2 \nonumber \\
    &\geq \frac{f(x^*)}{q}\frac{d}{m} + \frac{q-1}{q} \|\Exs[A\hat{x}_1] - Ax^* \|_2^2 \nonumber \\
    &\geq \frac{f(x^*)}{q}\frac{d}{m}
\end{align}
where the inequality on the third line follows from Lemma \ref{lem:error_lb_singlesketch_biased}.
\end{IEEEproof}

\begin{IEEEproof} [Proof of Theorem \ref{thm:IHS_error_decay}]
The update rule for distributed IHS is given as
\begin{align}
    x_{t+1} = x_t - \mu \frac{1}{q} \sum_{k=1}^q (A^TS_{t,k}^TS_{t,k}A)^{-1}A^T(Ax_t-b).
\end{align}
Let us decompose $b$ as $b=Ax^* + b^\perp$ and note that $A^T b^\perp = 0$ which gives us:
\begin{align}
    x_{t+1} = x_t - \mu \frac{1}{q} \sum_{k=1}^q (A^TS_{t,k}^TS_{t,k}A)^{-1}A^TA e_t.
\end{align}
Subtracting $x^*$ from both sides, we obtain an equation in terms of the error vector $e_t$ only:
\begin{align*}
    e_{t+1} &= e_t - \mu \frac{1}{q} \sum_{k=1}^q (A^TS_{t,k}^TS_{t,k}A)^{-1}A^TA e_t \\
    &= \left(I - \mu \frac{1}{q} \sum_{k=1}^q (A^TS_{t,k}^TS_{t,k}A)^{-1}A^TA \right) e_t.
\end{align*}
Let us multiply both sides by $A$ from the left and define $Q_{t,k} \coloneqq A(A^TS_{t,k}^TS_{t,k}A)^{-1}A^T$ and we will have the following equation:
\begin{align*}
    e^A_{t+1} &= \left(I - \mu \frac{1}{q} \sum_{k=1}^q Q_{t,k} \right) e^A_t.
\end{align*}

We now analyze the expectation of $\ell_2$ norm of $e^A_{t+1}$: 
\begin{align} \label{expected_error_expanded}
    \Exs [ \|e^A_{t+1}\|_2^2 ] &= \Exs \left[ \left\Vert \frac{1}{q} \sum_{k=1}^q (I - \mu Q_{t,k}) e^A_t \right\Vert_2^2 \right] \nonumber \\
    &= \frac{1}{q^2} \Exs \left[ \sum_{k=1}^q \sum_{l=1}^q \langle (I - \mu Q_{t,k}) e^A_t, (I - \mu Q_{t,l}) e^A_t \rangle \right] \nonumber \\
    &= \frac{1}{q^2} \sum_{k=1}^q \sum_{l=1}^q \Exs \left[ \langle (I - \mu Q_{t,k}) e^A_t, (I - \mu Q_{t,l}) e^A_t \rangle \right].
\end{align}
The contribution for $k \neq l$ in the double summation of \eqref{expected_error_expanded} is equal to zero because for $k \neq l$, we have
\begin{align*}
    \Exs \left[ \langle (I - \mu Q_{t,k}) e^A_t, (I - \mu Q_{t,l}) e^A_t \rangle \right] &= \langle \Exs [ (I - \mu Q_{t,k}) e^A_t], \Exs[ (I - \mu Q_{t,l}) e^A_t] \rangle \\
    & = \langle \Exs [ (I - \mu Q_{t,k}) e^A_t], \Exs[ (I - \mu Q_{t,k}) e^A_t] \rangle \\
    & = \left\Vert \Exs [ (I - \mu Q_{t,k}) e^A_t] \right\Vert_2^2.
\end{align*}
The term in the last line above is zero for $\mu=\frac{1}{\theta_1}$:
\begin{align*}
    \Exs [ (I - \mu Q_{t,k}) e^A_t] &= \Exs [(I - \mu A(A^TS_{t,k}^TS_{t,k}A)^{-1}A^T) e^A_t ] \\
    &= (I - \mu \theta_1 A(A^TA)^{-1} A^T) e^A_t  \\
    &= (I - \mu \theta_1 UU^T)e^A_t \\
    &= (I - UU^T)e^A_t = 0
\end{align*}
where we used $A=U\Sigma V^T$.
For the rest of the proof, we assume that we set $\mu = 1/\theta_1$. Now that we know the contribution from terms with $k \neq l$ is zero, the expansion in \eqref{expected_error_expanded} can be rewritten as
\begin{align*}
    \Exs [\|e^A_{t+1}\|_2^2 ] &= \frac{1}{q^2} \sum_{k=1}^q \Exs \left[ \langle (I - \mu Q_{t,k}) e^A_t, (I - \mu Q_{t,k}) e^A_t \rangle \right] \\
    &= \frac{1}{q^2} \sum_{k=1}^q \Exs[\| (I - \mu Q_{t,k}) e^A_t \|_2^2] \\
    &= \frac{1}{q} \Exs[\| (I - \mu Q_{t,1}) e^A_t \|_2^2] \\
    &= \frac{1}{q} \left( \|e^A_t\|_2^2 + \mu^2 \Exs[ \|Q_{t,1} e^A_t\|_2^2] - 2\mu (e^A_t)^T \Exs[Q_{t,1}]e^A_t \right) \\
    &= \frac{1}{q} \left( \mu^2 \Exs[ \|Q_{t,1} e^A_t\|_2^2] - \|e^A_t\|_2^2 \right) \\
    &= \frac{1}{q} \left( \mu^2 (e^A_t)^T \Exs[Q_{t,1}^T Q_{t,1}] e^A_t - \|e^A_t\|_2^2 \right) .
\end{align*}
The term $\Exs[Q_{t,1}^T Q_{t,1}]$ can be simplified using SVD decomposition $A=U\Sigma V^T$. This gives us $Q_{t,k} = U(U^TS_{t,k}^TS_{t,k}U)^{-1}U^T$ and furthermore we have:
\begin{align*}
    \Exs [Q_{t,1}^T Q_{t,1}] &= \Exs[ U(U^TS_{t,1}^TS_{t,1}U)^{-1}U^T U(U^TS_{t,1}^TS_{t,1}U)^{-1}U^T ] \\
    &=  \Exs[ U(U^TS_{t,1}^TS_{t,1}U)^{-1} (U^TS_{t,1}^TS_{t,1}U)^{-1}U^T ] \\
    &= U \Exs[(U^TS_{t,1}^TS_{t,1}U)^{-2}] U^T \\
    &= \theta_2 UU^T.
\end{align*}
Plugging this in, we obtain:
\begin{align*}
    \Exs [ \|e^A_{t+1}\|_2^2 ] &= \frac{1}{q} \left( \theta_2 \mu^2 (e^A_t)^T  UU^T e^A_t - \|e^A_t\|_2^2 \right) \\
    % &= \frac{1}{q} \left( \theta_2 \mu^2 \|U^T e^A_t\|_2^2 - \|e^A_t\|_2^2 \right) \\
    &= \frac{1}{q} \left( \theta_2 \mu^2 \|e^A_t\|_2^2 - \|e^A_t\|_2^2 \right) \\
    &= \frac{\theta_2 \mu^2 - 1}{q}  \|e^A_t\|_2^2 \\
    &= \frac{1}{q} \left( \frac{\theta_2}{\theta_1^2} - 1 \right) \|e^A_t\|_2^2 \,.
\end{align*}
\end{IEEEproof}

\begin{IEEEproof} [Proof of Corollary \ref{ihs_main_corollary}]
Taking the expectation with respect to the sketching matrices $S_{t,k}$, $k=1,\dots,q$ of both sides of the equation given in Theorem \ref{thm:IHS_error_decay}, we obtain
\begin{align*}
    \Exs [\|e^A_{t+1}\|_2^2] = \frac{1}{q} \left(\frac{\theta_2}{\theta_1^2}-1\right) \Exs [\|e^A_t\|_2^2] \, .
\end{align*}
This gives us the relationship between the initial error (when we initialize $x_0$ to be the zero vector) and the expected error in iteration $t$:
\begin{align*}
    \Exs [\|e^A_t\|_2^2] = \frac{1}{q^t} \left(\frac{\theta_2}{\theta_1^2}-1\right)^t \|Ax^*\|_2^2.
\end{align*}
It follows that the expected error reaches $\epsilon$-accuracy with respect to the initial error at iteration $T$ where:
\begin{align*}
    \frac{1}{q^T} \left(\frac{\theta_2}{\theta_1^2}-1\right)^T &= \epsilon \\
    q^T \left(\frac{\theta_2}{\theta_1^2}-1\right)^{-T} &= \frac{1}{\epsilon} \\
    T\left(\log(q)-\log\left(\frac{\theta_2}{\theta_1^2}-1\right) \right) &= \log(1/\epsilon) \\
    T &= \frac{\log(1/\epsilon)}{\log(q)-\log\left(\frac{\theta_2}{\theta_1^2}-1\right)}.
\end{align*}

Each iteration requires communicating a $d$-dimensional vector for every worker, and we have $q$ workers and the algorithm runs for $T$ iterations, hence the communication load is $Tqd$.

The computational load per worker at each iteration involves the following numbers of operations:
\begin{itemize}
    \item Sketching $A$: $mnd$ multiplications
    \item Computing $\tilde{H}_{t,k}$: $md^2$ multiplications
    \item Computing $g_t$: $O(nd)$ operations
    \item Solving $\tilde{H}_{t,k}^{-1} g_t$: $O(d^3)$ operations.
\end{itemize}
\end{IEEEproof}

\begin{IEEEproof} [Proof of Lemma \ref{expectation_inverse_regularization_2}]
In the following, we assume that we are in the regime where $n$ approaches infinity. 

The expectation term $\Exs [ (U^TS^TSU+\lambda_2 I)^{-1}]$ is equal to the identity matrix times a scalar (i.e. $cI_d$) because it is signed permutation invariant, which we show as follows. Let $P \in \mathbb{R}^{d\times d}$ be a permutation matrix and $D \in \mathbb{R}^{d\times d}$ be an invertible diagonal sign matrix ($-1$ and $+1$'s on the diagonals). A matrix $M$ is signed permutation invariant if $(DP) M (DP)^T = M$. We note that the signed permutation matrix is orthogonal: $(DP)^T(DP) = P^TD^TDP = P^TP = I_d$, which we later use in the sequel. 
\begin{align*}
    (DP)\Exs_S [ (U^TS^TSU+\lambda_2 I)^{-1} ](DP)^T &= \Exs_S [ (DP)(U^TS^TSU+\lambda_2 I)^{-1} (DP)^T ] \\
    &= \Exs_S [ ((DP)^T U^TS^TSU(DP)+\lambda_2 I)^{-1} ] \\
    &= \Exs_{SUPD} [ \Exs_S [ ((DP)^T U^TS^TSU(DP)+\lambda_2 I)^{-1} | SUPD ]] \\
    &= \Exs_{SUPD} [ ((DP)^T U^TS^TSU(DP)+\lambda_2 I)^{-1} ] \\
    &= \Exs_{SU^\prime} [ ({U^\prime}^T S^TSU^\prime +\lambda_2 I)^{-1} ]
\end{align*}
where we made the variable change $U^\prime = UDP$ and note that $U^\prime$ has orthonormal columns because $DP$ is an orthogonal transformation. $SUPD$ and $SU$ have the same distribution because $PD$ is an orthogonal transformation and $S$ is a Gaussian matrix. This shows that $\Exs [ (U^TS^TSU+\lambda_2 I)^{-1}]$ is signed permutation invariant.

Now that we established that $\Exs [ (U^TS^TSU+\lambda_2 I)^{-1}]$ is equal to the identity matrix times a scalar, we move on to find the value of the scalar. We use the identity $\Exs_{DP}[(DP) Q (DP)^T] = \frac{\tr Q}{d}I_d$ for $Q\in \mathbb{R}^{d\times d}$ where the diagonal entries of $D$ are sampled from the Rademacher distribution and $P$ is sampled uniformly from the set of all possible permutation matrices. We already established that $\Exs [ (U^TS^TSU+\lambda_2 I)^{-1}]$ is equal to $(DP)\Exs_S [ (U^TS^TSU+\lambda_2 I)^{-1} ](DP)^T$ for any signed permutation matrix of the form $DP$. It follows that
\begin{align*}
    \Exs [ (U^TS^TSU+\lambda_2 I)^{-1}] &= (DP)\Exs_S [ (U^TS^TSU+\lambda_2 I)^{-1} ](DP)^T \\
    &= \frac{1}{|R|}\sum_{DP \in R} (DP)\Exs_S [ (U^TS^TSU+\lambda_2 I)^{-1} ](DP)^T \\
    &= \Exs_{DP} [ (DP)\Exs_S [ (U^TS^TSU+\lambda_2 I)^{-1} ](DP)^T ] \\
    &= \frac{1}{d} \tr (\Exs_S [ (U^TS^TSU+\lambda_2 I)^{-1} ]) I_d 
\end{align*}
where we define $R$ to be the set of all possible signed permutation matrices $DP$ in going from line 1 to line 2. 

By Lemma \ref{expectation_inverse_regularization}, the trace term is equal to $d \times \theta_3(d/m, \lambda_2)$, which concludes the proof.
\end{IEEEproof}

\begin{IEEEproof} [Proof of Theorem \ref{thm:opt_lambda_2_newton_LS}]
Closed form expressions for the optimal solution and the output of the $k$'th worker are as follows:
\begin{align*}
    x^* &= (A^TA + \lambda_1 I_d)^{-1} A^Tb, \\
    \hat{x}_k &= (A^TS_k^TS_kA + \lambda_2 I_d)^{-1}A^TS_k^TS_kb.
\end{align*}
Equivalently, $x^*$ can be written as:
\begin{align*}
    x^* = \arg\min \left\Vert \begin{bmatrix} A \\ \sqrt{\lambda_1}I_d \end{bmatrix} x - \begin{bmatrix} b \\ 0_d \end{bmatrix} \right\Vert_2^2.
\end{align*}
This allows us to decompose $\begin{bmatrix} b \\ 0_d \end{bmatrix}$ as 
\begin{align*}
    \begin{bmatrix} b \\ 0_d \end{bmatrix} = \begin{bmatrix} A \\ \sqrt{\lambda_1}I_d \end{bmatrix} x^* + b^\perp
\end{align*}
where $b^\perp = \begin{bmatrix} b_1^\perp \\ b_2^\perp \end{bmatrix}$ with $b_1^\perp \in \mathbb{R}^{n}$ and $b_2^\perp \in \mathbb{R}^d$. From the above equation we obtain $b^\perp_2 = - \sqrt{\lambda_1}x^*$ and $\begin{bmatrix} A^T & \sqrt{\lambda_1}I_d \end{bmatrix} b^\perp = A^Tb_1^\perp +\sqrt{\lambda_1}b_2^\perp = 0$. 

The bias of $\hat{x}_k$ is given by (omitting the subscript $k$ in $S_k$ for simplicity)
\begin{align} 
    \Exs[ A (\hat{x}_k - x^*)] &= \Exs[A(A^TS^TSA + \lambda_2 I_d)^{-1}A^TS^TSb - Ax^*] \nonumber\\
    %&= \Exs[ U(U^TS^TSU+\lambda_2 \Sigma^{-2})^{-1}U^TS^TSb] - Ax^* \nonumber\\
    &= \Exs[ U(U^TS^TSU+\lambda_2 \Sigma^{-2})^{-1}U^TS^TS(Ax^*+b_1^\perp)] - Ax^* \nonumber\\
    %&= \Exs[ U(U^TS^TSU+\lambda_2 \Sigma^{-2})^{-1}U^TS^TSU\Sigma V^Tx^* + U(U^TS^TSU+\lambda_2 \Sigma^{-2})^{-1}U^TS^TSb_1^\perp] - Ax^* \nonumber\\
    % &= \Exs[ U(U^TS^TSU+\lambda_2 \Sigma^{-2})^{-1}(U^TS^TSU+\lambda_2 \Sigma^{-2}-\lambda_2 \Sigma^{-2})\Sigma V^Tx^* \nonumber \\
    % &\quad + U(U^TS^TSU+\lambda_2 \Sigma^{-2})^{-1}U^TS^TSb_1^\perp] - Ax^* \nonumber\\
    &= \Exs[ - \lambda_2 U(U^TS^TSU+\lambda_2 \Sigma^{-2})^{-1}\Sigma^{-1}V^Tx^*] + \Exs[U(U^TS^TSU+\lambda_2 \Sigma^{-2})^{-1}U^TS^TSb_1^\perp]\,. \nonumber
\end{align}
By the assumption $\Sigma=\sigma I_d$, the bias becomes
\begin{align} \label{eq:bias_expect_decomp}
    \Exs &[ A (\hat{x}_k - x^*)] = \Exs[ - \lambda_2 \sigma^{-1} U(U^TS^TSU+\lambda_2 \sigma^{-2} I_d)^{-1}V^Tx^*] + \Exs[U(U^TS^TSU+\lambda_2 \sigma^{-2} I_d)^{-1}U^TS^TSb_1^\perp]\, .
\end{align}
The first expectation term of \eqref{eq:bias_expect_decomp} can be evaluated using Lemma \ref{expectation_inverse_regularization_2} (as $n$ goes to infinity):
\begin{align} \label{eq:first_expect_term}
    \Exs [ - \lambda_2 \sigma^{-1} U&(U^TS^TSU+\lambda_2 I_d)^{-1}V^Tx^*] = -\lambda_2 \sigma^{-1} \theta_3(d/m, \lambda_2 \sigma^{-2}) UV^Tx^*.
\end{align}
To find the second expectation term in \eqref{eq:bias_expect_decomp}, let us first consider the full SVD of $A$ given by $A = \begin{bmatrix} U & U^\perp \end{bmatrix} \begin{bmatrix} \Sigma \\ 0_{(n-d) \times d} \end{bmatrix} V^T$ where $U \in \mathbb{R}^{n\times d}$ and $U^\perp \in \mathbb{R}^{n\times (n-d)}$. The matrix $\begin{bmatrix} U & U^\perp \end{bmatrix}$ is an orthogonal matrix, which implies $UU^T + U^\perp (U^\perp)^T = I_d$. If we insert $UU^T + U^\perp (U^\perp)^T = I_d$ between $S$ and $b_1^\perp$, the second term of \eqref{eq:bias_expect_decomp} becomes
\begin{align*}
    &\Exs[U(U^TS^TSU+\lambda_2 \sigma^{-2} I_d)^{-1}U^TS^TSb_1^\perp] = \\
    %&=\Exs[U(U^TS^TSU+\lambda_2 \sigma^{-2} I_d)^{-1}U^TS^TS(UU^T + U^\perp (U^\perp)^T)b_1^\perp] \\
    &= \Exs[U(U^TS^TSU+\lambda_2 \sigma^{-2} I_d)^{-1}U^TS^TS UU^T b_1^\perp] + \Exs[U(U^TS^TSU+\lambda_2 \sigma^{-2} I_d)^{-1}U^TS^TS U^\perp (U^\perp)^T b_1^\perp] \\
    &= \Exs[U(U^TS^TSU+\lambda_2 \sigma^{-2} I_d)^{-1}U^TS^TS UU^T b_1^\perp] \\
    % &= \Exs[U(U^TS^TSU+\lambda_2 \sigma^{-2} I_d)^{-1}(U^TS^TS U + \lambda_2 \sigma^{-2} I_d - \lambda_2 \sigma^{-2} I_d) U^T b_1^\perp] \\
    &= U (I_d - \lambda_2 \sigma^{-2} \Exs [ (U^TS^TSU+\lambda_2 \sigma^{-2} I_d)^{-1} ]) U^T b_1^\perp \\
    &= (1 - \lambda_2 \sigma^{-2} \theta_3(d/m, \lambda_2\sigma^{-2})) UU^T b_1^\perp \,.
\end{align*}
In these derivations, we used $\Exs_S [ U(U^TS^TSU+\lambda_2 \sigma^{-2} I_d)^{-1}U^TS^TS U^\perp (U^\perp)^T b_1^\perp ] =\Exs_{SU} [\Exs_S  [ U(U^TS^TSU+\lambda_2 \sigma^{-2} I_d)^{-1}U^TS^TS U^\perp (U^\perp)^T b_1^\perp | SU]] = 0$. This follows from $\Exs_S[SU^\perp|SU] = 0$ as $U$ and $U^\perp$ are orthogonal. The last line follows from Lemma \ref{expectation_inverse_regularization_2}, as $n$ goes to infinity.

We note that $U^T b_1^\perp = \lambda_1 \Sigma^{-1}V^Tx^*$ and for $\Sigma = \sigma I_d$, this becomes $U^T b_1^\perp = \lambda_1 \sigma^{-1} V^Tx^*$.
Bringing the pieces together, we have the bias equal to (as $n$ goes to infinity):
\begin{align*}
    \Exs [ A (\hat{x}_k - x^*)] &= -\lambda_2 \sigma^{-1} \theta_3(d/m, \lambda_2 \sigma^{-2}) UV^Tx^*  + \lambda_1 \sigma^{-1} (1 - \lambda_2 \sigma^{-2} \theta_3(d/m, \lambda_2\sigma^{-2})) UV^Tx^* \\
    &= \sigma^{-1} (\lambda_1 - \lambda_2 \theta_3(d/m, \lambda_2 \sigma^{-2}) ( 1 + \lambda_1 \sigma^{-2})) UV^Tx^*.
\end{align*}
If there is a value of $\lambda_2 > 0$ that satisfies $\lambda_1 - \lambda_2 \theta_3(d/m, \lambda_2 \sigma^{-2}) ( 1 + \lambda_1 \sigma^{-2}) = 0$, then that value of $\lambda_2$ makes $\hat{x}_k$ an unbiased estimator. Equivalently,
\begin{align} \label{LHS_equation}
    %\left( \frac{-\lambda_2\sigma^{-2}+d/m-1 + \sqrt{(-\lambda_2\sigma^{-2} + d/m-1)^2+4\lambda_2\sigma^{-2}d/m}}{2\sigma^{-2} d/m} \right) &= \frac{\lambda_1}{1+\lambda_1\sigma^{-2}} \\
    -\lambda_2\sigma^{-2}+\frac{d}{m}-1 + &\sqrt{(-\lambda_2\sigma^{-2} + \frac{d}{m}-1)^2+4\lambda_2\sigma^{-2}\frac{d}{m}} = 2 \frac{d}{m \sigma^2} \frac{\lambda_1}{1+\lambda_1\sigma^{-2}},
\end{align}
where we note that the LHS is a monotonically increasing function of $\lambda_2$ in the regime $\lambda_2 \geq 0$ and it attains its minimum in this regime at $\lambda_2=0$. Analyzing this equation using these observations, for the cases of $m > d$ and $m \leq d$ separately, we find that for the case of $m \leq d$, we need the following to be satisfied for zero bias:
\begin{align*}
    2\frac{d}{m\sigma^2}\frac{\lambda_1}{1+\lambda_1/\sigma^2} &\geq 2\left(\frac{d}{m}-1\right),\, \mbox{or more simply,} \\
    \lambda_1 &\geq \sigma^2 \left( \frac{d}{m} - 1\right) ,
\end{align*}
whereas there is no additional condition on $\lambda_1$ for the case of $m > d$.

The value of $\lambda_2$ that will lead to zero bias can be computed by solving the equation \eqref{LHS_equation} where the expression for the inverse of the left-hand side is given by $LHS^{-1}(y)=\frac{y\sigma^{-2}d/m - d/m+1}{y^{-1}-\sigma^{-2}}$. We evaluate the inverse at $y=\lambda_1 / (1+\lambda_1/\sigma^2)$ and obtain the following expression for $\lambda_2^*$: 
\begin{align*}
    \lambda_2^* = \lambda_1 - \frac{d}{m}\frac{\lambda_1}{1+\lambda_1/\sigma^2}.
\end{align*}
\end{IEEEproof}

%%%%%%%%%%%%%%%%%%%%%%%%%%%

\subsection{Proofs of Theorems and Lemmas in Section \ref{sec:privacy}}
The proof of the privacy result mainly follows due to Theorem \ref{thm:diff_priv_cited} (stated below for completeness) which is a result from \cite{sheffet15privacy_SM}.

\begin{IEEEproof} [Proof of Theorem \ref{thm:diff_priv_result}]
For some $\varepsilon$, $\delta$ and matrix $A_c$, if there exist values for $m$ such that the smallest singular value of $A_c$ satisfies $\sigma_{\min}(A_c) \geq w$, then using Theorem \ref{thm:diff_priv_cited}, we find that the sketch size $m$ has to satisfy the following for $(\varepsilon, \delta)$-differential privacy: 
\begin{align} \label{eq:sketch_size_privacy}
    m &\leq  \frac{1}{8\ln(4/\delta)} \left(\left(\frac{\sigma_{\min}^2}{B^2} - 1\right) \frac{1}{\frac{1}{\varepsilon}+\frac{1}{\ln(4/\delta)}} - 2\ln(4/\delta) \right)^2 \nonumber \\
    &= \frac{1}{8\beta} \left(\left(\frac{\sigma_{\min}^2}{B^2} - 1\right) \frac{\varepsilon \beta}{\varepsilon + \beta} - 2\beta \right)^2,
\end{align}
where we have set $\delta = 4/e^\beta$ in the second line. For the first line to follow from Theorem \ref{thm:diff_priv_cited}, we also need the condition $\frac{\sigma_{\min}^2}{B^2} \geq 3 + 2\frac{\beta}{\varepsilon}$ to be satisfied. Note that the rows of $A_c$ have bounded $\ell_2$-norm of $B_0\sqrt{d+1}$. We now substitute $B=B_0 \sqrt{d+1}$ and $\sigma_{\min} = \sigma_0 \sqrt{n}$ to obtain the simplified condition
\begin{align*}
    \frac{n}{d+1} \geq \left(3+2 \frac{\beta}{\varepsilon} \right) \frac{B_0^2}{\sigma_0^2},
\end{align*}
where $B_0$ and $\sigma_0$ are constants. Assuming this condition is satisfied, then we pick the sketch size $m$ as \eqref{eq:sketch_size_privacy} which can also be simplified:
\begin{align*}
    m = O\left(\beta \frac{n^2}{(d+1)^2} \frac{\varepsilon^2}{(\varepsilon + \beta)^2} \right).
\end{align*}
Note that the above arguments are for the privacy of a single sketch (i.e., $S_kA_c$). In the distributed setting where the adversary can attack all of the sketched data $S_1A_c, \dots, S_qA_c$, we can consider all of the sketched data to be a single sketch with size $mq$. Based on this argument, we can pick the sketch size as
\begin{align*}
    m = O\left(\frac{\beta}{q} \frac{n^2}{(d+1)^2} \frac{\varepsilon^2}{(\varepsilon + \beta)^2} \right).
\end{align*}
\end{IEEEproof}

% we should have d+1 instead of d
%%%%%%

% If we fix $\delta = 4/e^\beta$, then $\ln(4/\delta) = \beta$. Then, $m$ can be picked:
% \begin{align*}
%     m &\leq \frac{1}{q} \frac{1}{24} \left(\left(\frac{\sigma_{\min}^2}{B^2} - 1\right) \frac{3\varepsilon}{3+\varepsilon} - 6 \right)^2 \\
%     &= O\left(\frac{\varepsilon^2}{q} \right)
% \end{align*}

% Picking $m = O\left(\frac{\varepsilon^2}{q} \right)$ and assuming $d << m$, we obtain
% \begin{align*}
%     \frac{\Exs[f(\bar{x})]-f(x^*)}{f(x^*)} &= \frac{1}{q} \frac{d}{m-d-1} \\
%     &= O\left(\frac{1}{\varepsilon^2}\right)
% \end{align*}
%\qed

% Theorem 3.1. of https://arxiv.org/pdf/1507.00056v2.pdf
\begin{theorem}[Differential privacy for random projections \cite{sheffet15privacy_SM}] \label{thm:diff_priv_cited}
Fix $\varepsilon > 0$ and $\delta \in (0, 1/e)$. Fix $B > 0$. Fix a positive integer $m$ and let $w$ be such that
\begin{align} \label{eq:def_w}
    w^2 = B^2 \left( 1 + \frac{1+\frac{\varepsilon}{\ln(4/\delta)}}{\varepsilon} \left( 2 \sqrt{2m \ln(4/\delta)} + 2\ln(4/\delta) \right)\right).
\end{align}
Let $A$ be an $(n\times d)$-matrix with $d < m$ and where each row of $A$ has bounded $\ell_2$-norm of $B$. Given that $\sigma_{\min}(A) \geq w$, the algorithm that picks an $(m \times n)$-matrix $R$ whose entries are iid samples from the normal distribution $\mathcal{N}(0,1)$ and publishes the projection $RA$ is $(\varepsilon, \delta)$-differentially private.
\end{theorem}

% We note that Remark 1 is a direct consequence of of Theorem \ref{thm:diff_priv_result} and Theorem \ref{Thm:AvgGaussian}.

%%%%%%%

%%%%%%%%%%%%%%%%%%%%%%%%%
\subsection{Proofs of Theorems and Lemmas in Section \ref{sec:other_sketching_matrices}}

\begin{IEEEproof} [Proof of Lemma \ref{expected_obj_val_diff}]
The expectation of the difference between the costs $f(\bar{x})$ and $f(x^*)$ is given by
\begin{align} \label{eq:expected_diff}
    \Exs[f(\bar{x})]-f(x^*)% &= \Exs[\|A\bar{x}-b\|_2^2]-f(x^*) \\
    &= \Exs[\|A(\bar{x}-x^*)+Ax^*-b\|_2^2]-f(x^*) \nonumber \\
    & = \Exs [\|A(\bar{x}-x^*)\|_2^2 + \|Ax^*-b\|_2^2]-f(x^*) \nonumber \\
    & = \Exs [\|A(\bar{x}-x^*)\|_2^2]\,,
\end{align}
where we have used the orthogonality property of the optimal least squares solution $x^*$ given by the normal equations $A^T (Ax^*-b) = 0$. Next, we have
\begin{align*}
   \Exs[\|A(\bar{x}-x^*)\|_2^2] &= \Exs\left[ \left \|\frac{1}{q}\sum_{k=1}^q (A\hat{x}_k - Ax^*) \right \|_2^2 \right] \\
    &= \frac{1}{q^2} \Exs \left[ \sum_{k=1}^q\sum_{l=1}^q \langle A\hat{x}_k - Ax^*, A\hat{x}_l-Ax^* \rangle \right] \\
    &= \frac{1}{q^2} \sum_{k=1}^q \Exs \left[ \|A\hat{x}_k-Ax^* \|_2^2 \right]  + \frac{1}{q^2} \sum_{k\neq l\,, 1 \leq k,l \leq q} \Exs \left[ \langle A\hat{x}_k - Ax^*, A\hat{x}_l-Ax^* \rangle \right] \\
    %&= \frac{1}{q} \Exs \left[ \|A\hat{x}_1-Ax^* \|_2^2 \right] + \frac{q^2-q}{q^2} \Exs \left[ \langle A\hat{x}_1 - Ax^*, A\hat{x}_2-Ax^* \rangle \right] \\
    %&= \frac{1}{q} \Exs \left[ \|A\hat{x}_1-Ax^* \|_2^2 \right] + \frac{q-1}{q} \Exs[A\hat{x}_1 - Ax^*]^T \Exs[A\hat{x}_2 - Ax^*] \\
    &= \frac{1}{q} \Exs \left[ \|A\hat{x}-Ax^* \|_2^2 \right] + \frac{q-1}{q} \| \Exs[A\hat{x}] - Ax^* \|_2^2 \,.
\end{align*}
\end{IEEEproof}

\begin{IEEEproof} [Proof of Lemma \ref{norm_of_bias}]
The bias of the single sketch estimator can be expanded as follows:
\begin{align*}
    \| \Exs [A\hat{x}]-Ax^* \|_2 &= \|\Exs[A(A^TS^TSA)^{-1}A^TS^TS(Ax^*+b^\perp)] - Ax^* \|_2 \\
    % &= \|\Exs[A(A^TS^TSA)^{-1}A^TS^TS b^\perp ] \|_2 \\
    &= \|U \Exs[(U^TS^TSU)^{-1}U^TS^TSb^\perp \|_2 \\
    &= \| \Exs[(U^TS^TSU)^{-1}U^TS^TSb^\perp \|_2 \\
    &= \|\Exs[Qz] \|_2,
\end{align*}
where we define $Q := (U^TS^TSU)^{-1} $ and $z := U^TS^TSb^\perp$. The term $\|\Exs[Qz] \|_2^2$ can be upper bounded as follows when conditioned on the event $E$:
\begin{align*}
    \|\Exs[Qz] \|_2^2 &= \Exs[Qz]^T \Exs[Qz] = \Exs_{S}[Qz]^T \Exs_{S^\prime}[Q^\prime z^\prime] \\
    &= \Exs_{S_k} \Exs_{S_k^\prime} [z^T Q Q^\prime z^\prime ] \\
    &= \frac{1}{2} \Exs_{S} \Exs_{S^\prime} [(z+z^\prime)^TQQ^\prime (z+z^\prime) - z^TQQ^\prime z - {z^\prime}^T QQ^\prime z^\prime] \\
    &\leq \frac{1}{2} \Exs_{S} \Exs_{S^\prime} [ \|z+z^\prime\|_2^2 (1+\epsilon)^2 - (\|z\|_2^2 + \|z^\prime \|_2^2 )(1-\epsilon)^2  ] \\
    &= \Exs_{S} \Exs_{S^\prime} [ \left( \|z\|_2^2 2\epsilon + \|z^\prime\|_2^2 2\epsilon + z^Tz^\prime (1+\epsilon)^2 \right)] \\
    %&= 2\epsilon \Exs_{S_k} \Exs_{S_k^\prime} [\|z\|_2^2 ]  +  2\epsilon \Exs_{S_k} \Exs_{S_k^\prime} [\|z^\prime \|_2^2]  +  (1+\epsilon)^2 \Exs_{S_k} \Exs_{S_k^\prime} [z^Tz^\prime ] \\
    &= 4\epsilon \Exs [\|z\|_2^2 ]  +  (1+\epsilon)^2 \|\Exs [z ]\|_2^2 \,,
\end{align*}
where the inequality follows from the inequality $(1-\epsilon)I_d \preceq Q \preceq (1+\epsilon)I_d$ and some simple bounds for the minimum and maximum eigenvalues of the product of two positive definite matrices. Furthermore, the expectation of $z$ is equal to zero because $\Exs[z] = \Exs[U^TS^TSb^\perp]=U^T\Exs[S^TS]b^\perp=U^Tb^\perp=0$. Hence we obtain the claimed bound $
    \| \Exs [A\hat{x}|E]-Ax^* \|_2 \leq \sqrt{4\epsilon \Exs[\|z\|_2^2|E]}. $
% By Lemma \ref{expectation_upper_bound} and Lemma \ref{bound_z_times_indicator}, we have:
% \begin{align*}
%     \|\Exs[Qz] \|_2^2 &\leq 4\epsilon \mprob[E] \Exs_S [\|z\|_2^2]  +  (1+\epsilon)^2 \mprob[E]^2 \Exs [\|z\|_2^2] \\
%     &= (4\epsilon \mprob[E] + (1+\epsilon)^2 \mprob[E]^2) \Exs_S [\|z\|_2^2] \\
%     %&\leq (4\epsilon + (1+\epsilon)^2) \Exs_S [\|z\|_2^2] \\
%     &\leq (1+7\epsilon) \Exs_S [\|z\|_2^2]
% \end{align*}
% Finally, we obtain the upper bound on the norm of the bias as follows:
% \begin{align*}
%     \| \Exs [A\tilde{x}]-Ax^* \|_2 &\leq  \sqrt{(1+7\epsilon) \Exs_S [\|z\|_2^2]} + \|Ax^*\|_2 \mprob[E^c]
% \end{align*}
\end{IEEEproof}

\begin{IEEEproof} [Proof of Lemma \ref{bound_z_norm_ROS}]
For the randomized Hadamard sketch (ROS), the term $\Exs[ \| z\|_2^2]$ can be expanded as follows. We will assume that all the expectations are conditioned on the event $E$, which we defined earlier as $(1-\epsilon)I_d \preceq Q \preceq (1+\epsilon) I_d$.
\begin{align*}
\Exs [\|z\|_2^2] &= \Exs \left[ {b^\perp}^T \frac{1}{m}\sum_{i=1}^m s_i s_i^T UU^T \frac{1}{m}\sum_{j=1}^m s_j s_j^T b^{\perp} \right] \\
&= \Exs \left[ \frac{1}{m^2} \sum_{i=1}^m \sum_{j=1}^m {b^\perp}^T s_i s_i^T UU^T s_j s_j^T b^{\perp} \right]\\
&= \frac{1}{m^2} \sum_{1\le i=j\le m} {b^\perp}^T \Exs \left[ s_i s_i^T UU^T s_j s_j^T\right] b^{\perp} + \frac{1}{m^2} \sum_{i\neq j,\,1\le i,j\le m}  {b^\perp}^T \Exs [s_i s_i^T] UU^T \Exs [s_j s_j^T] b^{\perp}\\
&= \frac{m}{m^2} \, {b^\perp}^T \Exs \left[ s_1 s_1^T UU^T s_1 s_1^T \right] b^{\perp} + \frac{1}{m^2} \sum_{i\neq j,\,1\le i,j\le m}  {b^\perp}^T I_n UU^T I_n b^{\perp} \\
&= \frac{1}{m} \, {b^\perp}^T \Exs \left[ s_1 s_1^T UU^T s_1 s_1^T \right] b^{\perp},
\end{align*}
where we have used the independence of $s_i$ and $s_j$, $i \neq j$. This is true because of the assumption that the matrix $P$ corresponds to sampling with replacement.
\begin{align*}
    {b^\perp}^T \Exs \left[ s_1 s_1^T UU^T s_1 s_1^T \right] b^{\perp} &= \Exs [ (s_1^TUU^Ts_1) (s_1^Tb^\perp {b^\perp}^Ts_1) ] \\
    &= \Exs [ (s_1^TUU^Ts_1) ({b^\perp}^Ts_1)^2 ] \\
    &= \frac{1}{n} \sum_{i=1}^n \Exs[(h_i^TDUU^TDh_i)({b^\perp}^TDh_i)^2],
\end{align*}
where the row vector $h_i^T$ corresponds to the $i$'th row of the Hadamard matrix $H$. We also note that the expectation in the last line is with respect to the randomness of $D$.

Let us define $r$ to be the column vector containing the diagonal entries of the diagonal matrix $D$, that is, $r := [D_{11}, D_{22}, \dots, D_{nn}]^T$. Then, the vector $Dh_i$ is equivalent to $\diag(h_i)r$ where $\diag(h_i)$ is the diagonal matrix with the entries of $h_i$ on its diagonal.
\begin{align*}
    \frac{1}{n} \sum_{i=1}^n \Exs[(h_i^TDUU^TDh_i)({b^\perp}^TDh_i)^2] &= \frac{1}{n} \sum_{i=1}^n \Exs[(r^T\diag(h_i)UU^T\diag(h_i)r)({b^\perp}^T\diag(h_i)r)^2] \\
    &= \frac{1}{n} \sum_{i=1}^n \Exs[{b^\perp}^T \diag(h_i) r (r^TPr)  r^T \diag(h_i) b^\perp] \\
    &= \frac{1}{n} \sum_{i=1}^n {b^\perp}^T \diag(h_i) \Exs[ r (r^TPr) r^T ]\diag(h_i) b^\perp ,
\end{align*}
where we have defined $P := \diag(h_i) UU^T \diag(h_i)$. It follows that $\Exs[ r (r^TPr) r^T ] = 2P - 2 \diag(P) + \tr( P) I_n$. Here, $\diag(P)$ is used to refer to the diagonal matrix with the diagonal entries of $P$ as its diagonal. 
The trace of $P$ can be easily computed using the cyclic property of matrix trace as $\tr (P) = \tr (\diag(h_i) UU^T \diag(h_i)) = \tr(U^T \diag(h_i) \diag(h_i) U) = \tr(U^TU) = \tr(I_d) = d$.
%
% -------
% \begin{align*}
% \mbox{\red{ (*) this step can be improved}}
% &\le\frac{1}{m} \, {b^\perp}^T \frac{1}{n} \sum_{i=1}^n s_i  s_i^T b^{\perp} \max_{i=1,\dots,n} (s_i^T UU^T s_i) \\
% &=\frac{1}{m} \, {b^\perp}^T b^{\perp} \max_{i=1,\dots,n} (s_i^T UU^T s_i) \\
% &=\frac{1}{m} \| {b^\perp}\|_2^2 \max_{i=1,\dotsn} \|U^TDH p_i\|_2^2
% \end{align*}
% %
% \mbox{\red{(*)}}'s expectation can be computed as
% %
% \begin{align}
%     \sum_i \Exs ({b^\perp}' D h_i)^2 (h_i^T D  UU^T D h_i)  &= \sum_i \Exs   ({b^\perp}' diag(h_i) r)^2 (r diag(h_i)  UU^T diag(h_i) r \\
%     &= \sum_i \Exs \sum_{ijkl} r_i r_j r_k r_l a_i a_j q_{ij}\\
%     & = \sum_i \Exs \sum_{j} a_j a_j q_{jj} + \Exs \sum_{kl} a_k a_k q_{ll} + \dots.
% \end{align}
% %
% %
% where $r = vec(D)$, $a_i = b^\perp diag(h)$, $q_ij=diag(h_i) UU^T diag(h_i)$.\\
% %
% Another way to bound the term in \mbox{\red{(*)}}: Let $r = vec(D)$ i.i.d. $\pm 1$. $\Exs r r^T r^T Q r = I \trace Q  + Q + Q^T - 2diag(Q)$ for any matrix Q. We set $Q=Diag(h_i) UU^T Diag(h_i)$, which satisfies $\trace Q = \trace U^T Diag(h_i)Diag(h_i) U = \trace U^T U = \trace I_d = d$ using the cyclic property of matrix trace.
%
Next, we note that the term $\diag(P)$ can be simplified as $\diag(P)_{jj}=\|\tilde{u}_j \|_2^2$ where $\tilde{u}_j^T$ is the $j$'th row of $U$. This leads to
\begin{align*}
    {b^\perp}^T \diag(P) b^\perp &= \sum_{j=1}^n (b_j^\perp)^2 \|\tilde{u}_j \|_2^2 \geq \sum_{j=1}^n (b_j^\perp)^2 \min_i \|\tilde{u}_i \|_2^2 = \|b^\perp \|_2^2 \min_i \|\tilde{u}_i \|_2^2 \,.
\end{align*}
Going back to $\Exs [\|z \|_2^2]$, we have
\begin{align*}
    \Exs [\|z \|_2^2] &= \frac{1}{mn} {b^\perp}^T  n \diag(h_i) \left(2P-2\diag(P) + \tr(P)I_n\right)\diag(h_i) b^\perp \\
    &= \frac{d}{m} \|b^\perp \|_2^2 - \frac{2}{m} {b^\perp}^T \diag(P) b^\perp \\
    &\leq \frac{d}{m} \|b^\perp \|_2^2 - \frac{2}{m} \|b^\perp \|_2^2 \min_i \|\tilde{u}_i \|_2^2 \\
    &= \frac{1}{m} \|b^\perp \|_2^2 (d-2\min_i \|\tilde{u}_i \|_2^2) = \frac{d}{m} \left(1-\frac{2\min_i \|\tilde{u}_i \|_2^2}{d} \right) f(x^*).
\end{align*}
%
%where $m=$ \todo{what exactly do we set m to?}
%\todo{we can set $m=d \log(2n)/\epsilon)$}.
%
% Next taking expectations over the sign variables $D$ we obtain
% %
% \begin{align*}
% \Exs_D \Exs \|z\|_2^2 &\le \frac{1}{m} \| {b^\perp}\|_2^2 \rank(A) \log(2n)\\
% &\le \epsilon f(x^*).
% \end{align*}
% %
% The above inequality follows from the fact that $\sum_{j} H_{ij}D_{jj} U_{j}$ is sub-Gaussian with parameter $\|U\|_F^2 = \trace UU^T = rank(A)$.
% \qed
\end{IEEEproof}

\begin{IEEEproof} [Proof of Lemma \ref{bound_z_norm_uniform}]
We will assume that all the expectations are conditioned on the event $E$, which we defined earlier as $(1-\epsilon)I_d \preceq Q \preceq (1+\epsilon) I_d$.
% Note that $\Exs[S^TS] = \frac{1}{m} \sum_{i=1}^m s_i s_i^T = I_n$ where $s_i \in \mathbb{R}^n$ is the column vector corresponding to the $i$'th row of $S_k$ scaled by $1/\sqrt{m}$. Note that because of this particular scaling, $\Exs[s_is_i^T] = I_n$ holds.
For uniform sampling with replacement, we have
\begin{align*}
    \Exs [\|z\|_2^2] &= \frac{1}{m^2}\, \Exs\left[ {b^\perp}^T \sum_{i=1}^m s_i s_i^T UU^T \sum_{j=1}^m s_j s_j^T b^{\perp} \right] \\
    &= \frac{1}{m^2}\, \Exs\left[ \sum_{i=1}^m \sum_{j=1}^m {b^\perp}^T s_i s_i^T UU^T s_j s_j^T b^{\perp} \right] \\
    &= \frac{1}{m^2} \sum_{1\le i=j\le m} {b^\perp}^T \Exs [s_i s_i^T UU^T s_j s_j^T] b^{\perp} + \frac{1}{m^2}\sum_{i\neq j,\,1\le i,j\le m}  {b^\perp}^T \Exs [s_i s_i^T] UU^T \Exs [s_j s_j^T] b^{\perp} \\
    &= \frac{1}{m}\, {b^\perp}^T \Exs [s_1 s_1^T UU^T s_1 s_1^T] b^{\perp} + \frac{1}{m^2}\, \sum_{i\neq j,\,1\le i,j\le m}  {b^\perp}^T I_n UU^T I_n b^{\perp} \\
    &= \frac{1}{m}\, {b^\perp}^T \Exs [s_1 s_1^T UU^T s_1 s_1^T] b^{\perp} \\
    &= \frac{1}{m}\, {b^\perp}^T n^2 \frac{1}{n} \sum_{i=1}^n e_ie_i^TUU^Te_ie_i^T b^\perp \\
    % &= \frac{n}{m}\, {b^\perp}^T \sum_{i=1}^n \|\tilde{u}_i\|_2^2 e_i e_i^T b^\perp \\
    %&= \frac{n}{m}\, {b^\perp}^T Diag( \|\tilde{u}_i \|_2^2 ) b^\perp \\
    &= \frac{n}{m}\, \sum_{i=1}^n {b_i^\perp}^2 \|\tilde{u}_i \|_2^2 \\
    &\leq \frac{n}{m}\, \sum_{i=1}^n {b_i^\perp}^2 \max_j \|\tilde{u}_j \|_2^2 = \frac{\mu d}{m} f(x^*)  \,.
\end{align*}
Next, for uniform sampling without replacement, the rows $s_i$ and $s_j$ are not independent which can be seen by noting that given $s_i$, we know that $s_j$ will have its nonzero entry at a different place than $s_i$. Hence, differently from uniform sampling with replacement, the following term will not be zero:
\begin{align*}
    \frac{1}{m^2}\sum_{i\neq j,\,1\le i,j\le m}  {b^\perp}^T \Exs [s_i s_i^T UU^T s_j s_j^T] b^{\perp} &= \frac{m^2-m}{m^2} {b^\perp}^T \Exs[s_1s_1^TUU^Ts_2s_2^T] b^\perp \\
    &= \frac{m-1}{m} {b^\perp}^T n^2 \frac{1}{n^2-n} \sum_{i\neq j, 1 \leq i,j \leq n} e_ie_i^TUU^Te_je_j^Tb^\perp \\
    &= \frac{m-1}{m} \frac{n}{n-1} {b^\perp}^T \sum_{i\neq j, 1 \leq i,j \leq n} e_i\tilde{u}_i^T \tilde{u}_j e_j^Tb^\perp \\
    &= \frac{m-1}{m} \frac{n}{n-1} {b^\perp}^T (UU^T - \diag(\|\tilde{u}_i \|_2^2) b^\perp \\
    % &= \frac{m-1}{m} \frac{n}{n-1} (0 - \sum_{i=1}^n {b_i^\perp}^2 \|\tilde{u}_i \|_2^2) \\
    &= -\frac{m-1}{m} \frac{n}{n-1}\sum_{i=1}^n {b_i^\perp}^2 \|\tilde{u}_i \|_2^2 \,.
\end{align*}
It follows that for uniform sampling without replacement, we obtain
\begin{align*}
    \Exs[\|z\|_2^2] &= \left (\frac{n}{m} - \frac{m-1}{m} \frac{n}{n-1} \right) \sum_{i=1}^n {b_i^\perp}^2 \|\tilde{u}_i \|_2^2 \\
    &= \frac{n}{m} \frac{n-m}{n-1} \sum_{i=1}^n {b_i^\perp}^2 \|\tilde{u}_i \|_2^2 \\
    &\leq \frac{n}{m} \frac{n-m}{n-1} f(x^*) \max_i \|\tilde{u}_i \|_2^2 \\
    &= \frac{\mu d}{m} \frac{n-m}{n-1} f(x^*) \,.
\end{align*}
\end{IEEEproof}

\begin{IEEEproof} [Proof of Lemma \ref{bound_z_norm_leverage}]
We consider leverage score sampling with replacement. The rows $s_i$, $s_j$ $i \neq j$ are independent because sampling is with replacement. We will assume that all the expectations are conditioned on the event $E$, which we defined earlier as $(1-\epsilon)I_d \preceq Q \preceq (1+\epsilon) I_d$. For leverage score sampling, the term $\Exs [\|z\|_2^2]$ is upper bounded as follows:
\begin{align*}
    \Exs [\|z\|_2^2] &= \frac{1}{m^2}\, \Exs\left[ {b^\perp}^T \sum_{i=1}^m s_i s_i^T UU^T \sum_{j=1}^m s_j s_j^T b^{\perp} \right] \\
    &= \frac{1}{m^2}\, \Exs\left[ \sum_{i=1}^m \sum_{j=1}^m {b^\perp}^T s_i s_i^T UU^T s_j s_j^T b^{\perp} \right] \\
    &= \frac{1}{m^2} \sum_{1\le i=j\le m} {b^\perp}^T \Exs [s_i s_i^T UU^T s_j s_j^T] b^{\perp} + \frac{1}{m^2} \sum_{i\neq j,\,1\le i,j\le m}  {b^\perp}^T \Exs [s_i s_i^T] UU^T \Exs [s_j s_j^T] b^{\perp} \\
    % &= \frac{1}{m}\, {b^\perp}^T \Exs [s_1 s_1^T UU^T s_1 s_1^T] b^{\perp} + \frac{m^2-m}{m^2} {b^\perp}^T \Exs [s_1 s_1^T] UU^T \Exs [s_1 s_1^T] b^{\perp} \\
    &= \frac{1}{m}\, {b^\perp}^T \sum_{i=1}^n \frac{\ell_i}{d} \frac{d}{\ell_i}e_ie_i^TUU^T \frac{d}{\ell_i}e_ie_i^T b^\perp +  \frac{m^2-m}{m^2} {b^\perp}^T I_n UU^T I_n b^{\perp}\\
    &= \frac{1}{m}\, {b^\perp}^T \sum_{i=1}^n \frac{d}{\ell_i} \ell_i e_ie_i^T b^\perp  \\
    % &= \frac{d}{m}\, {b^\perp}^T \sum_{i=1}^n e_ie_i^T b^\perp  \\
    &= \frac{d}{m}\, \|b^\perp \|_2^2 = \frac{d}{m} f(x^*) \,.
\end{align*}
\end{IEEEproof}

\begin{IEEEproof} [Proof of Theorem \ref{thm:errorbounds_other_sketches}]
We will begin by using the results from Table 5 of \cite{mahoney2018averaging} for selecting the sketch size:
\begin{align} \label{eq:sketchsizes_other_sketches}
    m &= O\Big(\frac{d+\log(n)}{\epsilon^2} \log(d/\delta) \Big) \mbox{ for rand. Hadamard sketch} \nonumber \\
    m &= O\Big(\frac{\mu d}{\epsilon^2} \log(d/\delta) \Big) \mbox{ for uniform sampling} \nonumber \\
    m &= O\Big(\frac{d}{\epsilon^2} \log(d/\delta) \Big) \mbox{ for leverage score sampling}
\end{align}
where $\mu$ is the row coherence of $U$ as defined before. When the sketch sizes are selected according to the formulas above, the subspace embedding property given by
\begin{align}
    \|U^TS^TSU - I_d\|_2 \leq \epsilon
\end{align}
is satisfied with probability at least $1-\delta$.
Note that the subspace embedding property can be rewritten as
\begin{align}
    (1-\epsilon)I_d \preceq U^TS^TSU \preceq (1+\epsilon)I_d \,.
\end{align}
This implies the following relation for the inverse matrix $(U^TS^TSU)^{-1}$:
\begin{align}
    \frac{1}{1+\epsilon}I_d \preceq (U^TS^TSU)^{-1} \preceq \frac{1}{1-\epsilon}I_d \,.
\end{align}
Observe that $1-\epsilon \leq \frac{1}{1+\epsilon}$ for any $\epsilon$ and that $\frac{1}{1-\epsilon} \leq 1+2\epsilon$ for $0 \leq \epsilon \leq 0.5$. Assuming that $0 \leq \epsilon \leq 0.5$, we have
\begin{align}
    (1-2\epsilon) I_d \preceq (U^TS^TSU)^{-1} \preceq (1+2\epsilon) I_d \,.
\end{align}
We can rescale $\epsilon \leftarrow \epsilon / 2$ so that
\begin{align}
    (1-\epsilon) I_d \preceq (U^TS^TSU)^{-1} \preceq (1+\epsilon) I_d
\end{align}
and the effect of this rescaling will be hidden in the $O$ notation for the sketch size formulas.

We will define the events $E_k$, $k=1,\dots,q$ as $(1-\epsilon) I_d \preceq (U^TS_k^TS_kU)^{-1} \preceq (1+\epsilon) I_d$ and it follows that when the sketch sizes are selected according to \eqref{eq:sketchsizes_other_sketches}, we will have
\begin{align}
    P(E_k) \geq 1 - \delta \,, k=1,\dots,q \,.
\end{align}

Now, we find a simpler expression for the error of the single-sketch estimator:
\begin{align}
    A\hat{x} - Ax^* &= A(A^TS^TSA)^{-1}A^TS^TS^Tb - Ax^* \nonumber \\
    &= A(A^TS^TSA)^{-1}A^TS^TS^T(b^\perp +Ax^*) - Ax^* \nonumber \\
    &= A(A^TS^TSA)^{-1}A^TS^TS^Tb^\perp \nonumber \\
    % &= U\Sigma V^T V\Sigma^{-1}(U^TS^TSU)^{-1} \Sigma^{-1} V^TV\Sigma U^T S^TSb^\perp \nonumber \\
    &= U (U^TS^TSU)^{-1} U^T S^TSb^\perp \nonumber \\
    &= U Q z
\end{align}
where we defined $Q := (U^TS^TSU)^{-1}$, and $z := U^TS^TSb^\perp$.
The variance term is equal to
\begin{align}
    \Exs[\|A\hat{x} - Ax^*\|_2^2] &= \Exs[\|UQz\|_2^2] = \Exs[\|Qz\|_2^2] \,.
\end{align}
Conditioned on the event $E$, we can bound the expectation $\Exs[\|Qz\|_2^2 | E]$ as follows:
\begin{align}
    \Exs[\|Qz\|_2^2 | E] &= \Exs[z^TQ^TQz | E] \nonumber \\
    &\leq \Exs[(1+\epsilon)^2 \|z\|_2^2 | E] \nonumber \\
    &= (1+\epsilon)^2\Exs[ \|z\|_2^2 | E]
\end{align}
where we have used $Q^TQ \preceq (1+\epsilon)^2 I_d$, which follows from all eigenvalues of $Q$ being less than $(1+\epsilon)$ when conditioned on the event $E$. 
Next, from Lemma \ref{expected_obj_val_diff}, we obtain the following bound for the expected error
\begin{align}
    \Exs[f(\bar{x})-f(x^*) | E_1, \dots, E_q] &= \frac{1}{q} \Exs \left[ \|A\hat{x}-Ax^* \|_2^2 | E \right] + \frac{q-1}{q} \| \Exs[A\hat{x} - Ax^* | E] \|_2^2 \nonumber \\
    &\leq \frac{1}{q} (1+\epsilon^2) \Exs[ \|z\|_2^2 | E] + \frac{q-1}{q}4\epsilon \Exs[ \|z\|_2^2 | E] \,.
\end{align}
Using Markov's inequality gives us the following probability bound:
\begin{align}
    P(f(\bar{x}) - f(x^*) \geq \gamma | E_1,\dots,E_q) &\leq \frac{1}{\gamma} \Exs[f(\bar{x})-f(x^*) | E_1,\dots,E_q] \leq \frac{1}{\gamma q} \big((1+\epsilon)^2 + 4\epsilon(q-1)\big) \Exs[\|z\|_2^2 | E] \,.
\end{align}
The unconditioned probability can be computed as
\begin{align}
    P(f(\bar{x}) - f(x^*) \geq \gamma) &= P(f(\bar{x}) - f(x^*) \geq \gamma | E_1, \dots, E_q)P(\cap_{k=1}^q E_k) + P(f(\bar{x}) - f(x^*) \geq \gamma | \cup_{k=1}^q E_k^C)P(\cup_{k=1}^q E_k^C) \nonumber \\
    &\leq \frac{1}{\gamma q} \big((1+\epsilon)^2 + 4\epsilon(q-1)\big) \Exs[\|z\|_2^2 | E] + q\delta
\end{align}
where we have used $P(\cup_{k=1}^q E_k^C) \leq q\delta $ and $P(\cap_{k=1}^q E_k) \leq 1$.
We can obtain a bound for the relative error as follows:
\begin{align}
    P&\Big(\frac{f(\bar{x})}{f(x^*)} \geq 1 + \frac{\gamma}{f(x^*)} \Big) \leq \frac{1}{\gamma q} \big((1+\epsilon)^2 + 4\epsilon(q-1)\big) \Exs[\|z\|_2^2 | E] + q\delta \,.
\end{align}
Scaling $\gamma \leftarrow \gamma f(x^*)$ leads to
\begin{align}
    P&\Big(\frac{f(\bar{x})}{f(x^*)} \leq 1 + \gamma \Big) \geq 1 - q\delta - \frac{1}{\gamma f(x^*) q} \big((1+\epsilon)^2 + 4\epsilon(q-1)\big) \Exs[\|z\|_2^2 | E] \,.
\end{align}

We now plug the bounds for $\Exs[\|z\|_2^2 | E]$ from Lemma \ref{bound_z_norm_ROS}, \ref{bound_z_norm_uniform}, \ref{bound_z_norm_leverage} in the above bound and obtain:
\begin{itemize}
    \item \textbf{Randomized Hadamard sketch:} 
    \begin{align}
        &P\Big(\frac{f(\bar{x})}{f(x^*)} \leq 1 + \gamma \Big) \geq 1 - q\delta - \frac{d}{q\gamma m} \big((1+\epsilon)^2 + 4\epsilon(q-1)\big) \,.
    \end{align}
    
    \item \textbf{Uniform sampling with replacement:}
    \begin{align}
        &P\Big(\frac{f(\bar{x})}{f(x^*)} \leq 1 + \gamma \Big) \geq 1 - q\delta - \frac{\mu d}{q\gamma m} \big((1+\epsilon)^2 + 4\epsilon(q-1)\big) \,.
    \end{align}
    
    \item \textbf{Uniform sampling without replacement:}
    \begin{align}
        &P\Big(\frac{f(\bar{x})}{f(x^*)} \leq 1 + \gamma \Big) \geq 1 - q\delta - \frac{\mu d}{q\gamma m} \frac{n-m}{n-1} \big((1+\epsilon)^2 + 4\epsilon(q-1)\big) .
    \end{align}
    
    \item \textbf{Leverage score sampling:}
    \begin{align}
        &P\Big(\frac{f(\bar{x})}{f(x^*)} \leq 1 + \gamma \Big) \geq 1 - q\delta - \frac{d}{q\gamma m} \big((1+\epsilon)^2 + 4\epsilon(q-1)\big) \,.
    \end{align}
\end{itemize}
\end{IEEEproof}

\subsection{Proofs of Theorems and Lemmas in Section \ref{sec:nonlinear_problems}}
\begin{IEEEproof} [Proof of Theorem \ref{thm:dist_newton_sketch}]
The optimal update direction is given by
\begin{align*}
    \Delta_t^* =((H_t^{1/2})^T H_t^{1/2})^{-1}g_t = H_t^{-1} g_t
\end{align*}
and the estimate update direction due to a single sketch is given by
\begin{align*}
    \hat{\Delta}_{t,k} = \alpha_s ((H_t^{1/2})^TS_{t,k}^TS_{t,k}H_t^{1/2})^{-1}g_t.
\end{align*}
where $\alpha_s \in \mathbb{R}$ is the step size scaling factor to be determined.
Letting $S_{t,k}$ be a Gaussian sketch, the bias can be written as
\begin{align*}
    \Exs [H_t^{1/2}(\hat{\Delta}_{t,k} - \Delta_t^*)] &= \Exs[\alpha_s H_t^{1/2}((H_t^{1/2})^TS_{t,k}^TS_{t,k}H_t^{1/2})^{-1}g_t - H_t^{1/2}H_t^{-1}g_t] \\
    &= \alpha_s H_t^{1/2} \Exs[((H_t^{1/2})^TS_{t,k}^TS_{t,k}H_t^{1/2})^{-1}]g_t - H_t^{1/2}H_t^{-1}g_t \\
    &= \alpha_s \theta_1 H_t^{1/2}((H_t^{1/2})^TH_t^{1/2})^{-1}g_t - H_t^{1/2}H_t^{-1}g_t \\
    &= \left( \alpha_s \theta_1 - 1 \right) H_t^{1/2}H_t^{-1}g_t \,.
\end{align*}
In the third line, we plug in the mean of $((H_t^{1/2})^TS_{t,k}^TS_{t,k}H_t^{1/2})^{-1}$ which is distributed as inverse Wishart distribution (see Lemma \ref{theta1_theta2_lemma_TEMP}). This calculation shows that the single sketch estimator gives an unbiased update direction for $\alpha_s = 1/\theta_1$.

The variance analysis is as follows:
\begin{align*}
    \Exs [\|H_t^{1/2}(\hat{\Delta}_{t,k} - \Delta_t^*)\|_2^2] &=\Exs[\hat{\Delta}_{t,k}^TH_t\hat{\Delta}_{t,k} + {\Delta_t^*}^TH_t\Delta_t^* - 2{\Delta_t^*}^TH_t\hat{\Delta}_{t,k}] \\
    %&= \alpha_s^2 g_t^T \Exs[ ((H_t^{1/2})^TS_{t,k}^TS_{t,k}H_t^{1/2})^{-1} H_t ((H_t^{1/2})^TS_{t,k}^TS_{t,k}H_t^{1/2})^{-1} ]g_t + g_t^T H_t^{-1} g_t \\ 
    %&\quad - 2\alpha_s g_t^T \Exs[((H_t^{1/2})^TS_{t,k}^TS_{t,k}H_t^{1/2})^{-1}] g_t \\
    &= \alpha_s^2 g_t^T \Exs[ ((H_t^{1/2})^TS_{t,k}^TS_{t,k}H_t^{1/2})^{-1} H_t ((H_t^{1/2})^TS_{t,k}^T S_{t,k}H_t^{1/2})^{-1} ]g_t + \left(1 - 2 \alpha_s \theta_1 \right) g_t^T H_t^{-1} g_t .
\end{align*}
Plugging $H_t^{1/2}=U\Sigma V^T$ into the first term and assuming $H_t^{1/2}$ has full column rank, the expectation term becomes
\begin{align*}
    \Exs[ ((H_t^{1/2})^TS_{t,k}^TS_{t,k}H_t^{1/2})^{-1} H_t ((H_t^{1/2})^TS_{t,k}^TS_{t,k}H_t^{1/2})^{-1} ] &= V\Sigma^{-1} \Exs[ (U^TS_{t,k}^TS_{t,k}U)^{-2} ] \Sigma^{-1}V^T \\
    &= V\Sigma^{-1} (\theta_2 I_d) \Sigma^{-1}V^T \\
    &= \theta_2 V \Sigma^{-2} V^T ,
\end{align*}
% \begin{align*}
%     &g_t^T \Exs[ ((H_t^{1/2})^TS_{t,k}^TS_{t,k}H_t^{1/2})^{-1} H_t ((H_t^{1/2})^TS_{t,k}^TS_{t,k}H_t^{1/2})^{-1} ]g_t = \\
%     &=g_t^T V\Sigma^{-1} \Exs[ (U^TS_{t,k}^TS_{t,k}U)^{-2} ] \Sigma^{-1}V^T g_t \\
%     &= g_t^T V\Sigma^{-1} (\theta_2 I_d) \Sigma^{-1}V^T g_t \\
%     &= \theta_2 g_t^T V \Sigma^{-2} V^T g_t \,,
% \end{align*}
%
where the third line follows due to Lemma \ref{theta1_theta2_lemma_TEMP}. Because $H_t^{-1} = V\Sigma^{-2}V^T$, the variance becomes:
\begin{align*}
    \Exs[\|H_t^{1/2}(\hat{\Delta}_{t,k} - \Delta_t^*)\|_2^2] &= \left( \alpha_s^2 \theta_2 + 1 - 2 \alpha_s \theta_1\right) g_t^T V\Sigma^{-2}V^T g_t = \left( \alpha_s^2 \theta_2 + 1 - 2 \alpha_s \theta_1\right) \| \Sigma^{-1}V^Tg_t \|_2^2 \,.
\end{align*}
% If we assume the gradient $g_t$ is of the form $g_t = (H_t^{1/2})^Th_t$, then the variance can be further simplified into the following expression:
% \begin{align*}
%     \Exs[\|H_t^{1/2}(\hat{\Delta}_{t,k} - \Delta_t^*)\|_2^2] &= \left( \alpha_s^2 \theta_2 + 1 - 2\alpha_s \theta_1 \right) \|U^Th_t\|_2^2.
% \end{align*}
%
It follows that the variance is minimized when $\alpha_s$ is chosen as $\alpha_s = \theta_1 / \theta_2$.
\end{IEEEproof}

\begin{lemma} [\cite{lacotte2019faster_SM}] \label{theta1_theta2_lemma_TEMP}
For the Gaussian sketch matrix $S \in \mathbb{R}^{m\times n}$ with i.i.d. entries distributed as $\mathcal{N}(0,1/\sqrt{m})$ where $m \geq d$, and for $U\in \mathbb{R}^{n\times d}$ with $U^TU=I_d$, the following are true:
\begin{align}
    \Exs[(U^TS^TSU)^{-1}] &= \theta_1 I_d, \nonumber \\
    \Exs[(U^TS^TSU)^{-2}] &= \theta_2 I_d,
\end{align}
where $\theta_1$ and $\theta_2$ are defined as
\begin{align} \label{eq:theta_1_2_definitions_TEMP}
    \theta_1 &\coloneqq \frac{m}{m-d-1}, \nonumber \\
    \theta_2 &\coloneqq \frac{m^2(m-1)}{(m-d)(m-d-1)(m-d-3)}.
\end{align}
\end{lemma}

\begin{IEEEproof} [Proof of Theorem \ref{opt_lambda_2_newton}]
In the following, we omit the subscripts in $S_{t,k}$ for simplicity. Using the SVD decomposition of $H_t^{1/2} = U\Sigma V^T$, the bias can be written as
\begin{align*}
    \Exs[ H_t^{1/2}(\hat{\Delta}_{t,k} - \Delta_t^*)] &= U\Exs[ (U^TS^TSU + \lambda_2 \Sigma^{-2})^{-1}]\Sigma^{-1}V^Tg_t - U (I_d + \lambda_1 \Sigma^{-2})^{-1} \Sigma^{-1}V^Tg_t \,.
\end{align*}
By the assumption that $\Sigma = \sigma I_d$, the bias term can be simplified as
\begin{align*}
    \Exs[ H_t^{1/2} &(\hat{\Delta}_{t,k} - \Delta_t^*)] = \sigma^{-1} U\Exs[ (U^TS^TSU + \lambda_2 \sigma^{-2} I_d)^{-1}]V^Tg_t  - \sigma^{-1} (1+\lambda_1\sigma^{-2})^{-1} U V^Tg_t .
\end{align*}
By Lemma \ref{expectation_inverse_regularization_2}, as $n$ goes to infinity, we have
\begin{align*}
    \Exs[ H_t^{1/2}(\hat{\Delta}_{t,k} - \Delta_t^*)] &=\sigma^{-1} \left( \theta_3(d/m, \lambda_2\sigma^{-2}) - \frac{1}{1+\lambda_1\sigma^{-2}} \right) UV^Tg_t \\
    &= \frac{-\frac{\lambda_2}{\sigma^{2}}+\frac{d}{m}-1 + \sqrt{(-\frac{\lambda_2}{\sigma^{2}} + \frac{d}{m}-1)^2+4\frac{\lambda_2}{\sigma^{2}}\frac{d}{m}}}{2\lambda_2\sigma^{-1} d/m}UV^Tg_t - \frac{\sigma^{-1}}{1+\lambda_1\sigma^{-2}} UV^Tg_t  .
\end{align*}
The bias becomes zero for the value of $\lambda_2$ that satisfies the following equation:
\begin{align} \label{eq:lhs_defn}
    %\frac{-\lambda_2\sigma^{-2}+d/m-1 + \sqrt{(-\lambda_2\sigma^{-2} + d/m-1)^2+4\lambda_2\sigma^{-2}d/m}}{2\lambda_2\sigma^{-2} d/m} = \frac{1}{1+\lambda_1\sigma^{-2}} \nonumber \\
     &\sqrt{\left(-\sigma^{-2} + \frac{1}{\lambda_2}\left(\frac{d}{m}-1 \right) \right)^2 + 4\sigma^{-2}\frac{d}{m\lambda_2}} -\sigma^{-2} +\frac{1}{\lambda_2} \left(\frac{d}{m}-1\right) = 2\sigma^{-2}\frac{d}{m} \frac{1}{1+\lambda_1\sigma^{-2}} \,.
\end{align}
In the regime where $\lambda_2 \geq 0$, the LHS of \eqref{eq:lhs_defn} is always non-negative and is monotonically decreasing in $\lambda_2$.The LHS approaches zero as $\lambda_2 \rightarrow \infty$. We now consider the following cases:
\begin{itemize}
    \item Case 1: $m \leq d$. Because $d/m-1 \geq 0$, as $\lambda_2 \rightarrow 0$, the LHS goes to infinity. Since the LHS can take any values between $0$ and $\infty$, there is an appropriate $\lambda_2^*$ value that makes the bias zero for any $\lambda_1$.

    \item Case 2: $m > d$. In this case, $d/m-1<0$. The maximum of LHS in this case is reached as $\lambda_2 \rightarrow 0$ and it is equal to $2\sigma^{-2}\frac{d}{m-d}$. % used l'hospital's rule to obtain this
    As long as $2\sigma^{-2}\frac{d}{m} \frac{1}{1+\lambda_1\sigma^{-2}} \leq 2\sigma^{-2}\frac{d}{m-d}$ is true, then we can drive the bias down to zero. More simply, this corresponds to $\lambda_1\sigma^{-2} \geq -d/m$, which is always true. Therefore in the case of $m > d$ as well, there is a $\lambda_2^*$ value for any value of $\lambda_1$ that will drive the bias down to zero.
\end{itemize}
To sum up, for any given value for the regularization parameter $\lambda_1$, it is possible to find a $\lambda_2^*$ value to make the sketched update direction unbiased. The optimal value for $\lambda_2$ is given by $LHS^{-1}(2\sigma^{-2}\frac{d}{m} \frac{1}{1+\lambda_1\sigma^{-2}})$ where $LHS^{-1}(y)=\frac{4\sigma^{-2}y^{-1}d/m+2(d/m-1)}{y+2\sigma^{-2}}$, which simplifies to the following expression:
% \begin{align*}
%     \lambda_2^* = \frac{\lambda_1+\sigma^2\frac{d}{m}}{1 + \frac{d}{m} \frac{1}{1+(\lambda_1/\sigma^2)}}.
% \end{align*}
\begin{align*}
    \lambda_2^* = \left(\lambda_1 + \frac{d}{m}\sigma^2 \right) \left(1 - \frac{d/m}{1 + \lambda_1 \sigma^{-2} + d/m} \right).
\end{align*}
\end{IEEEproof}

\section{Additional Numerical Results}
In this section, we present additional experimental results.

\subsection{Scalability of the Serverless Implementation}
\begin{figure}%[!t]
\begin{center}
\centerline{\includegraphics[width=0.35\columnwidth]{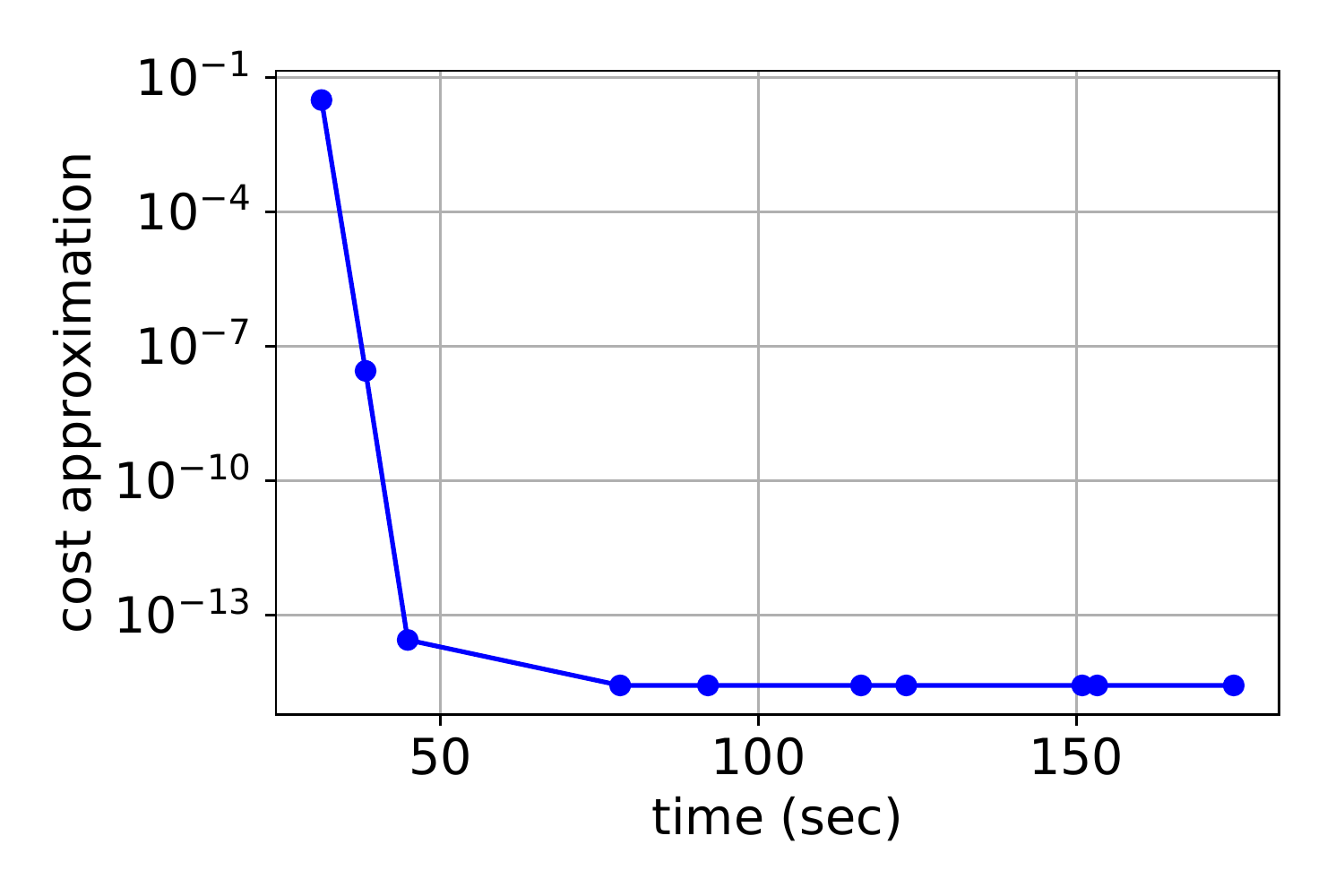}}
\caption{Cost approximation vs time when we solve the problem given in \eqref{log_barrier_opt_prob_2} for a large scale randomly generated dataset ($44$ GB sized) on AWS Lambda. Circles correspond to times that the iterates $x_t$ are computed.  Problem parameters are as follows: $n=200000$, $d=30000$, $\lambda_1=1$, $m=2000$, $q=100$, $\lambda=10$.}
\label{fig:log_barrier_awslambda_scalability}
\end{center}
\end{figure}

Figure \ref{fig:log_barrier_awslambda_scalability} shows the cost against time when we solve the problem given in \eqref{log_barrier_opt_prob_2} for large scale data on AWS Lambda using the distributed Newton sketch algorithm. The setting in this experiment is such that each worker node has access to a different subset of data, and there is no additional sketching applied. The dataset used is randomly generated and the goal here is to demonstrate the scalability of the algorithm and the serverless implementation. The size of the data matrix $A$ is $44$ GB.

In the serverless computing implementation, we reuse the serverless functions during the course of the algorithm, meaning that the same $q=100$ functions are used for every iteration. We note that every iteration requires two rounds of communication with the master node. The first round is for the communication of the local gradients, and the second round is for the approximate update directions. The master node, also a serverless function, is also reused across iterations. Figure \ref{fig:log_barrier_awslambda_scalability} illustrates that each iteration takes a different amount of time and iteration times can be as short as $5$ seconds. The reason for some iterations taking longer times is what is referred to as the straggler problem, which is a phenomenon commonly encountered in distributed computing. More precisely, the iteration time is determined by the slowest of the $q=100$ nodes and nodes often slow down for a variety of reasons causing stragglers. A possible solution to the issue of straggling nodes is to use error correcting codes to insert redundancy to computation and hence to avoid waiting for the outputs of all of the worker nodes \cite{Lee2018, bartan2019polar}. We identify that implementing straggler mitigation for solving large scale problems via approximate second order optimization methods such as distributed Newton sketch is a promising direction.

\subsection{Experiments on UCI Datasets}
In the case of large datasets and limited computing resources of worker nodes such as memory and lifetime, most of the standard sketches are computationally too expensive as discussed in the main body of the paper. This is the reason why we limited the scope of the large scale experiments to uniform sampling, SJLT, and hybrid sketch. In this section we present some additional experimental results on smaller datasets to empirically verify the theoretical results of the paper.

\begin{figure}%[!t]
\hspace{2cm}
\begin{minipage}[b]{0.3\linewidth}
  \centering
  \centerline{\includegraphics[width=\columnwidth]{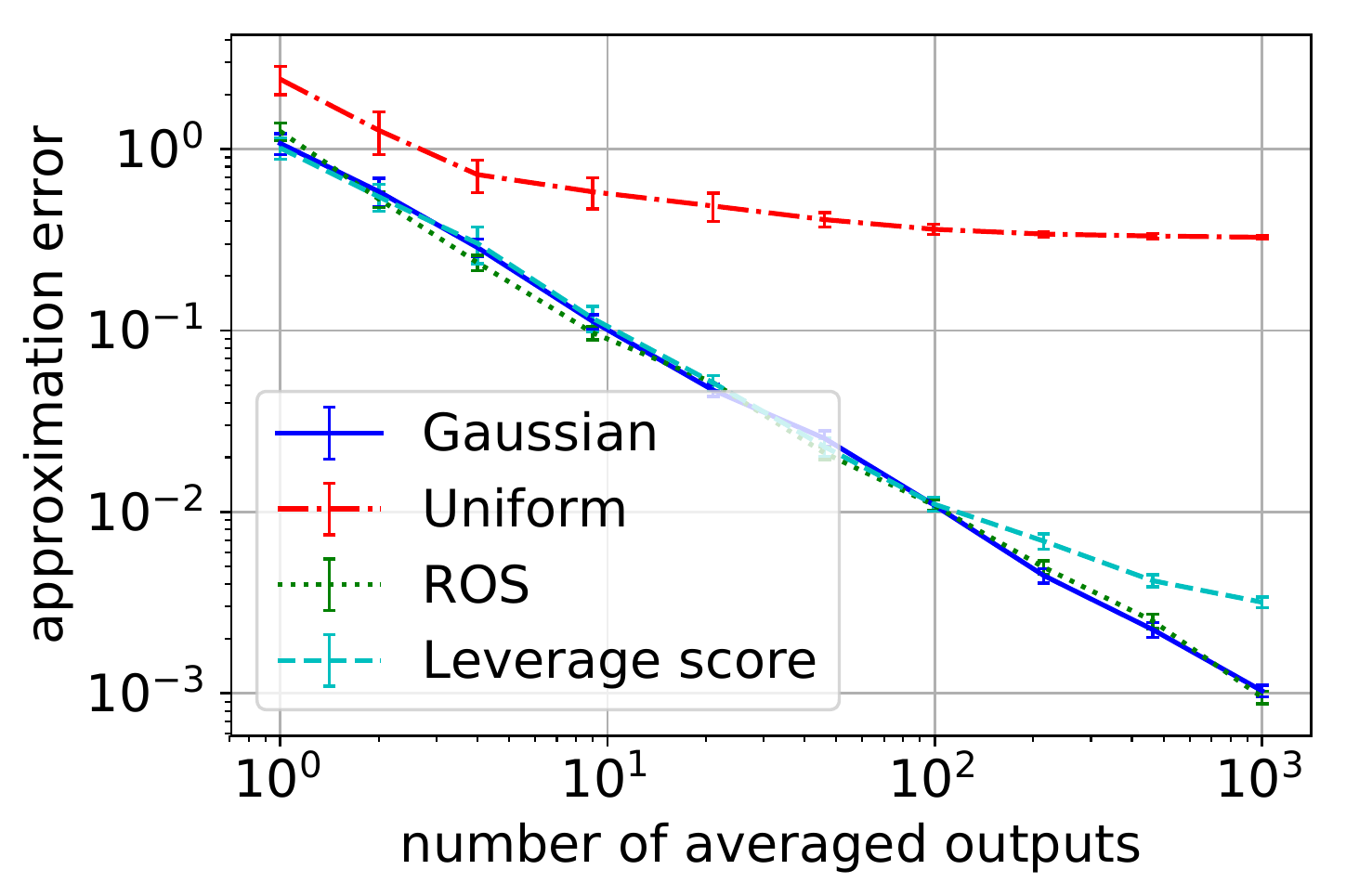}}
  \centerline{a) echocardiogram}\medskip
\end{minipage}
\hspace{2cm}
\begin{minipage}[b]{0.3\linewidth}
  \centering
  \centerline{\includegraphics[width=\columnwidth]{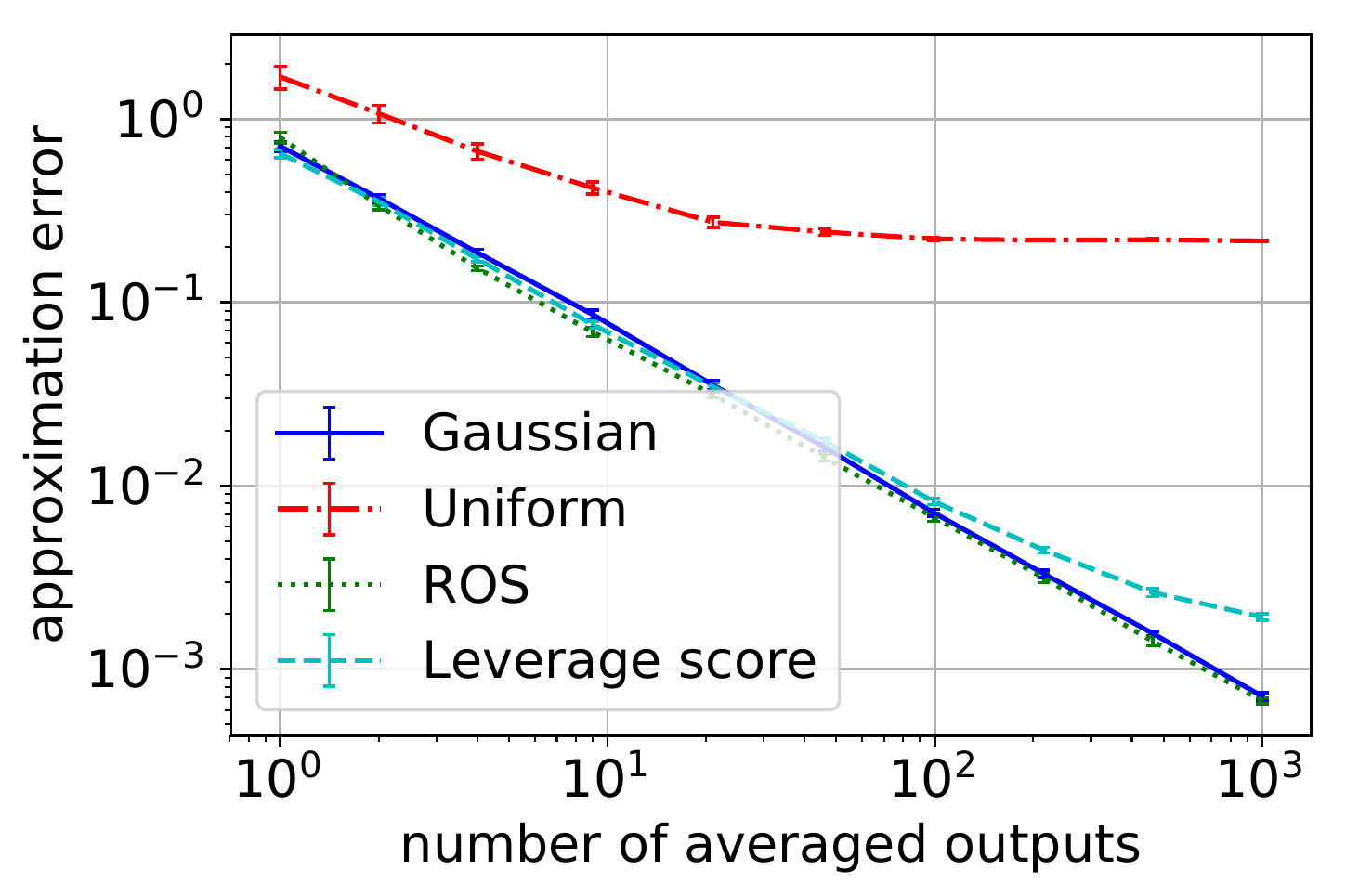}}
  \centerline{b) oocytes merluccius nucleus}\medskip
\end{minipage}
\caption{Approximation error against the number of averaged outputs in log-log scale for various sketching methods on two UCI datasets. All of the curves have been averaged over $25$ independent trials and the vertical error bars show the standard error. The parameters are as follows. Plot a: $n=131$, $d=10$, $m=20$. Plot b: $n=1022$, $d=41$, $m=100$.}
\label{supplement_uci_datasets}
\end{figure}

We present results on two UCI datasets in Figure \ref{supplement_uci_datasets} comparing the performances of the sketches we discussed in the paper.

Figure \ref{supplement_uci_datasets} shows that Gaussian and ROS sketches lead to unbiased estimators in the experiments because the corresponding curves appear linear in the log-log scale plots. These experiment results suggest that the upper bound that we have found for the bias of the ROS sketch may not be tight. We see that the estimates for uniform sampling and leverage score sampling are biased. The observation that the Gaussian sketch estimator is unbiased in the experiments is perfectly consistent with our theoretical findings. Furthermore, we observe that the approximation error is the highest in uniform sampling, which is also in agreement with the theoretical upper bounds that we have presented.

% \section*{Acknowledgment}

\ifCLASSOPTIONcaptionsoff
  \newpage
\fi



%\bibliographystyle{IEEEtran}
% \end{thebibliography}
\bibliographystyle{IEEEtran}
\bibliography{refs}

\end{document}